\documentclass[12pt]{article}
\usepackage{theorem}
\usepackage{amssymb}
\usepackage{graphicx}
\usepackage[mathscr]{eucal}
\usepackage{amsbsy}
\usepackage{amsmath}
\usepackage{stmaryrd}
\newtheorem{prop}{}[section]

{\theorembodyfont{\upshape} \newtheorem{rema}[prop]{}}
\newcommand{\boma}[1]{{\mbox{\boldmath $#1$} }}

\begin{document}
\newcommand{\uper}[1]{\stackrel{\barray{c} {~} \\ \mbox{\footnotesize{#1}}\farray}{\longrightarrow} }
\newcommand{\nop}[1]{ \|#1\|_{\piu} }
\newcommand{\no}[1]{ \|#1\| }
\newcommand{\nom}[1]{ \|#1\|_{\meno} }
\newcommand{\UU}[1]{e^{#1 \AA}}
\newcommand{\UD}[1]{e^{#1 \Delta}}
\newcommand{\bb}[1]{\mathbb{{#1}}}
\newcommand{\HO}[1]{\bb{H}^{{#1}}}
\newcommand{\Hz}[1]{\bb{H}^{{#1}}_{\zz}}
\newcommand{\Hs}[1]{\bb{H}^{{#1}}_{\ss}}
\newcommand{\Hg}[1]{\bb{H}^{{#1}}_{\gg}}
\newcommand{\HMs}[1]{\bb{H}^{{#1}}}
\newcommand{\CMs}[1]{\bb{C}^{{#1}}}
\newcommand{\EMs}[1]{\bb{E}^{{#1}}}
\newcommand{\EMnus}[1]{\bb{E}^{{#1}}_{\nu}}
\newcommand{\HM}[1]{\bb{H}^{{#1}}_{\so}}
\newcommand{\CM}[1]{\bb{C}^{{#1}}_{\so}}
\newcommand{\EM}[1]{\bb{E}^{{#1}}_{\so}}
\newcommand{\EMnu}[1]{\bb{E}^{{#1}}_{\so\,\nu}}
\newcommand{\vers}[1]{\widehat{#1}}
\def\Ax{D}
\def\ia{r}
\def\ja{s}
\def\ella{t}
\def\aa{a}
\def\Ss{\mathscr{S}}
\def\Nn{\mathscr{N}}
\def\Xa{\Theta}
\def\Xb{\Phi}
\def\Xc{\Psi}
\def\r{\ell}
\def\rr{\hat{\ell}}
\def\TA{A\,\,}
\def\TB{B\,\,}
\def\TC{C\,\,}
\def\TD{D\,\,}
\def\TE{E\,\,}
\def\TF{F\,\,}
\def\Gm{G^{[-]}}
\def\Xpd{X^{+}_d}
\def\Xmd{X^{-}_d}
\def\Ypd{Y^{+}_d}
\def\Ymd{Y^{-}_d}
\def\Yy{Y}
\def\II{\mathcal I}
\def\cc{\mathfrak{c}}
\def\A{\lambda}
\def\B{\mu}
\def\F{\nu}
\def\X{\boma{\omega}}
\def\Xx{\boma{\hat{\omega}}}
\def\S{\mathscr{H}}
\def\KP{K^{+}}
\def\KPP{K^{(+)}}
\def\GP{G^{+}}
\def\GPP{G^{(+)}}
\def\Kk{{\frak K}}
\def\Gk{{\frak G}}
\def\Km{{\mathscr K}}
\def\dK{\delta {\mathscr K}}
\def\dzK{\delta_0 {\mathscr K}}
\def\DK{\Delta {\mathscr K}}
\def\Cgot{{\mathfrak C}}
\def\DG{\Delta {\mathscr G}}
\def\Pis{P''}
\def\pp{p'}
\def\ps{p''}
\def\Ed{\hat{E}}
\def\Dd{\hat{D}}
\def\um{u_{-}}
\def\up{\mathfrak{u}}
\def\wp{w'}
\def\ws{w''}
\def\Wp{W'}
\def\Ws{W''}
\def\Pp{{\mathcal P}'}
\def\Ps{{\mathcal P}''}
\def\dG{\delta {\mathscr G}}
\def\Gg{{\mathscr G}}
\def\lam{\lambda}
\def\Lam{\Lambda}
\def\vh{\vers{h}}
\def\vk{\vers{k}}
\def\Do{\mathscr{E}}
\def\te{\vartheta}
\def\teta{H}
\def\we{\wedge}
\def\op{\,\mbox{\scriptsize{or}}\,}
\def\Sd{\Sfe^{d-1}}
\def\St{\Sfe^{2}}
\def\uz{u_{0}}
\def\vz{v_{0}}
\def\vi{v}
\def\ef{\psi}
\def\fun{\mathcal{F}}
\def\fun{{\tt f}}
\def\tvainf{\vspace{-0.4cm} \barray{ccc} \vspace{-0,1cm}{~}
\\ \vspace{-0.2cm} \longrightarrow \\ \vspace{-0.2cm} \scriptstyle{T \vain + \infty} \farray}
\def\De{F}
\def\er{\epsilon}
\def\erd{\er_0}
\def\Tn{T_{\star}}
\def\Tc{T_{\tt{c}}}
\def\Tb{T_{\tt{b}}}
\def\Tl{\mathscr{T}}
\def\Tm{T}
\def\pa{p_{\tt{a}}}
\def\Ta{T_{\tt{a}}}
\def\ua{u_{\tt{a}}}
\def\Tg{T_{G}}
\def\Tgg{T_{I}}
\def\Tw{T_{w}}
\def\Ts{T_{\Ss}}
\def\Tr{\Tl}
\def\Sp{\Ss'}
\def\Tsp{T_{\Sp}}
\def\vsm{\vspace{-0.1cm}\noindent}
\def\comple{\scriptscriptstyle{\complessi}}
\def\ug{u_G}
\def\nume{0.407}
\def\numerob{0.00724}
\def\deln{7/10}
\def\delnn{\dd{7 \over 10}}
\def\e{c}
\def\p{p}
\def\z{z}
\def\symd{{\mathfrak S}_d}
\def\del{\omega}
\def\Del{\delta}
\def\Di{\Delta}
\def\Ss{{\mathscr{S}}}
\def\Ww{{\mathscr{W}}}
\def\mmu{\hat{\mu}}
\def\rot{\mbox{rot}\,}
\def\curl{\mbox{curl}\,}
\def\Mm{\mathscr M}
\def\XS{\boma{x}}
\def\TS{\boma{t}}
\def\Lam{\boma{\eta}}
\def\DS{\boma{\rho}}
\def\KS{\boma{k}}
\def\LS{\boma{\lambda}}
\def\PR{\boma{p}}
\def\VS{\boma{v}}
\def\ski{\! \! \! \! \! \! \! \! \! \! \! \! \! \!}
\def\h{L}
\def\EMP{M'}
\def\R{R}
\def\Aa{{\mathscr{A}}}
\def\Rr{{\mathscr{R}}}
\def\Zz{{\mathscr{Z}}}
\def\Jj{{\mathscr{J}}}
\def\E{E}
\def\FFf{\mathscr{F}}
\def\Xim{\Xi_{\meno}}
\def\Ximn{\Xi_{n-1}}
\def\lan{\lambda}
\def\om{\omega}
\def\Om{\Omega}
\def\Sim{\Sigm}
\def\Sip{\Delta \Sigm}
\def\Sigm{{\mathscr{S}}}
\def\Ki{{\mathscr{K}}}
\def\Hi{{\mathscr{H}}}
\def\zz{{\scriptscriptstyle{0}}}
\def\ss{{\scriptscriptstyle{\Sigma}}}
\def\gg{{\scriptscriptstyle{\Gamma}}}
\def\so{\ss \zz}
\def\Dv{\bb{\DD}'}
\def\Dz{\bb{\DD}'_{\zz}}
\def\Ds{\bb{\DD}'_{\ss}}
\def\Dsz{\bb{\DD}'_{\so}}
\def\Dg{\bb{\DD}'_{\gg}}
\def\Ls{\bb{L}^2_{\ss}}
\def\Lg{\bb{L}^2_{\gg}}
\def\bF{{\bb{V}}}
\def\Fz{\bF_{\zz}}
\def\Fs{\bF_\ss}
\def\Fg{\bF_\gg}
\def\Pre{P}
\def\UUU{{\mathcal U}}
\def\fiapp{\phi}
\def\PU{P1}
\def\PD{P2}
\def\PT{P3}
\def\PQ{P4}
\def\PC{P5}
\def\PS{P6}
\def\Q{P6}
\def\Xp{Q3}
\def\Vi{V}
\def\bVi{\bb{V}}
\def\K{V}
\def\Ks{\bb{\K}_\ss}
\def\Kz{\bb{\K}_0}
\def\KM{\bb{\K}_{\, \so}}
\def\HGG{\bb{H}^\G}
\def\HG{\bb{H}^\G_{
\so}}
\def\EG{{\mathfrak{P}}^{\G}}
\def\G{G}
\def\de{\delta}
\def\esp{\sigma}
\def\dd{\displaystyle}
\def\LP{\mathfrak{L}}
\def\dive{\mbox{div}}
\def\la{\langle}
\def\ra{\rangle}
\def\um{u_{\meno}}
\def\uv{\mu_{\meno}}
\def\Fp{ {\textbf F_{\piu}} }
\def\Ff{ {\textbf F} }
\def\Fm{ {\textbf F_{\meno}} }
\def\Eb{ {\textbf E} }
\def\piu{\scriptscriptstyle{+}}
\def\meno{\scriptscriptstyle{-}}
\def\omeno{\scriptscriptstyle{\ominus}}
\def\Tt{ {\mathscr T} }
\def\Ee{ {\textbf E} }
\def\VP{{\mbox{\tt VP}}}
\def\CP{{\mbox{\tt CP}}}
\def\cp{$\CP(f_0, t_0)\,$}
\def\cop{$\CP(f_0)\,$}
\def\copn{$\CP_n(f_0)\,$}
\def\vp{$\VP(f_0, t_0)\,$}
\def\vop{$\VP(f_0)\,$}
\def\vopn{$\VP_n(f_0)\,$}
\def\vopdue{$\VP_2(f_0)\,$}
\def\leqs{\leqslant}
\def\geqs{\geqslant}
\def\mat{{\frak g}}
\def\tG{t_{\scriptscriptstyle{G}}}
\def\tN{t_{\scriptscriptstyle{N}}}
\def\TK{t_{\scriptscriptstyle{K}}}
\def\CK{C_{\scriptscriptstyle{K}}}
\def\CN{C_{\scriptscriptstyle{N}}}
\def\CG{C_{\scriptscriptstyle{G}}}
\def\CCG{{\mathscr{C}}_{\scriptscriptstyle{G}}}
\def\tf{{\tt f}}
\def\ti{{\tt t}}
\def\ta{{\tt a}}
\def\tc{{\tt c}}
\def\tF{{\tt R}}
\def\C{{\mathscr C}}
\def\P{{\mathscr P}}
\def\V{{\mathscr V}}
\def\TI{\tilde{I}}
\def\TJ{\tilde{J}}
\def\Lin{\mbox{Lin}}
\def\Hinfc{ H^{\infty}(\reali^d, \complessi) }
\def\Hnc{ H^{n}(\reali^d, \complessi) }
\def\Hmc{ H^{m}(\reali^d, \complessi) }
\def\Hac{ H^{a}(\reali^d, \complessi) }
\def\Dc{\DD(\reali^d, \complessi)}
\def\Dpc{\DD'(\reali^d, \complessi)}
\def\Sc{\SS(\reali^d, \complessi)}
\def\Spc{\SS'(\reali^d, \complessi)}
\def\Ldc{L^{2}(\reali^d, \complessi)}
\def\Lpc{L^{p}(\reali^d, \complessi)}
\def\Lqc{L^{q}(\reali^d, \complessi)}
\def\Lrc{L^{r}(\reali^d, \complessi)}
\def\Hinfr{ H^{\infty}(\reali^d, \reali) }
\def\Hnr{ H^{n}(\reali^d, \reali) }
\def\Hmr{ H^{m}(\reali^d, \reali) }
\def\Har{ H^{a}(\reali^d, \reali) }
\def\Dr{\DD(\reali^d, \reali)}
\def\Dpr{\DD'(\reali^d, \reali)}
\def\Sr{\SS(\reali^d, \reali)}
\def\Spr{\SS'(\reali^d, \reali)}
\def\Ldr{L^{2}(\reali^d, \reali)}
\def\Hinfk{ H^{\infty}(\reali^d, \KKK) }
\def\Hnk{ H^{n}(\reali^d, \KKK) }
\def\Hmk{ H^{m}(\reali^d, \KKK) }
\def\Hak{ H^{a}(\reali^d, \KKK) }
\def\Dk{\DD(\reali^d, \KKK)}
\def\Dpk{\DD'(\reali^d, \KKK)}
\def\Sk{\SS(\reali^d, \KKK)}
\def\Spk{\SS'(\reali^d, \KKK)}
\def\Ldk{L^{2}(\reali^d, \KKK)}
\def\Knb{K^{best}_n}
\def\sc{\cdot}
\def\k{\mbox{{\tt k}}}
\def\x{\mbox{{\tt x}}}
\def\g{ {\textbf g} }
\def\QQQ{ {\textbf Q} }
\def\AAA{ {\textbf A} }
\def\gr{\mbox{gr}}
\def\sgr{\mbox{sgr}}
\def\loc{\mbox{loc}}
\def\PZ{{\Lambda}}
\def\PZAL{\mbox{P}^{0}_\alpha}
\def\epsilona{\epsilon^{\scriptscriptstyle{<}}}
\def\epsilonb{\epsilon^{\scriptscriptstyle{>}}}
\def\lgraffa{ \mbox{\Large $\{$ } \hskip -0.2cm}
\def\rgraffa{ \mbox{\Large $\}$ } }
\def\restriction{\upharpoonright}
\def\M{\mu}
\def\m{m}
\def\Fre{Fr\'echet~}
\def\I{{\mathcal N}}
\def\ap{{\scriptscriptstyle{ap}}}
\def\fiap{\varphi_{\ap}}
\def\dfiap{{\dot \varphi}_{\ap}}
\def\DDD{ {\mathfrak D} }
\def\BBB{ {\textbf B} }
\def\EEE{ {\textbf E} }
\def\GGG{ {\textbf G} }
\def\TTT{ {\textbf T} }
\def\KKK{ {\textbf K} }
\def\HHH{ {\textbf K} }
\def\FFi{ {\bf \Phi} }
\def\GGam{ {\bf \Gamma} }
\def\sc{ {\scriptstyle{\bullet} }}
\def\a{a}
\def\ep{\epsilon}
\def\c{\kappa}
\def\parn{\par \noindent}
\def\elle{L}
\def\ro{\rho}
\def\al{\alpha}
\def\si{\sigma}
\def\be{\beta}
\def\ga{\gamma}
\def\ch{\chi}
\def\et{\eta}
\def\complessi{{\bf C}}
\def\len{{\bf L}}
\def\reali{{\bf R}}
\def\interi{{\bf Z}}
\def\Z{{\bf Z}}
\def\naturali{{\bf N}}
\def\Sfe{ {\bf S} }
\def\To{ {\bf T} }
\def\Td{ {\To}^d }
\def\Tt{ {\To}^3 }
\def\Zd{ \interi^d }
\def\Id{{\bf{I}}^d }
\def\Zt{ \interi^3 }
\def\Zet{{\mathscr{Z}}}
\def\Ze{\Zet^d}
\def\T1{{\textbf To}^{1}}
\def\es{s}
\def\ee{{E}}
\def\FF{\mathcal F}
\def\FFu{ {\textbf F_{1}} }
\def\FFd{ {\textbf F_{2}} }
\def\GG{{\mathcal G} }
\def\EE{{\mathcal E}}
\def\KK{{\mathcal K}}
\def\PP{{\mathcal P}}
\def\PPP{{\mathscr P}}
\def\PN{{\mathcal P}}
\def\PPN{{\mathscr P}}
\def\QQ{{\mathcal Q}}
\def\J{J}
\def\Np{{\hat{N}}}
\def\Lp{{\hat{L}}}
\def\Jp{{\hat{J}}}
\def\Pp{{\hat{P}}}
\def\Pip{{\hat{\Pi}}}
\def\Vp{{\hat{V}}}
\def\Ep{{\hat{E}}}
\def\Gp{{\hat{G}}}
\def\Kp{{\hat{K}}}
\def\Ip{{\hat{I}}}
\def\Tp{{\hat{T}}}
\def\Mp{{\hat{M}}}
\def\La{\Lambda}
\def\Ga{\Gamma}
\def\Si{\Sigma}
\def\Upsi{\Upsilon}
\def\Gam{\Gamma}
\def\Gag{{\check{\Gamma}}}
\def\Lap{{\hat{\Lambda}}}
\def\Upsig{{\check{\Upsilon}}}
\def\Kg{{\check{K}}}
\def\ellp{{\hat{\ell}}}
\def\j{j}
\def\jp{{\hat{j}}}
\def\BB{{\mathcal B}}
\def\LL{{\mathcal L}}
\def\MM{{\mathcal U}}
\def\SS{{\mathcal S}}
\def\DD{D}
\def\Dd{{\mathcal D}}
\def\VV{{\mathcal V}}
\def\WW{{\mathcal W}}
\def\OO{{\mathcal O}}
\def\RR{{\mathcal R}}
\def\TT{{\mathcal T}}
\def\AA{{\mathcal A}}
\def\CC{{\mathcal C}}
\def\JJ{{\mathcal J}}
\def\NN{{\mathcal N}}
\def\HH{{\mathcal H}}
\def\XX{{\mathcal X}}
\def\XXX{{\mathscr X}}
\def\YY{{\mathcal Y}}
\def\ZZ{{\mathcal Z}}
\def\CC{{\mathcal C}}
\def\cir{{\scriptscriptstyle \circ}}
\def\circa{\thickapprox}
\def\vain{\rightarrow}
\def\salto{\vskip 0.2truecm \noindent}
\def\spazio{\vskip 0.5truecm \noindent}
\def\vs1{\vskip 1cm \noindent}
\def\fine{\hfill $\square$ \vskip 0.2cm \noindent}
\def\ffine{\hfill $\lozenge$ \vskip 0.2cm \noindent}
\newcommand{\rref}[1]{(\ref{#1})}
\def\beq{\begin{equation}}
\def\feq{\end{equation}}
\def\beqq{\begin{eqnarray}}
\def\feqq{\end{eqnarray}}
\def\barray{\begin{array}}
\def\farray{\end{array}}
\makeatletter \@addtoreset{equation}{section}
\renewcommand{\theequation}{\thesection.\arabic{equation}}
\makeatother
\begin{titlepage}
{~}
\vspace{1cm}
\begin{center}
{\Large \textbf{New results on the constants in some inequalities}}
\vskip 0.2cm
{~}
\hskip -0.4cm
{\Large\textbf{for the Navier-Stokes quadratic nonlinearity}}
\end{center}
\vspace{0.5truecm}
\begin{center}
{\large
Carlo Morosi$\,{}^a$, Mario Pernici$\,{}^b$, Livio Pizzocchero$\,{}^c$({\footnote{Corresponding author}})} \\
\vspace{0.5truecm} ${}^a$ Dipartimento di Matematica, Politecnico di Milano,
\\ P.za L. da Vinci 32, I-20133 Milano, Italy \\
e--mail: carlo.morosi@polimi.it \\
${}^b$ Istituto Nazionale di Fisica Nucleare, Sezione di Milano, \\
Via Celoria 16, I-20133 Milano, Italy \\
e--mail: mario.pernici@mi.infn.it
\\
${}^c$ Dipartimento di Matematica, Universit\`a di Milano\\
Via C. Saldini 50, I-20133 Milano, Italy\\
and Istituto Nazionale di Fisica Nucleare, Sezione di Milano, Italy \\
e--mail: livio.pizzocchero@unimi.it
\end{center}
\begin{abstract}
We give fully explicit upper and lower bounds for the constants in two
known inequalities related to the quadratic nonlinearity of the
incompressible (Euler or) Navier-Stokes
equations on the torus $\Td$. These inequalities
are ``tame'' generalizations (in the sense of Nash-Moser)
of the ones analyzed in the previous works [Morosi and Pizzocchero:
CPAA 2012, Appl.Math.Lett. 2013].
\end{abstract}
\vspace{1cm} \noindent
\hfill \parn
\par \vspace{0.05truecm} \noindent \textbf{AMS 2000 Subject classifications:} 76D05, 26D10, 46E35.
\end{titlepage}
\section{Introduction}
Let us consider the homogeneous incompressible Navier-Stokes (NS) equations on a torus
$\Td = (\reali/2 \pi \interi)^d$ of arbitrary dimension; the nonlinear part of these
equations is governed
by the bilinear map $\PPP$ sending two sufficiently regular
vector fields $v, w : \Td \vain \reali^d$ into
\beq \PPP(v,w) := \LP(v \sc \partial w)~. \label{mapp} \feq
In the above $v \sc \partial w: \Td \vain \reali^d$ is
the vector field of components $(v \sc \partial w)_s := \sum_{r=1}^d v_r \partial_r w_s$
and $\LP$ is the Leray projection onto the space of divergence free vector fields
(see Section \ref{notations} for more details). Of course the NS equations read
\beq {\partial u \over \partial t}  = \nu \Delta u - \PPP(u,u) + f ,
\label{eulnu} \feq
where: $u = u(x,t)$ is the divergence free velocity field, depending on $x \in \Td$ and
on time $t$; $\nu \geqs 0$ is the kinematic viscosity, $\Delta$ is the Laplacian
of $\Td$; $f = f(x,t)$ is the (Leray projected) external force per unit mass.
In the inviscid case $\nu=0$, \rref{eulnu} become the Euler equations.\par
In this paper we focus the attention on certain inequalities fulfilled by $\PPP$ in the
framework of Sobolev spaces. For any real $n$, we denote with $\HM{n}$ the Sobolev space formed
by the (distributional) vector fields $v$ on $\Td$ with vanishing divergence
and mean, such that $\sqrt{-\Delta}^{\,n} v$ is in $L^2$; this carries
the inner product $\la v | w \ra_n := \la \sqrt{-\Delta}^{\,n} v | \sqrt{-\Delta}^{\,n} w \ra_{L^2}$
and the norm $\| v \|_n := \sqrt{\la v | v \ra_n}$ (see the forthcoming Eqs.\,\rref{defhn} \rref{definner}).
Let $p, n$ be real numbers; it is known that
\beq n > d/2,~ v \in \HM{n},~ w \in \HM{n+1} \quad \Rightarrow \quad \PPP(v,w) \in \HM{n} \label{known} \feq
and that there are positive real constants $K_n$, $G_n$, $K_{p n}$, $G_{p n}$ such that:
\beq \| \PPP(v, w) \|_n \leqs K_n \| v \|_n \| w \|_{n+1} \qquad \mbox{for $n > d/2$,
$v \in \HM{n}$, $w \in \HM{n+1}$}~, \label{basineq} \feq
\beq | \la \PPP(v, w) | w \ra_n | \leqs G_n \| v \|_n \| w \|^2_n
\qquad \mbox{for $n > d/2 + 1$, $v \in \HM{n}$,
$w \in \HM{n+1}$}~, \label{katineq} \feq
\beq \| \PPP(v, w) \|_p \leqs {1 \over 2} K_{p n} ( \| v \|_p \| w \|_{n+1} + \| v \|_n \| w \|_{p+1})
\label{basineqa} \feq
$$ \qquad \mbox{for $p \geqs n > d/2$, $v \in \HM{p}$, $w \in \HM{p+1}$}~, $$
\beq | \la \PPP(v, w) | w \ra_p | \leqs
{1 \over 2} G_{p n} (\| v \|_p \| w \|_n + \| v \|_n \| w \|_p)\| w \|_p \label{katineqa} \feq
$$ \mbox{for $p \geqs n > d/2 + 1$, $v \in \HM{p}$, $w \in \HM{p+1}$}~. $$
Statements \rref{known} \rref{basineq} indicate that $\PPP$ maps continuously
$\HM{n} \times \HM{n+1}$ to $\HM{n}$ if $n > d/2$. Eq.\,\rref{basineqa} with $p=n$ implies Eq.\,\rref{basineq}, with $K_n := K_{n n}$;
similarly, \rref{katineqa} with $p=n$ gives \rref{katineq} with $G_n := G_{n n}$.
\par
Eq.\,\rref{basineq}
is closely related to the basic norm inequalities about
multiplication in Sobolev spaces, and \rref{katineq} is due to Kato \cite{Kato};
for these reasons, in \cite{cog} \cite{cok} the inequalities
\rref{basineq} and \rref{katineq}
are referred to, respectively, as the ``basic'' and
``Kato'' inequalities for $\PPP$
({\footnote{Due to a remark of \cite{cog}, we could write the inequality
\rref{katineq} and its extension \rref{katineqa}
using, in place of $\PPP(v,w) =
\LP(v \sc w)$, the vector field
(with non zero divergence) $v \sc \partial w$.
The cited reference considers the Sobolev space
$\Hz{n}$ of vector fields $v$ on $\Td$ with vanishing
mean and $\sqrt{-\Delta}^{\,n} v$ in $L^2$,
with the inner product
$\la v | w \ra_n := \la \sqrt{-\Delta}^{\,n} v | \sqrt{-\Delta}^{\,n} \ra_{L^2}$;
for any $n > d/2$ and $v \in \HM{n}$, $w \in \HM{n+1}$ one has
$v \sc \partial w \in \Hz{n}, \PPP(v,w) \in \HM{n}$ and
$\la \PPP(v,w) | w \ra_n = \la v \sc \partial w | w \ra_n$. However,
these considerations will play no role in the present paper.}}).
Eqs.\,\rref{basineqa}\,\rref{katineqa}
are tame refinements of \rref{basineq} \rref{katineq}
(in the general sense given to tameness
in studies on the Nash-Moser implicit function theorem \cite{Ham}).
We remark that inequalities
very similar to \rref{katineqa} are used by
Temam in \cite{Tem}, Beale-Kato-Majda in \cite{BKM}
and Robinson-Sadowski-Silva in the recent work \cite{RSS}.
\par
From here to the end of the paper we intend
\beq K_{n}, G_n, K_{p n}, G_{p n} := \mbox{the sharp constants in
\rref{basineq} \rref{katineq} \rref{basineqa} \rref{katineqa}} \feq
(i.e., the minimum constants fulfilling these inequalities). In the previous
papers \cite{cog} \cite{cok}, explicit upper and lower bounds
were provided for $K_n$ and $G_n$. In the present work we generalize
the cited results deriving upper and lower bounds for $K_{p n}$ and
$G_{p n}$, for all real $p,n$ as in Eqs.\,\rref{basineqa}\,\rref{katineqa}.
Our derivations of the upper bounds also give, as byproducts,
simple and self-consistent proofs of the related inequalities;
the approach proposed follows ideas from Temam \cite{Tem} and Constantin-Foias
\cite{CoFo}, making them more quantitative. The lower bounds are
obtained substituting suitable trial vector fields in Eqs.\,\rref{basineqa}\,\rref{katineqa}.
\par
The relevance of a quantitative information on
the constants $K_{p n}, G_{p n}$ is
pointed out, e.g., in \cite{smo}.
In the cited work, the inequalities
\rref{basineq}\,-\,\rref{katineqa} and the constants therein
are used to give bounds on the exact $C^\infty$ solution of the NS Cauchy problem
with smooth initial data (including the Euler case $\nu=0$)
via the \emph{a posteriori} analysis of an approximate solution;
these estimates concern the interval of existence of the
exact solution and its Sobolev distance of any order from
the approximate solution. Paper \cite{smo} uses systematically the known fact
that the space of $C^\infty$ vector fields on $\Td$ with vanishing divergence
and mean coincides with $\cap_{p \in \reali} \HM{p}$;
the tame structure of the inequality
\rref{katineqa} is essential for an efficient implementation
of the a posteriori analysis since, after fixing a basic
order $n > d/2 + 1$, it induces simple estimates in terms
of the Sobolev norms of arbitrary order $p \geqs n$.
The setting of \cite{smo} is in fact a $C^\infty$ variant
of the framework introduced in \cite{appeul} (and
inspired by Chernyshenko \emph{et al.} \cite{Che}), where the exact
and approximate NS solutions live in a Sobolev space of
a given finite order, and the a posteriori analysis
is based only on the inequalities \rref{basineq} \rref{katineq}.
For some applications of the general schemes of \cite{appeul} \cite{smo},
in addition to these papers we wish to mention \cite{padova} \cite{reylarge}
\cite{apprey}.
\vskip 0.2cm \noindent
\textbf{Organization and main results of the paper.}
Section \ref{notations} reviews some basic notations and presents a number
of elementary facts about the bilinear map $\PPP$;
one of these facts is proved in Appendix \ref{appehel}.
The subsequent Sections \ref{secbasic} and \ref{seckato}
present our upper bounds $K^{+}_{p n}, G^{+}_{p n}$ for
the sharp constants \rref{basineqa} and \rref{katineqa}, respectively;
these are described by Theorems \ref{maink}, \ref{maing}
and have the form
\beq K^{+}_{p n} = {1 \over (2 \pi)^{d/2}} \sqrt{ \sup_{k \in \Zd \setminus \{0 \}}
\KK_{p n}(k)}~, \qquad G^{+}_{p n} = {1 \over (2 \pi)^{d/2}} \sqrt{ \sup_{k \in \Zd \setminus \{0 \}}
\GG_{p n}(k)} \label{ineq} \feq
where $\KK_{p n}, \GG_{p n}: \Zd \setminus \{0 \} \vain [0,+\infty)$ are
explicitly given, bounded functions. For each $k$, $\KK_{p n}(k)$
and $\GG_{p n}(k)$ are infinite (zeta-type) sums
over the lattice $\Zd$ or, to be precise, on $\Zd \setminus \{0 ,k\}$:
see Eqs.\,\rref{kknd}\,\rref{cantake}\,\rref{ggnd}\,\rref{cantakeg}.
Sections \ref{secbasic} and \ref{seckato} also propose, as preliminary
results, some elementary upper bounds on the sups in Eq.\,\rref{ineq};
these imply elementary upper bounds $K^{\la + \ra}_{p n}, G^{\la + \ra}_{p n}$
for $K_{p n}$ and $G_{p n}$, much rougher than $K^{+}_{p n}$ and $G^{+}_{p n}$.
\par
The next step along these lines
is the accurate computation of the functions $\KK_{p n},
\GG_{p n}$ and of their sups. After some preliminaries
presented in Section \ref{proveg}, this subject
is discussed in detail in Sections \ref{sectk} and \ref{sectg};
the basic idea is to approximate
the infinite sums $\KK_{p n}(k), \GG_{p n}(k)$
with finite sums over the integer points of
suitable balls, giving accurate reminder estimates;
in the same spirit, the sups of $\KK_{p n}$
and $\GG_{p n}$ over $\Zd$ are approximated
with sups over the integer points of a ball,
giving again error estimates. This construction
finally produces precise upper approximants
$\KPP_{p n}$, $\GPP_{p n}$ for $\KP_{p n}, \GP_{p n}$.
\par
The procedures of Sections \ref{sectk} and \ref{sectg} are suitable for automatic
computations. Indeed, we have implemented
such procedures writing a {\tt{C}} program for
the computation of the previously mentioned
approximants via finite sums,
and using {\tt{Mathematica}} for some related
symbolic and numerical calculations
({\footnote{Throughout the paper,
an expression like $r= a. b c d...~$ means the following: computation
of real number $r$ via {\tt{C}} or {\tt{Mathematica}}
produces as an output $a.b c d$, followed by other digits not reported for brevity.
As indicated in Section \ref{sec345},
some of the {\tt{C}} computations
have been validated using the {\tt{Arb}} library
\cite{arb}, that gives certified roundoff errors.}}).
In Section \ref{sec345} we
give some details on the overall
procedures, on their computer
implementation and, in particular,
on the calculations of the following bounds: $\KPP_{p n}$ for
\beq d = 3~, \qquad n=2~, \qquad p = 2,3,...,10~, \label{casesdue} \feq
\beq d = 3~, \qquad n=3~, \qquad p = 3,4,...,10~~ \label{casestre} \feq
and $\GPP_{p n}$ for the cases \rref{casestre}.
\par
In Sections \ref{seclowk} and \ref{seclowg} we
derive some lower bounds
$K^{(-)}_{p n}$, $G^{(-)}_{p n}$ for $K_{p n}, G_{p n}$, respectively;
as anticipated, these are obtained substituting
for $v$ and $w$ in Eqs.\,\rref{basineqa}\,\rref{katineqa}
suitable trial vector fields (which are relatively simple,
since they have finitely many nonzero
Fourier components). The chosen vector fields often depend
on one or more parameters, so the best lower bounds arising from them are obtained
by optimization with respect to the parameters. Both
Section \ref{seclowk} and \ref{seclowg}
exhibit the numerical values of the
above mentioned lower bounds or,
to be precise, of some lower approximants for them, in the cases
\rref{casesdue} or \rref{casestre} already
considered in connection with the upper bounds.
For the reader's convenience, hereafter we summarize
in Tables \TA\!\!,\,\TB the numerical values of the bounds
$K^{(\pm)}_{p n}$, $G^{(\pm)}_{p n}$ computed in Sections
\ref{sec345}, \ref{seclowk}, \ref{seclowg}, together with
the ratios of the lower to the upper bounds. \par
Section \ref{seclowk} and \ref{seclowg} also give rougher lower bounds
$K^{\la - \ra}_{p n}, G^{\la - \ra}_{p n}$. These can be combined
with the rough upper bounds $K^{\la + \ra}_{p n}, G^{\la + \ra}_{p n}$
of Sections \ref{sectk} and \ref{sectg} to prove the following statement on
the sharp constants $K_{p n} , G_{p n}$: for fixed ($d$ and) $n$,
\beq (K_{p n})^{1/p}, (G_{p n})^{1/p} \vain 2 \qquad \mbox{for $p \vain + \infty$}~; \feq
this concluding result is the subject of Section \ref{secinf}.
\vskip 0.5cm \noindent
\vbox{
\noindent
\textbf{Table \TA. Upper and lower bounds $\boma{K^{(\pm)}_{p n}}$ on the
constants $\boma{K_{p n}}$, with their ratios, in the cases \rref{casesdue}
\rref{casestre}}
\vskip 0.5cm \noindent
$$ \begin{tabular}{|c|| c|c|c|}
\hline
$(p,n)$ & $K^{(-)}_{p n}$ & $\KPP_{p n}$ & $K^{(-)}_{p n}/\KPP_{p n}$ \\[0.1cm]
\hline \hline
$(2,2)$ & 0.126 & 0.335 & 0.376...  \\ \hline
$(3,2)$ & 0.179 & 0.643 & 0.278...  \\ \hline
$(4,2)$ & 0.264 & 0.831 & 0.317... \\ \hline
$(5,2)$ & 0.463 & 1.16 & 0.339...  \\ \hline
$(6,2)$ & 0.793 & 1.94 & 0.408...  \\ \hline
$(7,2)$ & 1.33 & 3.02 & 0.440...  \\ \hline
$(8,2)$ & 2.20 & 5.07 & 0.433...  \\ \hline
$(9,2)$ & 3.60 & 8.54 & 0.421...  \\ \hline
$(10,2)$ & 5.83 & 14.5 & 0.402...  \\ \hline
$(3,3)$ & 0.179 & 0.320 & 0.559...  \\ \hline
$(4,3)$ & 0.253 & 0.539 & 0.469... \\ \hline
$(5,3)$ & 0.418 & 0.909 & 0.459...  \\ \hline
$(6,3)$ & 0.732 & 1.54 & 0.475...  \\ \hline
$(7,3)$ & 1.25 & 2.58 & 0.484...  \\ \hline
$(8,3)$ & 2.10 & 4.28 & 0.490...  \\ \hline
$(9,3)$ & 3.48 & 7.04 & 0.494...  \\ \hline
$(10,3)$ & 5.69  & 11.5 & 0.494...  \\ \hline
\end{tabular} $$
}
\vfill \eject \noindent
{~}
\vskip -2cm \noindent
\vbox{
\noindent
\textbf{Table \TB. Upper and lower bounds $\boma{G^{(\pm)}_{p n}}$ on the
constants $\boma{G_{p n}}$, with their ratios, in the cases \rref{casestre}}
\vskip 0.5cm \noindent
$$ \begin{tabular}{|c|| c|c|c|}
\hline
$(p,n)$ & $G^{(-)}_{p n}$ & $\GPP_{p n}$ & $G^{(-)}_{p n}/\GPP_{p n}$ \\[0.1cm]
\hline \hline
$(3,3)$ & 0.121 & 0.438 & 0.276...  \\ \hline
$(4,3)$ & 0.235 & 1.03 & 0.228... \\ \hline
$(5,3)$ & 0.408 & 1.26 & 0.323...  \\ \hline
$(6,3)$ & 0.674 & 2.06 & 0.327...  \\ \hline
$(7,3)$ & 1.08 & 3.58 & 0.301...  \\ \hline
$(8,3)$ & 1.74 & 5.68 & 0.306...  \\ \hline
$(9,3)$ & 2.77 & 9.64 & 0.287...  \\ \hline
$(10,3)$ & 4.40 & 16.4 & 0.268...  \\ \hline
\end{tabular} $$
}
\section{Preliminaries}
\label{notations}
Throughout the paper we work
in any dimension $d \in \{2,3,...\}$.
\salto
\textbf{Some notations.} For $a, b \in \complessi^d$ we write
$a \sc b := \sum_{r=1}^d a_r b_r$ , $\overline{a} := (\overline{{a_r}})_{r=1,...,d}$
and $|a| := \sqrt{\overline{a} \sc a}$.
We often consider the torus $\Td := (\reali/2 \pi \interi)^d$ and the
lattice $\Zd$, associated to it in Fourier analysis. For $\ell, k$ in $\Zd$ we put
\beq \Zd_\ell := \Zd \setminus \{\ell \}~; \qquad
\Zd_{\ell k} := \Zd \setminus \{\ell,k \}~. \feq
\salto
\textbf{Function spaces.}
When working on $\Td$, we often use the Fourier basis $(e_k)_{k \in \Zd}$, where
\beq e_k(x) := {e^{i k \sc x} \over (2 \pi)^{d/2}} \qquad \mbox{for $x \in \Td$}~. \label{furba} \feq
The space of $\reali^d$-valued distributions on $\Td$ is denoted with
\beq D'(\Td, \reali^d) \equiv \mathbb{D}'(\Td) \equiv
\mathbb{D}'~; \feq
each $v \in \mathbb{D}'$ has a weakly convergent
Fourier expansion $v = \sum_{k \in \Zd} v_k e_k$,
with Fourier coefficients $v_k = \overline{v_{-k}} \in \complessi^d$. \par
We have a divergence operator $\dive : \mathbb{D}' \vain D'$, $v \mapsto \dive v$
where $D' \equiv D'(\Td, \reali)$ is the space of real distributions;
this has the Fourier representation, of obvious meaning,
$(\dive v)_k = k \sc v_k$.
For $v \in \mathbb{D}'$ the mean value $\la v \ra \in \reali^d$ is, by definition, the action of $v$
on the constant test function $(2 \pi)^{-d}$, and $\la v \ra = (2 \pi)^{-d/2} v_0$.
If $\mathbb{X}$ is the space $\mathbb{D}'$
or any vector subspace of it we put
\beq \mathbb{X}_\ss := \{ v \in \mathbb{X}~|~\dive v = 0 \} =
\{ v \in \mathbb{X}~|~k \sc v_k = 0~~\mbox{for $k \in \Zd$} \}~, \label{xss} \feq
\beq \mathbb{X}_\zz := \{ v \in \mathbb{X}~|~\la v \ra = 0 \}
= \{ v \in \mathbb{X}~|~v_0 = 0 \}~, \feq
\beq \mathbb{X}_{\so} := \mathbb{X}_{\ss} \cap \mathbb{X}_{\zz}~. \label{xso} \feq
The Laplacian is an operator $\Delta : \mathbb{D}' \vain \mathbb{D}'_\zz$,
$v \mapsto \Delta v$ with the Fourier representation $(\Delta v)_k = -|k|^2 v_k$.
If $v \in \mathbb{D}'_\zz$ and $n \in \reali$, we define $\sqrt{-\Delta}^n v$
to be the element of $\mathbb{D}'_\zz$ with Fourier coefficients $(\sqrt{-\Delta}^n v)_k =
|k|^n v_k$ for all $k \in \Zd_0$. \par
In the sequel we consider the spaces
\beq L^p(\Td, \reali^d) \equiv \mathbb{L}^p(\Td) \equiv \mathbb{L}^p~, \feq
most frequently in the Hilbertian case $p=2$.
The notations $\mathbb{L}^p_\ss$, $\mathbb{L}^p_\zz$,
$\mathbb{L}^p_{\so}$ are intended according to Eqs.\,\rref{xss}-\rref{xso}.
For any $n \in \reali$, we introduce the Sobolev space
\beq {~} \hspace{-0.5cm} \HM{n}(\Td) \equiv \HM{n} :=
\{ v \in \mathbb{D}'~|~~\dive v = 0,~ \la v \ra = 0, ~
\sqrt{-\Delta}^{\,n} v \in \bb{L}^2~\}  \label{defhn} \feq
$$ = \{ v \in \mathbb{D}'~|~~k \sc v_k = 0 ~\forall k \in \Zd,~~ v_0 = 0, ~
\sum_{k \in \Zd_0} |k|^{2 n} | v_k |^2 < + \infty \}~; $$
this is equipped with
the inner product and with the induced norm
\beq \la v | w \ra_n := \la \sqrt{-\Delta}^{\,n} v |  \sqrt{-\Delta}^{\,n} w \ra_{L^2}
= \sum_{k \in \Zd_0} |k|^{2 n} \overline{v_k} \sc w_k, ~~
\| v \|_n := \sqrt{\la v | v \ra_n}~. \label{definner} \feq
Let $n,n', m \in \reali$. One has $\HM{n'} \subset \HM{n}$ and $\| ~\|_n \leqs \|~\|_{n'}$ if $n \leqs n'$;
moreover $\sqrt{-\Delta}^{\,m} \HM{m+n} = \HM{n}$.
By the standard Sobolev lemma, $\HM{n}$ is embedded continuously
in $\mathbb{L}^\infty_{\so}$ if $n > d/2$.
\salto
\textbf{Leray projection.} This is the map
\beq \LP : \mathbb{D}' \vain \mathbb{D}'_\ss~, \qquad v \mapsto \LP v \feq
defined via the Fourier representation
\parn
\vbox{
\beq (\LP v)_k := \LP_k v_k~\mbox{for all $v \in \mathbb{D}', k \in \Zd$}, \feq
\beq \LP_k : \complessi^d \mapsto \complessi^d~\mbox{the orthogonal projection
of $\complessi^d$ onto $k^{\perp}$}. \label{deflpk} \feq
}
Of course, orthogonality in $\complessi^d$ is defined in terms
of the inner product sending $a, b \in \complessi^d$ into $\overline{a}
\sc b$; $k^{\perp}$ is
the orthogonal complement of $k$, i.e., $k^{\perp} := \{ a \in \complessi^d~|~k \sc a = 0 \}$.
If $c \in \complessi^d$, one has
\beq \LP_ k c = c - {k \sc c \over |k|^2}\, k \quad \mbox{for all $k \in \Zd_0$}~, \qquad
\LP_0 c = c~. \label{lpk} \feq
From the Fourier representation
it is evident that
\beq \LP \mathbb{D} = \mathbb{D}'_{\ss}~, \quad
\LP \restriction \mathbb{D}'_{\ss} = {\boma{1}}_{{\mathbb{D'}_{\ss}}}~, \quad
\LP \mathbb{D}'_{\zz} = \mathbb{D}'_{\so}~, \quad \LP \mathbb{L}^2 = \mathbb{L}^2_\ss~, \quad
\LP \mathbb{L}^2_{\zz} = \mathbb{L}^2_{\so}~, \feq
\beq \| \LP v \|_{L^2} \leqs \| v \|_{L^2}~,
\qquad \la \LP v | w \ra_{L^2} =  \la v | \LP w \ra_{L^2} \qquad
\mbox{for $v, w \in \mathbb{L}^2$}~. \feq
\salto
\textbf{The NS bilinear map $\boma{\PPP}$.}
We are now ready to define precisely the map \rref{mapp}. Let us
consider two vector fields
\beq v \in \mathbb{L}^2~, \qquad  w \in \mathbb{D}'~\mbox{such that
$\partial_r w \in \mathbb{L}^2$ for $r = 1,...,d$} \label{sucht} \feq
(which implies $w \in \mathbb{L}^2$). Then we can define
the vector field $v \sc \partial w$ of components
$(v \sc \partial w)_s := \sum_{r=1}^d v_r \partial_r w_s$, that
fulfills
\beq v \sc \partial w \in \mathbb{L}^1~; \feq
we can define as well
\beq \PPP(v, w) := \LP(v \sc \partial w) \in \LP \mathbb{L}^1~. \feq
The Fourier coefficients of these vector fields are obtained by elementary manipulations,
and are as follows (see, e.g., \cite{cog}):
\beq (v \sc \, \partial w)_k =  {i \over (2 \pi)^{d/2}}
\sum_{h \in \Zd} (v_{h} \sc \, (k - h)) w_{k - h}
 \label{dainf} \feq
\beq \PPP(v,w)_k =
\left(\LP(v \sc \partial w)\right)_k = {i \over (2 \pi)^{d/2}}
\sum_{h \in \Zd} (v_{h} \sc \, (k - h)) \LP_k w_{k - h}
\label{infert} \feq
for all $k \in \Zd$. (In the above sums over $h$,
one can replace $\Zd$ with $\Zd_{k} = \Zd \setminus \{k\}$;
if $v$ has zero mean, one can replace $\Zd$ with $\Zd_{0 k} = \Zd
\setminus \{0, k\}$.)
One also proves that
\beq \la v \sc \partial w \ra = \la \PPP(v,w) \ra = 0 \qquad \mbox{if $\dive v=0$} \label{dueto} \feq
(see again \cite{cog}, Lemma 2.1).
Of course, the maps sending $v, w$ as in Eq.\,\rref{sucht} into $v \sc \partial w$ and
$\PPP(v,w)$ are bilinear. Let us go on making the stronger assumption
\beq v \in \mathbb{L}^\infty, \quad w~ \mbox{as in \rref{sucht}}~; \label{suchtt} \feq
then
\beq v \sc \partial w \in \mathbb{L}^2, \quad \PPP(v, w) \in \mathbb{L}^2_{\ss} \label{220} \feq
and, on account of \rref{dueto},
\beq v \sc \partial w \in \mathbb{L}^2_\zz,~~\PPP(v, w) \in \mathbb{L}^2_{\so} \qquad
\mbox{if $\dive v = 0$}~.\label{onac} \feq
Let us also mention that
\beq \la v \sc \partial w | w \ra_{L^2} = 0 \qquad \mbox{if $\dive v = 0$}~, \label{treto} \feq
\beq \la \PPP(v,w) | w \ra_{L^2} = 0 \qquad \mbox{if $\dive v=0$, $\dive w = 0$}~. \label{quato} \feq
Concerning Eq.\,\rref{treto} see, e.g., Lemma 2.3 of \cite{cog}; once we have
\rref{treto}, assuming $\dive w=0$ we infer $w = \LP w$ and
$0 = \la v \sc \partial w | \LP w \ra_{L^2} =
\la \LP(v \sc \partial w) | w \ra_{L^2} = \la \PPP(v,w)|w \ra_{L^2}$,
whence Eq.\,\rref{quato}.
\salto
\textbf{The bilinear maps $\boma{P_{h \ell}}$ and their norms.} Eq.\,\rref{infert} contains the expression
$(v_{h} \sc \, (k-h)) \LP_{k} w_{k-h}$ which has the form
$(a \sc \ell) \LP_{h + \ell} b$ where $\ell := k-h \in \Zd$
and $a := v_h, b := w_{\ell}$; if $\dive v=0$ and
$\dive w=0$ we have $a \in h^{\perp}$,
$b \in \ell^{\perp}$. We fix the attention on the
normalized expression $\dd{a \sc \ell \over |\ell|} \LP_{h + \ell} b$
as a function of $a, b$; more precisely we consider, for $h, \ell \in \Zd_0$, the map
\beq P_{h \ell} : h^{\perp} \times \ell^{\perp} \vain ({h + \ell})^{\perp}, \quad
(a, b) \mapsto P_{h \ell}(a, b) := {a \sc \ell \over |\ell|} \, \LP_{h + \ell} \, b~. \label{phel} \feq
This is a bilinear map between the finite dimensional spaces
indicated above (all of them subspaces of $\complessi^d$);
of course, Eq.\,\rref{infert} (and the remarks that follow it) indicate that
\beq \PPP(v,w)_k =
{i \over (2 \pi)^{d/2}} \sum_{h \in \Zd_{0 k}} |k - h| P_{h, k - h}(v_{h}, w_{k - h})  \label{inferttt} \feq
$$ \mbox{if $v, \partial_1 w,..., \partial_d w \in \mathbb{L}^2$, $v_0=0$
and $\dive v = \dive w = 0$}~. $$
For arbitrary $h, \ell \in \Zd_0$ we can introduce the norm
\beq | P_{h \ell} | := \min \{ Q \in [0,+\infty)~|~|P_{h \ell}(a, b) |
\leqs Q |a| |b |~\mbox{for all $a \in h^\perp$, $b \in \ell^\perp$} \}~. \label{norm}\feq
The above norm can be computed explicitly, as shown in Appendix \ref{appehel}. Indeed,
denoting with $\te_{q r} \in [0,\pi]$ the convex angle between any two
vectors $q, r \in \reali^d \setminus \{0 \}$
({\footnote{Of course, $\cos \te_{q r} = \dd{q \sc r \over |q| |r|}$ and $\sin \te_{q r} =
\sqrt{1 - \dd{(q \sc r)^2 \over |q|^2 |r|^2}}$.}}), we have
\beq |P_{h \ell}| = \left\{ \barray{ll} \sin \te_{h \ell}  & \mbox{if $d \geqs 3$}\,, \\
\sin \te_{h \ell} \cos \te_{h + \ell, \ell} & \mbox{if $d = 2$}
\farray \right.
\label{eqnorm}
\feq
(where $\te_{h + \ell, \ell}$ indicates any angle in $[0,\pi]$ if
$h+\ell=0$; in this situation $\te_{h \ell}=\pi$, so $\sin \te_{h \ell} \cos \te_{h + \ell, \ell}=0$).
In any case we have the bounds
\beq |P_{h \ell}| \leqs \sin \te_{h \ell} \leqs 1~, \label{bou} \feq
to be used in the sequel according to convenience. Let us also
remark that Eq.\,\rref{eqnorm} implies
\beq | P_{\ell h} | = | P_{h \ell} | \quad \mbox{if $d \geqs 3$}~. \label{recall} \feq
\salto
\textbf{An obvious remark.} As already declared in the Introduction,
in this paper we are mainly interested
in $\PPP(v,w)$ for $v \in \HM{p}$, $w \in \HM{p+1}$ and $p > d/2$.
In this case all the conditions on $v, w$ appearing
in Eqs.\,\rref{sucht}\,\rref{suchtt} and \rref{onac}-\rref{quato}
are satisfied (in particular, $v \in \bb{L}^\infty$ by the Sobolev
embedding lemma); thus $v \sc \partial w$, $\PPP(v,w)$
are well defined and possess all the
properties listed in Eqs.\,\rref{dainf}-\rref{dueto} and \rref{220}-\rref{quato}.
This suffices to infer $\PPP(v,w) \in \mathbb{L}^2_{\so}$;
the stronger statement $\PPP(v,w) \in \HM{p}$ is proved
explicitly in the next section.
\section{The inequality \rref{basineqa}; upper bounds for its sharp constant $\boma{K_{p n}}$}
\label{secbasic}
In this section we systematically refer to the maps $P_{h \ell}$ defined
in \rref{phel}, and to their norms described by Eqs.\,\rref{norm}\,\rref{eqnorm}. Moreover,
we consider
\beq \mbox{$p, n \in \reali$ such that $p \geqs n > d/2$}~. \label{31} \feq
\begin{prop}
\label{lemk}
\textbf{Proposition.}
One can define a function
\beq \KK_{p n d} \equiv \KK_{p n} : \Zd_0 \vain (0,+\infty), \label{kknd} \feq
$$ k \mapsto  \KK_{p n}(k) := 4 |k|^{2 p}
\sum_{h \in \Zd_{0 k}} {|P_{h, k - h}|^2 \over (|h|^{p} |k-h|^{n} + |h|^n |k-h|^{p})^2}~; $$
in fact, the sum written above is finite for all $k \in \Zd_0$.
Moreover
\beq \sup_{k \in \Zd_0} \KK_{p n}(k) \leqs 2^{2 p + 2} \zeta_{2 n}~,
\qquad \zeta_{2 n} := \sum_{h \in \Zd_{0}} {1 \over |h|^{2 n}}  \label{zetad} \feq
(note that $\zeta_{2 n} < + \infty$, since $2 n > d$).
\end{prop}
\textbf{Proof.}
Let $k \in \Zd_0$. The sum in Eq.\,\rref{kknd}
and $\KK_{p n}(k)$ certainly exist as elements of $[0,+\infty]$;
the same can be said of the other sums appearing in the proof. \par
Hereafter we derive an upper bound on $\KK_{p n}(k)$ yielding Eq.\,\rref{zetad}
(and ensuring, \textsl{a fortiori}, the finiteness of $\KK_{p n}(k)$).
To this purpose we note the following: for all $h \in \Zd$ one has $k = (k-h) + h$, whence
$|k| \leqs |k - h| + |h|$ and
\beq |k|^{2 p} \leqs (|k - h| + |h|)^{2 p} \leqs 2^{2 p - 1} (|k-h|^{2 p} + |h|^{2 p}) \label{eqand} \feq
(in the last inequality we have used the fact that $(x + y)^{q } \leqs 2^{q-1} (x^{q} + y^{q})$
for $q \in [1, +\infty)$ and $x, y \in [0,+\infty)$).
Inserting \rref{eqand} and the bound
$|P_{h, k - h}| \leqs 1$ (see Eq.\,\rref{bou}) into the definition
\rref{kknd} of $\KK_{p n}(k)$ we get
$$ \KK_{p n}(k) \leqs 2^{2 p + 1}
\Big(\sum_{h \in \Zd_{0 k}} {|k-h|^{2 p}
\over (|h|^{p} |k-h|^{n} + |h|^n |k-h|^{p})^2} +
\sum_{h \in \Zd_{0 k}} {|h|^{2 p}
\over (|h|^{p} |k-h|^{n} + |h|^n |k-h|^{p})^2} \Big). $$
The second sum becomes the first one after a change of
variable $h \vain k-h$; thus
$$ \KK_{p n}(k) \leqs 2^{2 p + 2}
\sum_{h \in \Zd_{0 k}} {|k-h|^{2 p}
\over (|h|^{p} |k-h|^{n} + |h|^n |k-h|^{p})^2}  $$
and since $|h|^{p} |k-h|^{n} + |h|^n |k-h|^{p}
\geqs |h|^n |k-h|^{p}$ we get
\beq \KK_{p n}(k) \leqs 2^{2 p + 2}
\sum_{h \in \Zd_{0}} {1 \over |h|^{2 n}} = 2^{2 p + 2} \zeta_{2 n}~. \label{seeb} \feq
This concludes the proof. \fine
\begin{rema}
\label{remcit}
\textbf{Remarks.}
(i) Let us generalize a remark presented in \cite{cok}
about the constants $K_{p n}$ for $p=n$. To this purpose,
if $r \in \{1,...,d\}$ and $\sigma$ is any permutation of $\{1,...,d\}$,
we define the reflection operator $R_r$ and the permutation operator $P_\si$ setting
\beq R_r, P_\si : \reali^d \vain \reali^d ~, \label{refper} \feq
$$ R_r(k_1,.., k_r,...,k_d) := (k_1,...,-k_r,...,k_d)~, \qquad P_{\si}(k_1,...,k_d) :=
(k_{\si(1)},...,k_{\si(d)})~; $$
these are orthogonal operators (with respect to the inner product $\sc$ of $\reali^d$),
sending $\Zd_0$ into itself. One easily checks that the function $\KK_{p n}$
in \rref{kknd} fulfills
\beq \KK_{p n}(R_r k) = \KK_{p n}(k)~, \qquad \KK_{p n}(P_\sigma k) = \KK_{p n}(k)~
\qquad \mbox{for each $k \in \Zd_0$} \label{claim} \feq
(indeed, the norms $|k|,|h|,|k-h|,|P_{h, k-h}|$
in the definition of $\KK_{p n}(k)$ do not change
if an orthogonal operator is applied to $h$ and $k$).
Due to \rref{claim}, the computation of $\KK_{p n}(k)$
can always be reduced to the case $k_1 \geqs k_2 \geqs ... \geqs k_d \geqs 0$.
\parn
(ii) Typically, the bound \rref{zetad} on $\sup \KK_{p n}$ is very rough;
in the subsequent Sections \ref{sectk} and \ref{sec345} we present
much more accurate estimates on this sup, based on a
lengthy analysis of the function $\KK_{p n}$. As an example,
let $d=n=3$, $p=10$. Then the bound \rref{zetad} reads
$\sup_{k \in \Zt_0} \KK_{10, 3}(k) \leqs 3.53 \times 10^{7}$;
on the other hand, the methods of Section \ref{sectk} and
their numerical implementation
in Section \ref{sec345} give $\sup_{k \in \Zt_0} \KK_{10, 3}(k) \leqs 3.27 \times 10^4$.
Nevertheless, the bound \rref{zetad} is not useless; we return to it at the end of this
section (see Corollary \ref{corrough}) and, especially, in Section \ref{secinf}.
\end{rema}
\begin{prop}
\label{maink}
\textbf{Theorem.} Let $v \in \HM{p}$ and $w \in \HM{p+1}$. Then $\PPP(v,w) \in \HM{p}$ and
\beq \| \PPP(v, w) \|_p \leqs {1 \over 2} \KP_{p n} ( \| v \|_p \| w \|_{n+1} + \| v \|_n \| w \|_{p+1})
\label{basineqap}~, \feq
where
\beq \KP_{p n d} \equiv \KP_{p n}
:= {1 \over (2 \pi)^{d/2}} \sqrt{\sup_{k \in \Zd_0} \KK_{p n}(k)} \label{cantake} \feq
and $\KK_{p n}$ is the function defined by \rref{kknd}.
So, the inequality \rref{basineqa} holds and its sharp constant $K_{p n d} \equiv K_{p n}$ is such that
\beq K_{p n} \leqs \KP_{p n}~. \label{ubk} \feq
\end{prop}
\textbf{Proof.} Let us start from the relation \rref{inferttt}
$$ \PPP(v,w)_k =
{i \over (2 \pi)^{d/2}} \sum_{h \in \Zd_{0 k}} |k-h| P_{h, k - h}(v_{h}, w_{k - h})~; $$
we note that $\PPP(v,w)$ has zero mean due to \rref{dueto}, and is divergence free
by construction. Let us fix $k \in \Zd_0$; from
Eqs.\,\rref{inferttt}\,\rref{norm} we infer
\beq |\PPP(v,w)_k | \leqs
{1 \over (2 \pi)^{d/2}} \sum_{h \in \Zd_{0 k}} |k-h| | P_{h, k - h}| |v_{h}| |w_{k - h}|   \feq
$$ = {1 \over 2 (2 \pi)^{d/2}} \sum_{h \in \Zd_{0 k}}
{2 | P_{h, k - h}|  \over |h|^p |k-h|^{n} + |h|^n ||k-h|^{p}}
(|h|^p |k-h|^{n+1} + |h|^n |k-h|^{p+1}) |v_{h}| |w_{k - h}|~, $$
and the Cauchy inequality $\sum_h a_h b_h \leqs \sqrt{ \sum_h a_h^2} \sqrt{ \sum_h b^2_h}$
($a_h, b_h \in [0,+\infty)$) gives
\parn
\vbox{
\beq |\PPP(v,w)_k | \leqs  {1 \over 2 (2 \pi)^{d/2}}
\sqrt{\sum_{h \in \Zd_{0 k}}
{4 | P_{h, k - h}|^2  \over (|h|^p |k-h|^{n} + |h|^n |k-h|^{p})^2} } \feq
$$ \times \sqrt{
\sum_{h \in \Zd_{0 k}}
(|h|^p |k-h|^{n+1} + |h|^n ||k-h|^{p+1})^2 |v_{h}|^2 |w_{k - h}|^2
}~.
$$}
\noindent
Multiplying by $|k|^p$ and comparing with the definition
\rref{kknd} of $\KK_{p n}$ we see that
$$ |k|^p |\PPP(v,w)_k | \leqs  {\sqrt{\KK_{p n}(k)}\over 2 (2 \pi)^{d/2}}
\sqrt{ \sum_{h \in \Zd_{0 k}}
(|h|^p |k-h|^{n+1} + |h|^n ||k-h|^{p+1})^2 |v_{h}|^2 |w_{k - h}|^2}~,
$$
i.e.,
\parn
\vbox{
\beq |k|^p |\PPP(v,w)_k | \leqs  {\sqrt{\KK_{p n}(k)}\over 2 (2 \pi)^{d/2}}
\, \sqrt{ \sum_{h \in \Zd_{0 k}} (a_{k h} + b_{k h})^2}~, \feq
$$ a_{k h} := |h|^p |k-h|^{n+1} |v_{h}| |w_{k - h}|~, \qquad
b_{k h} := |h|^n |k-h|^{p+1} |v_{h}| |w_{k - h}|~. $$
}
But $\sqrt{ \sum_{h \in \Zd_{0 k}} (a_{k h} + b_{k h})^2}
\leqs \sqrt{ \sum_{h \in \Zd_{0 k}} a_{k h}^2} + \sqrt{ \sum_{h \in \Zd_{0 k}} b_{k h}^2}$, so
\par
\vbox{
\beq |k|^p |\PPP(v,w)_k | \leqs   {\sqrt{\KK_{p n}(k)}\over 2 (2 \pi)^{d/2}}
 \, (q_k + p_k)~, \label{deqk} \feq
$$ q_k := \!\!\!\sqrt{ \sum_{h \in \Zd_{0 k}} |h|^{2 p} |k-h|^{2(n+1)} |v_{h}|^2 |w_{k - h}|^2~} ,~
p_k := \!\!\!\sqrt{ \sum_{h \in \Zd_{0 k}} |h|^{2 n} |k-h|^{2(p+1)} |v_{h}|^2 |w_{k - h}|^2}. $$}
\noindent
To go on we note that, according to \rref{cantake},
$(2 \pi)^{-d/2} \sqrt{\KK_{p n}(k)} \leqs \KP_{p n}$; thus
\beq |k|^p |\PPP(v,w)_k | \leqs  {\KP_{p n} \over 2} \, (q_k + p_k)~, \feq
which implies
\beq {~} \hspace{-0.7cm} \sqrt{\sum_{k \in \Zd_0} |k|^{2 p} |\PPP(v,w)_k |^2} \leqs  {\KP_{p n} \over 2}
\sqrt{\sum_{k \in \Zd_0} (q_k + p_k)^2}
\leqs {\KP_{p n} \over 2} \Big( \sqrt{\sum_{k \in \Zd_0} q^2_k}  + \sqrt{\sum_{k \in \Zd_0} p_k^2}\,\Big).
\label{sumov} \feq
On the other hand, the definition of $q_k$ in \rref{deqk}
gives
\beq \sum_{k \in \Zd_0} q^2_k = \sum_{(k,h) \in \Zd_0 \times \Zd_0, k \neq h}
|h|^{2 p} |k-h|^{2(n+1)} |v_{h}|^2 |w_{k - h}|^2 \label{manip} \feq
$$ = {~} \hspace{-0.5cm} \sum_{(h,\ell) \in \Zd_0 \times \Zd_0, \ell \neq -h}
\hspace{-0.2cm} |h|^{2 p} |v_{h}|^2 |\ell|^{2(n+1)}  |w_{\ell}|^2
\leqs \hspace{-0.2cm} \sum_{(h,\ell) \in \Zd_0 \times \Zd_0} \hspace{-0.2cm}
|h|^{2 p}  |v_{h}|^2  |\ell|^{2(n+1)} |w_{\ell}|^2 = \| v \|^2_{p} \| w \|^2_{n+1}~.
$$
Similarly
\beq \sum_{k \in \Zd_0} p^2_k \leqs \| v \|^2_{n} \| w \|^2_{p+1}~; \feq
inserting these results into \rref{sumov} we get
\beq \sqrt{\sum_{k \in \Zd_0} |k|^{2 p} |\PPP(v,w)_k |^2} \leqs  {\KP_{p n} \over 2}
\big(\| v \|_{p} \| w \|_{n+1} +  \| v \|_{n} \| w \|_{p+1}\big)~. \feq
This proves that $\PPP(v,w) \in \HM{p}$ (we have already noted the vanishing
of the mean and divergence of this vector field); moreover
\rref{basineqap} is found to hold, with $\KP_{p n}$ as in
Eq.\,\rref{cantake}.
\fine
Theorem \ref{maink} gives an upper bound on $K_{p n}$ in terms of $\sup \KK_{p n}$;
we have anticipated that an accurate evaluation of this sup
requires a lengthy analysis, occupying
Section \ref{sectk}. However, at present we have the bound in Proposition \ref{lemk}
on $\sup \KK_{p n}$, whose roughness has been emphasized
in Remark \ref{remcit} (ii). Using this rough estimate,
we obtain the following from the cited propositions.
\begin{prop}
\label{corrough}
\textbf{Corollary.} The sharp constant $K_{p n}$ of \rref{basineqa} has the bound
\beq K_{p n} \leqs K^{\la+\ra}_{p n}~, \qquad
K^{\la+\ra}_{p n} := {2^{p + 1} \over (2 \pi)^{d/2}} \sqrt{\zeta_{2 n}}
\label{kpp} \feq
(with $\zeta_{2 n}$ as in \rref{zetad}; note that $(K^{\la+\ra}_{p n})^{1/p}
\vain 2$ for fixed $d, n$ and $p \vain + \infty$).
\end{prop}
\textbf{Proof}. In fact, Eqs.\,\rref{cantake}\,\rref{ubk} and \rref{zetad} give
$$ K_{p n} \leqs {1 \over (2 \pi)^{d/2}} \sqrt{ \sup_{k \in \Zd_0} \KK_{p n}(k)} \leqs
{1 \over (2 \pi)^{d/2}} \sqrt{2^{2 p + 2} \zeta_{2 n}} ~, $$
whence the thesis \rref{kpp}. \fine
In spite of its roughness, the upper bound \rref{kpp} on $K_{p n}$ has its own
theoretical interest; in fact, as shown in Section \ref{secinf}, the combination
of \rref{kpp} with a suitable lower bound can be used to evaluate $\lim_{p
\vain + \infty} (K_{p n})^{1/p}$.
\section{The inequality \rref{katineqa}; upper bounds for its sharp constant $\boma{G_{p n}}$}
\label{seckato}
The derivation of the generalized Kato inequality \rref{katineqa}
proposed hereafter is rather similar, in the special case $p=n$,
to the one given in \cite{cog}. Both in \cite{cog} and herein,
we refine and make a bit more quantitative some basic
ideas expressed by Temam \cite{Tem} and Constantin-Foias (see \cite{CoFo}, Chapter 10).
\par
Let us start from an elementary inequality, very similar to
some relations presented in \cite{CoFo} \cite{Tem}, whose proof is reported only for completeness.
\begin{prop}
\textbf{Lemma.} Consider a real $p \geqs 1$. Then
\beq \Big| |b|^p - |a|^p \Big| \leqs p \, |b-a| \max(|b|, |a|)^{p-1}~
\qquad \mbox{for $a, b \in \reali^d$}~. \label{inequab} \feq
\end{prop}
\textbf{Proof.} It suffices to prove the inequality \rref{inequab} with the
assumptions
\beq a, b \in \reali^d, \qquad  0 \not \in [a,b]~, \label{segment} \feq
where $[a,b]$ is the segment of $\reali^d$ with endpoints $a,b$;
the inequality is subsequently extended to the case $0 \in [a,b]$ by
elementary continuity considerations. Assuming \rref{segment}, let us
consider the function
\beq F_p : \reali^d \setminus \{0 \} \vain \reali~, \qquad u \mapsto F_p(u) := |u|^p~, \feq
which is $C^\infty$ with
\beq \mbox{grad} F_p(u) = p |u|^{p-2} u~. \feq
We have
$$  \Big| |b|^p - |a|^p \Big| = \Big| F_p(b) - F_p(a) \Big|
\leqs |b-a| \max_{u \in [a,b]} |\mbox{grad} F_p(u)| $$
$$ = p |b-a| \max_{u \in [a,b]} |u|^{p-1}
= p |b-a| \max(|b|, |a|)^{p-1}~. $$
\fine
From here to the end of the section we consider
\beq \mbox{$p, n \in \reali$ such that $p \geqs n > d/2 + 1$}~.  \label{41} \feq
\begin{prop} \textbf{Proposition.}
\label{lemg}
One can define a function
\beq \GG_{p n d} \equiv \GG_{p n} : \Zd_0 \vain (0,+\infty), \label{ggnd} \feq
$$ k \mapsto  \GG_{p n}(k) :=
4 \sum_{h \in \Zd_{0 k}}
{(|k|^p - |k - h|^p)^2 | P_{h, k - h}|^2  \over (|h|^p |k-h|^{n-1} + |h|^n |k-h|^{p-1})^2}~; $$
in fact, the sum written above is finite for all $k \in \Zd_0$. Moreover
\beq \sup_{k \in \Zd_0} \GG_{p n}(k) \leqs 2^{2 p} p^2 \zeta_{2 n - 2}~, \qquad
\zeta_{2 n - 2} := \sum_{h \in \Zd_{0}} {1 \over |h|^{2 n-2}} \label{zetadg} \feq
(note that $\zeta_{2 n - 2} < + \infty$, since $2 n - 2 > d$).
\end{prop}
\textbf{Proof.} Let $k \in \Zd_0$. The sum in Eq.\,\rref{ggnd}
and $\GG_{p n}(k)$ certainly exist as elements of $[0,+\infty]$; the same can be said of the other sums that
appearing in the proof. \par
Hereafter we derive an upper bound on $\GG_{p n}(k)$, yielding Eq.\,\rref{zetadg}
(and ensuring, \textsl{a fortiori}, the finiteness of $\GG_{p n}(k)$).
To this purpose we note that, for all $h \in \Zd$, the
inequality \rref{inequab} with $b =k$ and $a=k-h$ gives
\beq (|k|^p - |k - h|^p)^2 \leqs p^2 |h|^2 \max(|k|, |k-h|)^{2 p - 2} \feq
$$ \leqs  p^2 |h|^2 \max(|k-h| + |h|, |k-h||)^{2 p - 2} = p^2 |h|^2 (|k-h|+ |h||)^{2 p - 2} $$
$$ \leqs 2^{2 p - 3} p^2 |h|^2 (|k-h|^{2 p - 2} + |h|^{2 p - 2}) $$
(concerning the last inequality, see the comment after Eq.\,\rref{eqand}). Inserting
this inequality and the relation $|P_{h, k-h}| \leqs 1$ (see Eq.\,\rref{bou}) into
the definition \rref{ggnd} of $\GG_{p n}(k)$ we get
\parn
\vbox{
\beq \GG_{p n}(k) \leqs 2^{2 p - 1} p^2
\sum_{h \in \Zd_{0 k}}
{|h|^2 (|k-h|^{2 p - 2} + |h|^{2 p - 2})  \over (|h|^p |k-h|^{n-1} + |h|^n |k-h|^{p-1})^2} \feq
$$ = 2^{2 p - 1} p^2
\left(
\sum_{h \in \Zd_{0 k}}
{|k-h|^{2 p - 2} \over (|h|^{p-1} |k-h|^{n-1} + |h|^{n-1} |k-h|^{p-1})^2} \right. $$
$$ + \left. \sum_{h \in \Zd_{0 k}}
{|h|^{2 p - 2}  \over (|h|^{p-1} |k-h|^{n-1} + |h|^{n-1} |k-h|^{p-1})^2}
\right)~;$$}
\noindent
the second sum above becomes the first one after a change of
variable $h \vain k-h$, and thus
\beq \GG_{p n}(k) \leqs  2^{2 p} p^2
\sum_{h \in \Zd_{0 k}}
{|k-h|^{2 p - 2} \over (|h|^{p-1} |k-h|^{n-1} + |h|^{n-1} |k-h|^{p-1})^2}~. \feq
Since $|h|^{p-1} |k-h|^{n-1} + |h|^{n-1} |k-h|^{p-1}
\geqs |h|^{n-1} |k-h|^{p-1}$, we get
\beq \GG_{p n}(k)
\leqs 2^{2 p} p^2
\sum_{h \in \Zd_{0}} {1 \over |h|^{2 n-2}} = 2^{2 p} p^2 \zeta_{2 n - 2}~; \label{seebg} \feq
this concludes the proof. \fine
\begin{rema}
\label{remcitg}
\textbf{Remarks.} The forthcoming comments (i)(ii) are quite
similar to Remarks \ref{remcit} about the function $\KK_{p n}$ and its sup. \parn
(i) Let $r \in \{1,...,d\}$, and let $\sigma$ denote a permutation
of $\{1,...,d\}$; denoting with $R_r$ and $P_{\sigma}$ the reflection
and permutation operators \rref{refper}, we have
\beq \GG_{p n}(R_r k) = \GG_{p n}(k)~, \qquad \GG_{p n}(P_\sigma k) = \GG_{p n}(k)~
\qquad \mbox{for each $k \in \Zd_0$}~. \label{claimg} \feq
So, the computation of $\GG_{p n}(k)$ can be reduced to the case
$k_1 \geqs k_2 ... \geqs k_d \geqs 0$.
\parn
(ii) The bound \rref{zetadg} on $\sup \GG_{p n}$ is very rough;
in Sections \ref{sectg} and \ref{sec345} we present
much more accurate estimates on this sup, based on a
lengthy analysis of the function $\GG_{p n}$. As an example,
let $d=n=3$, $p=10$. Then the bound \rref{zetad} reads
$\sup_{k \in \Zt_0} \GG_{10, 3}(k) \leqs 1.74 \times 10^9$;
on the other hand, the methods of Section \ref{sectg} and
their numerical implementation
in Section \ref{sec345} give $\sup_{k \in \Zt_0} \GG_{10, 3}(k) \leqs 6.64 \times 10^4$.
\end{rema}
\begin{prop}
\label{maing}
\textbf{Theorem.} Let $v \in \HM{p}$ and $w \in \HM{p+1}$
(so that $\PPP(v, w) \in \HM{p}$). Then
\beq | \la \PPP(v, w) | w \ra_p | \leqs
{1 \over 2} \GP_{p n} (\| v \|_p \| w \|_n + \| v \|_n \| w \|_p)\| w \|_p~,
\label{katineqap} \feq
where
\beq \GP_{p n d} \equiv \GP_{p n} :=
{1 \over (2 \pi)^{d/2}} \sqrt{\sup_{k \in \Zd_0} \GG_{p n}(k)} \label{cantakeg} \feq
and $\GG_{p n}$ is the function defined by \rref{ggnd}.
Therefore, the inequality \rref{katineqa} holds and its sharp constant
$G_{p n d} \equiv G_{p n}$ fulfills
\beq G_{p n} \leqs \GP_{p n}~. \label{ubg} \feq
\end{prop}
\textbf{Proof.}
We fix $v \in \HM{p}, w \in \HM{p+1}$ and proceed in several steps. \parn
\textsl{Step 1. We have $\PPP(v,w)  \in \HM{p}$, $\sqrt{-\Delta}^{\, p} \PPP(v,w)
\in \bb{L}^{2}_{\so}$, $\sqrt{-\Delta}^{\, p} w \in \HM{1}$ and
$\PPP(v, \sqrt{-\Delta}^{\,p} w)  \in \bb{L}^{2}_{\so}$;
furthermore the vector field}
\beq z := \sqrt{-\Delta}^{\, p} \PPP(v, w) - \PPP(v, \sqrt{-\Delta}^{\,p} w) \in \bb{L}^{2}_{\so}
\label{defze}\feq
\textsl{fulfills the equality}
\beq  \la \PPP(v, w) |w \ra_{p} = \la z | \sqrt{-\Delta}^{\,p} w \ra_{L^2}~, \label{wehave} \feq
\textsl{which implies}
\beq | \la \PPP(v, w) |w \ra_{p} | \leqs \| z \|_{L^2} \| w \|_p~. \label{wehave2} \feq
The statement $\PPP(v,w)  \in \HM{p}$ is known after Theorem
\ref{maink}, and the statement $\sqrt{-\Delta}^{\, p} \PPP(v,w)
\in \bb{L}^{2}_{\so}$ is just a reformulation of it.
Our assumption $v \in \HM{p}$ implies $v \in \bb{L}^{\infty}_{\so}$, by the Sobolev embedding;
of course $\sqrt{-\Delta}^{\,p}$ sends $\HM{p+1}$ into $\HM{1}$,
thus $\sqrt{-\Delta}^{\,p} w \in \HM{1}$. Now, applying the
second result \rref{onac} with $w$ replaced
by $\sqrt{-\Delta}^{\, p} w$ we obtain
$\PPP(v, \sqrt{-\Delta}^{\,p} w)  \in \bb{L}^{2}_{\so}$.
To go on, we note that the definition of $\la~|~\ra_p$ gives
\beq \la \PPP(v,w) | w \ra_p = \la \sqrt{-\Delta}^{\, p} \PPP(v,w) | \sqrt{-\Delta}^{\, p} w \ra_{L^2} \label{sum1} \feq
and that Eq.\,\rref{quato} with $w$ replaced by $\sqrt{-\Delta}^{\, p} w$ gives
\beq 0 = \la \PPP(v, \sqrt{-\Delta}^{\, p} w) | \sqrt{-\Delta}^{\, p} w \ra_{L^2}~. \label{sum2} \feq
Subtracting Eq.\,\rref{sum2} from Eq.\,\rref{sum1} we obtain the thesis
\rref{wehave}, with $z$ given by
\rref{defze}. Eq.\,\rref{wehave} and the Schwartz inequality yield
$| \la \PPP(v,w) |w \ra_{p} |$ $\leqs \| z \|_{L^2} \| \sqrt{-\Delta}^{\,p} w \|_{L^2}$
$= \| z \|_{L^2} \| w \|_p$\,, as in \rref{wehave2}.
\parn
\textsl{Step 2. The vector field $z$ in \rref{defze} has
Fourier coefficients
\beq z_{k} = - {i \over (2 \pi)^{d/2}}
\sum_{h \in \Zd_{0 k}} (|k|^p - |k - h|^p) |k-h| P_{h, k - h}(v_h,  w_{k - h})
\quad \mbox{for all $k \in \Zd_0$}~.
\label{znk} \feq}
In fact
\beq {~} \hspace{-0.4cm} [\sqrt{-\Delta}^{\, p} \PPP(v,w)]_k = |k|^p \PPP(v,w)_k =
- {i \over (2 \pi)^{d/2}} |k|^{p} \sum_{h \in \Zd_{0 k}} |k-h| P_{h, k - h}(v_{h}, w_{k - h})\,;  \feq
the last equality follows from Eq.\,\rref{inferttt}. Using the same equation we get
\beq \PPP(v, \sqrt{-\Delta}^{\,p} w)_k =
- {i \over (2 \pi)^{d/2}} \sum_{h \in \Zd_{0 k}} |k-h| P_{h, k - h}(v_{h}, [\sqrt{-\Delta}^{\, p} w]_{k - h}) \feq
$$ = - {i \over (2 \pi)^{d/2}} \sum_{h \in \Zd_{0 k}} |k-h| P_{h, k - h}(v_{h}, |k-h|^p w_{k - h})~. $$
The last two equations and the definition of $z$ yield the thesis \rref{znk}.
\parn
\textsl{Step 3. Estimating the Fourier coefficients of $z$}.
Let $k \in \Zd_0$; Eq.\,\rref{znk} implies
\beq |z_k| \leqs {1 \over (2 \pi)^{d/2}}
\sum_{h \in \Zd_{0 k}} \Big||k|^p - |k - h|^p \Big| \, |k-h| \,|P_{h, k - h}|\,|v_{h}|\, |w_{k - h}|
\label{zznk} \feq
$$ = {1 \over 2 (2 \pi)^{d/2}} \sum_{h \in \Zd_{0 k}}
{2 \Big||k|^p - |k - h|^p \Big| | P_{h, k - h}|  \over |h|^p |k-h|^{n-1} + |h|^n ||k-h|^{p-1}}
(|h|^p |k-h|^{n} + |h|^n |k-h|^{p}) |v_{h}| |w_{k - h}|~. $$
Now, the Cauchy inequality $\sum_h a_h b_h \leqs \sqrt{ \sum_h a_h^2} \sqrt{ \sum_h b^2_h}$
($a_h, b_h \in [0,+\infty)$) gives
\beq | z_k | \leqs  {1 \over 2 (2 \pi)^{d/2}}
\sqrt{\sum_{h \in \Zd_{0 k}}
{4 (|k|^p - |k - h|^p)^2 | P_{h, k - h}|^2  \over (|h|^p |k-h|^{n-1} + |h|^n ||k-h|^{p-1})^2} } \feq
$$ \times \sqrt{
\sum_{h \in \Zd_{0 k}}
(|h|^p |k-h|^{n} + |h|^n ||k-h|^{p})^2 |v_{h}|^2 |w_{k - h}|^2
}~;
$$
comparing with the definition
\rref{ggnd} of $\GG_{p n}$ we see that
\beq | z_k | \leqs  {1 \over 2 (2 \pi)^{d/2}} \sqrt{\GG_{p n}(k)} \sqrt{
\sum_{h \in \Zd_{0 k}}
(|h|^p |k-h|^{n} + |h|^n ||k-h|^{p})^2 |v_{h}|^2 |w_{k - h}|^2 }~.
\label{bouzk}
\feq
The last inequality has the form
\beq |z _k | \leqs  {1 \over 2 (2 \pi)^{d/2}}
\sqrt{\GG_{p n}(k)}
\sqrt{ \sum_{h \in \Zd_{0 k}} (a_{k h} + b_{k h})^2}~, \feq
$$ a_{k h} := |h|^p |k-h|^{n} |v_{h}| |w_{k - h}|~, \qquad
b_{k h} := |h|^n |k-h|^{p} |v_{h}| |w_{k - h}|~. $$
But $\sqrt{ \sum_{h \in \Zd_{0 k}} (a_{k h} + b_{k h})^2}
\leqs \sqrt{ \sum_{h \in \Zd_{0 k}} a_{k h}^2} + \sqrt{ \sum_{h \in \Zd_{0 k}} b_{k h}^2}$, so
\par
\vbox{
\beq | z_k | \leqs  {1 \over 2 (2 \pi)^{d/2}}
\sqrt{\GG_{p n}(k)} (q_k + p_k)~, \label{deqgk} \feq
$$ q_k := \sqrt{ \sum_{h \in \Zd_{0 k}} |h|^{2 p} |k-h|^{2 n} |v_{h}|^2 |w_{k - h}|^2~} ~,~~
p_k := \sqrt{ \sum_{h \in \Zd_{0 k}} |h|^{2 n} |k-h|^{2 p} |v_{h}|^2 |w_{k - h}|^2}~. $$}
\noindent
To go on we note that, according to \rref{cantakeg},
$(2 \pi)^{-d/2} \sqrt{\GG_{p n}(k)} \leqs \GP_{p n}$; thus
\beq | z_k | \leqs  {\GP_{p n} \over 2} (q_k + p_k)~. \feq
\textsl{Step 4. Estimating $\| z \|_{L^2}$}.
With $q_k$ and $p_k$ as in Eq.\,\rref{deqgk}, we have
\beq \| z \|_{L^2} = \sqrt{\sum_{k \in \Zd_0} | z_k |^2} \leqs  {\GP_{p n} \over 2}
\sqrt{\sum_{k \in \Zd_0} (q_k + p_k)^2}
\leqs {\GP_{p n}\over 2} \Big( \sqrt{\sum_{k \in \Zd_0} q^2_k}  + \sqrt{\sum_{k \in \Zd_0} p_k^2}\,\Big).
\label{sumog} \feq
On the other hand, manipulations very similar to the ones in
Eq.\,\rref{manip} give
\beq \sum_{k \in \Zd_0} q^2_k \leqs \| v \|^2_{p} \| w \|^2_{n}~, \qquad
\sum_{k \in \Zd_0} p^2_k \leqs \| v \|^2_{n} \| w \|^2_{p}~; \feq
inserting this result into \rref{sumog} we get
\beq \| z \|_{L^2} \leqs  {\GP_{p n} \over 2}
\big(\| v \|_{p} \| w \|_{n} +  \| v \|_{n} \| w \|_{p}\big)~. \label{bzld} \feq
\textsl{Step 5. Conclusion of the proof.} We return to
Eq.\,\rref{wehave2}, and insert therein the bound
\rref{bzld} for $\| z \|_{L^2}$. This gives
the inequality \rref{katineqap}, with $\GP_{p n}$ as in \rref{cantakeg}. \fine
Theorem \ref{maink} gives an upper bound on $G_{p n}$ in terms of $\sup \GG_{p n}$
that, as anticipated, will be the subject of accurate estimates in
Section \ref{sectg}. For the moment, using the rough bound of Proposition \ref{lemg}
on $\sup \GG_{p n}$ we obtain:
\begin{prop}
\label{corroughg}
\textbf{Corollary.} The sharp constant $G_{p n}$ of \rref{katineqa} has the bound
\beq G_{p n} \leqs G^{\la+\ra}_{p n}~, \qquad
G^{\la+\ra}_{p n} := {2^{p} \, p \over (2 \pi)^{d/2}} \sqrt{\zeta_{2 n-2}}
\label{gpp} \feq
(with $\zeta_{2 n-2}$ as in \rref{zetadg}; note that $(G^{\la+\ra}_{p n})^{1/p}
\vain 2$ for fixed $d,n$ and $p \vain + \infty$).
\end{prop}
\textbf{Proof}. In fact, Eqs.\,\rref{cantakeg}\,\rref{ubg} and \rref{zetadg} give
$$ G_{p n} \leqs {1 \over (2 \pi)^{d/2}} \sqrt{ \sup_{k \in \Zd_0} G_{p n}(k)} \leqs
{1 \over (2 \pi)^{d/2}} \sqrt{2^{2 p} p^2 \zeta_{2 n-2}} ~, $$
whence the thesis \rref{gpp}. \fine
Similarly to the rough bound \rref{kpp} on $K_{p n}$, the present
bound \rref{gpp} wil be useful in Section \ref{secinf} to evaluate the
$p \vain + \infty$ limit of
$(G_{p n})^{1/p}$.
\section{Some tools preparing the analysis of the functions
$\boma{\KK_{p n}}$ and $\boma{\GG_{p n}}$}
\label{proveg}
As anticipated, in Sections \ref{sectk}-\ref{sec345}
we will show how to compute accurately the functions $\KK_{p n}$, $\GG_{p n}$
of Eqs. \rref{kknd} \rref{ggnd} and their sups; here we introduce
some tools devised for this purpose. \par
First of all we fix some notations, to be used throughout the rest of the paper.
\begin{prop}
\textbf{Definition.}
\label{versk}
(i) $\delta_{a b}$ is the Kronecker delta ($\delta_{a b} := 1$
if $a=b$ and $\delta_{a b} := 0$ if $a \neq b$). \parn
(ii) $\teta : \reali \vain \{0,1\}$
is the Heaviside function such that $\teta(z) := 0$
if $z < 0$ and $\teta(z) := 1$ if $z \geqs 0$. \parn
(iii) $\Gamma$ is the Euler Gamma function,
$\binom{\,\cdot\,}{\,\cdot\,}$ are the binomial coefficients. \parn
(iv) $\Sd$ denotes the unit spherical hypersurface in $\reali^d$, i.e.,
$\Sd := \{ u \in \reali^d~|~|u| = 1 \}$.
For each $q \in \reali^d \setminus \{0\}$, the versor of $q$ is
$\vers{q}  := \dd{q \over |q|} \in \Sd$.
\end{prop}
In the sequel we also maintain the following notation,
already introduced in Section \ref{notations}: for all $q, r \in \reali^d \setminus \{0 \}$,
$\te_{q r} \in [0,\pi]$ denotes
the convex angle between $q$ and $r$ (so that $\cos \te_{q r} = \vers{q} \sc \vers{r}$).
\begin{prop}
\label{lemmaf}
\textbf{Lemma.} For any function $f : \Zd_0 \vain \reali$ and $k \in \Zd_0$,
$\ro \in (1,+\infty)$, one has
\beq \sum_{h \in \Zd_{0 k}, |h| < \ro \op |k-h| < \ro} f(h)
= \sum_{h \in \Zd_{0 k},|h| < \ro} \big[ f(h) + \teta(|k-h| - \ro) f(k-h) \big]~.
\label{tesif} \feq
\end{prop}
\textbf{Proof.} See \cite{cok}. \fine
\begin{prop}
\label{lemks}
\textbf{Lemma.} For any $p,n \in \reali$ with $p \geqs n > 1$, the following holds. \parn
(i) Consider the function
\parn
\vbox{
\beq b_{p n} : [0,4] \times [0,1] \vain [0,+\infty)~, \label{debn} \feq
$$ b_{p n}(z,u) :=
\dd{2 z (4 - z) (1 - u)^{2 n} u^{2 n} (1 - z u + z u^2)^{p} \over
[(1 - u)^{2 n} + u^{2 n}] [ (1 - u)^p u^n + (1-u)^n u^p]^2} \qquad \mbox{if $u \in (0,1)$}, $$
$$ b_{p n}(z,0) := b_{p n}(z,1) := \dd{2 z (4 - z) \over 1 + 3 \, \delta_{p n}}~. $$}
This is well defined and continuous, which implies the existence of
\beq B_{p n} := \max_{z \in [0,4], u \in [0,1]} b_{p n}(z,u) > 0~. \label{debnk} \feq
(ii) Given $h, \ell \in \reali^d \setminus \{0 \}$, consider the convex angle $\te_{h \ell}$
and define $z \in [0,4]$, $u \in (0,1)$ through the equation
\beq \cos \te_{h \ell} = 1 - \dd{z \over 2}~,  \qquad |h | = \dd{u \over 1 - u} |\ell|~;
\label{asineq} \feq
then
\beq {|h+\ell|^{2 p} \sin^2 \te_{h \ell} \over (|h|^{p} |\ell|^{n} + |h|^{n} |\ell|^{p})^2}
= {b_{p n}(z,u) \over 8} \left( {1 \over |h|^{2 n}} + {1 \over |\ell|^{2 n}} \right)~. \label{onthesiok} \feq
This implies
\beq {|h+\ell|^{2 p}  \sin^2 \te_{h \ell} \over (|h|^{p} |\ell|^{n} + |h|^{n} |\ell|^{p})^2}
\leqs {B_{p n} \over 8} \left( {1 \over |h|^{2 n}} + {1 \over |\ell|^{2 n}} \right)~. \label{onthesik} \feq
\end{prop}
\textbf{Proof.} (i) Trivial (in particular it is not difficult to check that
$b_{p n}(z,0)$ $= \lim_{u \vain 0}$ $b_{p n}(z,u)$ and $b_{p n}(z,1) = \lim_{u \vain 1} b_{p n}(z,u)$). \parn
(ii) Consider the quantity
\beq {|h+\ell|^{2 p} \sin^2 \te_{h \ell} \over (|h|^{p} |\ell|^{n} + |h|^{n} |\ell|^{p})^2}
\left( {1 \over |h|^{2 n}} + {1 \over |\ell|^{2 n}} \right)^{-1} \label{rat2} \feq
and express it via the relations $\sin^2 \te_{h \ell} = 1 - \cos^2 \te_{h \ell}$ and
\beq |h + \ell| = \sqrt{|h|^2 + 2 |h| |\ell| \cos \te_{h \ell} + |\ell|^2}~; \label{esprel} \feq
subsequently, express $\cos \te_{h \ell}$ and $|h|$ via Eq.\,\rref{asineq}.
After tedious manipulations it is found that \rref{rat2} equals $\dd{b_{p n}(z,u) \over 8}$,
and Eq.\,\rref{onthesiok} is proved. Eq.\,\rref{onthesik} is an obvious consequence. \fine
\begin{rema}
\label{rembpn}
\textbf{Remark.} Let $b_{p n}$ be defined as in the previous lemma.
It is readily found that the derivatives $\partial b_{p n}/\partial z$,
$\partial b_{p n}/\partial u$ vanish at $(z,u)  = ({4 \over p+2}, {1 \over 2})$, and
\beq b_{p n}\left({4 \over p+2}, {1 \over 2} \right) =
{4^{p + 1} \over p+2} \left( {p + 1 \over p+2} \right)^{p+1}~; \feq
moreover, the hessian of $b_{p n}$ at $(z,u)  = ({4 \over p+2}, {1 \over 2})$ is positive
defined if $p \geqs n > 1$. Considering $B_{p n} := \max_{z \in [0,4], u \in [0,1]} b_{p n}(z,u)$,
we conjecture that
\beq B_{p n} = b_{p n}\left({4 \over p+2}, {1 \over 2} \right) =
{4^{p + 1} \over p+2} \left( {p + 1 \over p+2} \right)^{p+1} \label{eqcong} \feq
for all $(p,n)$ with $p \geqs n > 1$;
note that $B_{p n}$ does not depend on $n$ if the above statement holds.
For given $(p,n)$, statement \rref{eqcong} can be tested using a computer to plot
$b_{p n}$ or to maximize it numerically. In this way we have obtained
that \rref{eqcong} holds for all $(p,n)$ as in Eqs.\,\rref{casesdue}\rref{casestre}
(i.e., in all cases considered in the sequel to exemplify the evaluation
of $\KK_{p n}$).
\end{rema}
\begin{prop}
\label{lemgs}
\textbf{Lemma.} For any $p,n \in \reali$ with $p \geqs n > 1$, the following holds. \parn
(i) Consider the function
\parn
\vbox{
\beq c_{p n} : [0,4] \times [0,1] \vain [0,+\infty)~, \label{decn} \feq
$$ c_{p n}(z,u) :=
\dd{2 z (4 - z) (1 - u)^{2 n + 2} u^{2 n} [(1 - z u + z u^2)^{p/2} - (1 - u)^p  ]^2 \over
[(1 - u)^{2 n} u^2 + (1 - u)^2 u^{2 n}] [ (1 - u)^p u^n + (1-u)^n u^p]^2} \qquad \mbox{if $u \in (0,1)$}, $$
$$ c_{p n}(z,0) := \dd{p^2 z (4 - z) (2 - z)^2 \over 2 ( 1 + 3 \, \delta_{p n})}~,
\qquad c_{p n}(z,1) := \dd{2 z (4 - z) \over 1 + 3 \, \delta_{p n}}~. $$}
This is well defined and continuous, which implies the existence of
\beq C_{p n} := \max_{z \in [0,4], u \in [0,1]} c_{p n}(z,u) > 0~. \label{decin} \feq
(ii) Given $h, \ell \in \reali^d \setminus \{0 \}$
define $z \in [0,4]$, $u \in (0,1)$ as in Eq.\,\rref{asineq}; then
\beq {(|h+\ell|^p - |\ell|^{p})^2 \sin^2 \te_{h \ell} \over (|h|^{p} |\ell|^{n-1} + |h|^{n} |\ell|^{p-1})^2}
= {c_{p n}(z,u) \over 8} \left( {1 \over |h|^{2 n - 2}} + {1 \over |\ell|^{2 n - 2}} \right)~. \label{onthesio} \feq
This implies
\beq {(|h+\ell|^p - |\ell|^{p})^2 \sin^2 \te_{h \ell} \over (|h|^{p} |\ell|^{n-1} + |h|^{n} |\ell|^{p-1})^2}
\leqs {C_{p n} \over 8} \left( {1 \over |h|^{2 n - 2}} + {1 \over |\ell|^{2 n - 2}} \right)~. \label{onthesi} \feq
\end{prop}
\textbf{Proof.} (i) Trivial (in particular it is not difficult to check that
$c_{p n}(z,0)$ $= \lim_{u \vain 0}$ $ c_{p n}(z,u)$ and $c_{p n}(z,1) = \lim_{u \vain 1} c_{p n}(z,u)$). \parn
(ii) Let us consider the ratio
\beq {(|h+\ell|^p - |\ell|^{p})^2 \sin^2 \te_{h \ell} \over (|h|^{p} |\ell|^{n-1} + |h|^{n} |\ell|^{p-1})^2}
\left( {1 \over |h|^{2 n - 2}} + {1 \over |\ell|^{2 n - 2}} \right)^{-1} \label{rat1} \feq
and express it using Eq.\,\rref{esprel};
subsequently, write $\cos \te_{h \ell}$ and $|h|$ as in \rref{asineq}.
After tedious manipulations it is found that the ratio \rref{rat1} equals $\dd{c_{p n}(z,u) \over 8}$,
and Eq.\,\rref{onthesio} is proved. Eq.\,\rref{onthesi} is an obvious consequence. \fine
\begin{rema}
\textbf{Examples.} Let $c_{p n}, C_{p n}$ be defined as in the previous lemma.
For $n= 3$ and $p=3,4,5,10$ we have the following results,
obtained by numerical optimization via {\tt{Mathematica}}:
\beq C_{3 3}= c_{3 3}(0.696034..., 0.464530...) = 14.8144... , \label{valc2} \feq
$$ C_{4 3} = c_{4 4}(0.610279...,0.439178 ...) = 61.1705... , $$
$$ C_{5 3}= c_{5 3}(0.545364...,0.443863 ...) =  229.715... , $$
$$ C_{10, 3} = c_{10, 3}(0.332954...,0.489262 ...) = 1.36660...\times 10^5 \,. $$
\end{rema}
\begin{prop}
\label{deken}
\textbf{Lemma.} Let $p, n \in \reali$, $p \geqs n > 0$; then (i)(ii) hold. \parn
(i) Let us introduce the domain
\beq \Do := \{ (c,\xi) \in \reali^2~|~c \in [-1,1], \, \xi \in [0,+\infty), \, (c,\xi) \neq (1,1) \}
\label{dedo} \feq
and put
\parn
\vbox{
\beq E_{p n} : \Do \vain [0,+\infty)~, \label{deen} \feq
$$ E_{p n}(c,\xi) := {1 - c^2  \over [ (1 - 2 c \xi + \xi^2)^{p/2+ 1/2} +
\xi^{p-n} (1 - 2 c \xi + \xi^2)^{n/2 + 1/2}]^2}~. $$
}
Then the above function is well defined and continuous on $\Do$. \parn
(ii) Let $h, k \in \reali^d \setminus \{0\}$, $h \neq k$ and
consider the convex angles $\te_{h k}, \te_{h, k-h}$. Then
\beq
{ |k|^{2p} \sin^2 \te_{h, k - h}
\over (|h|^{p} |k-h|^{n} + |h|^{n} |k-h|^{p})^2}
= {1 \over |h|^{2 n}} \, E_{p n}\left(\cos \te_{h k}\,, {|h| \over |k|} \right)~.
\label{funk} \feq
\end{prop}
\textbf{Proof.} (i) Trivial. \parn
(ii) The parallelograms of sides $h$ and $k$, $h$ and $k-h$ have the same
area; thus $|h| |k| \sin \te_{h k} = |h| |k-h| \sin \te_{h, k - h}$, whence
\beq \sin \te_{h, k - h} = {|k| \over |k-h|} \sin \te_{h k} =
{|k| \over |k-h|} \sqrt{1 - \cos^2 \te_{h k}}~; \label{iden1} \feq
moreover,
\beq |k - h| = \sqrt{|k|^2 - 2 |k| |h| \cos \te_{h k} + |h|^2} =
|k| \sqrt{1 - 2 \, \cos \te_{h k} {|h| \over |k|}  + {|h|^2 \over |k|^2}}~. \label{iden2} \feq
Let us consider the function in the left hand side of \rref{funk}, and
reexpress it using the identities \rref{iden1} \rref{iden2};
in this way, after some manipulations we obtain Eq.\,\rref{funk}.
\par\fine
\begin{prop}
\label{deken2}
\textbf{Lemma.} Let $p,n \in \reali$, $p \geqs n >0$, and consider the function
$E_{p n} : \Do \vain \reali$ of Lemma \ref{deken}. Introduce the set
\beq \Gamma_{p n} := \{ \ia + (p-n) \ja \,|\, \ia,\ja \in \naturali\} \label{defgamma} \feq
and represent it as an increasing sequence:
\beq \Gamma_{p n} = \{ 0 = \gamma_{p n 0} < \gamma_{p n 1} < \gamma_{p n 2} < ... \}\,. \feq
There are two sequences of functions
\beq Q_{p n j} \in C([-1,1],\reali), ~~ c \mapsto Q_{p n j}(c) \qquad (j \in \naturali) \feq
\beq S_{p n j} \in C(\Do, \reali),~~(c, \xi) \mapsto S_{p n j}(c, \xi)
\qquad (j \in \naturali \setminus \{0 \}) \feq
uniquely determined by the following prescription: for each $m \in \naturali$
one has
\beq E_{p n}(c, \xi) = \sum_{j=0}^m Q_{p n j}(c) \, \xi^{\gamma_{p n j}} + S_{p n, m+1}(c, \xi)
\, \xi^{\gamma_{p n, m+1}} \quad \mbox{for all $(c, \xi) \in \Do$}~. \label{onehask} \feq
Moreover, each function $Q_{p n j}$ is of polynomial type.
\end{prop}
\textbf{Proof.} It suffices to show the following: \parn
(a) for each $m \in \naturali$, there is a unique family of functions $Q_{p n 0},...,Q_{p n m}
\in C([-1,1], \reali)$, $S_{p n, m+1} \in C(\Do, \reali)$ such that \rref{onehask} holds. Moreover,
the functions $Q_{p n j}$ ($j=0,...,m)$ are polynomials; \parn
(b) for $m < m' \in \naturali$, the family
$Q_{p n 0}, ... ,Q_{p n m}, S_{p n, m+1}$ of item (a)
and the family $Q'_{p n 0}, ... ,Q'_{p n m'}, S_{p n, m'+1}$ of item (a)
with $m$ replaced by $m'$ are such that $Q_{p n 0} = Q'_{p n 0}$, ... ,$Q_{p n m} = Q'_{p n m}$. \par
Let us first prove the uniqueness statement in (a), for a given $m \in \naturali$.
To this purpose we note that, given a family as in (a),
Eq.\,\rref{onehask} implies
\beq Q_{p n 0}(c) = E_{p n}(c, 0)~, \label{car0} \feq
$$ Q_{p n j}(c) = \lim_{\xi \vain 0} {1 \over \xi^{\gamma_{p n j}}}
\Big( E_{p n}(c, \xi) - \sum_{\ell=0}^{j-1} Q_{p n \ell} \xi^{\gamma_{p n \ell}} \Big) \quad
\mbox{for $j=1,...,m$}~; $$
this set of recursive relations determines uniquely the functions
$Q_{p n j}$ for  $j=0,...,m$. Once we have uniqueness for the sequence $(Q_{p n j})_{j=0,...,m}$,
uniqueness of $S_{p n,m+1}$ follows noting that
\rref{onehask} implies
\beq S_{p n, m + 1}(c, \xi) = {1 \over \xi^{\gamma_{p n, m+1}}}
\Big( E_{p n}(c, \xi) - \sum_{j=0}^m Q_{p n j}(c) \xi^{\gamma_{p n j}} \Big)
\quad \mbox{for $(c, \xi) \in \Do, \xi \neq 0$} \feq
and that, by the continuity requirement for $S_{p n, m+1}$,
$S_{p n, m+1}(c,0)$ is the $\xi \vain 0$ limit of the right hand
side in the above equation. \par
Now, let us prove statement (b) for given $m < m' \in \naturali$.
To this purpose we note that, besides the characterization
\rref{car0} for $Q_{p n 0}, ..., Q_{p n m}$ we have a similar
characterization for $Q'_{p n 0}, ..., Q'_{p n m'}$; these imply
$Q_{p n 0}(c) = E_{p n}(c,0) = Q'_{p n 0}(c)$,
\hbox{$Q_{p n 1}(c) = \lim_{\xi \vain 0} \xi^{-\gamma_{p n 1}}(E_{p n}(c) - Q_{p n 0}(c))$}
$= \lim_{\xi \vain 0} \xi^{-\gamma_{p n 1}} (E_{p n}(c) - Q'_{p n 0}(c))$
$= Q'_{p n 1}(c)$ and so on, up to $Q_{p n m}(c) = Q'_{p n m}(c)$. \par
Let us pass to prove, for any $m \in \naturali$,
the existence of the functions $Q_{p n 0},...,Q_{p n m}$,
$S_{p n, m+1}$ fulfilling the conditions in (a) and the polynomial nature of the functions
$Q_{p n j}$; for the sake of brevity we discuss the case $p > n$, leaving
to the reader the case $p=n$ which is even simpler.
Let us note that Eq.\,\rref{deen} has the form
\beq E_{p n}(c, \xi) = A_{p n}(c, \xi, \xi^{p - n})~, \label{efabk} \feq
where
\beq A_{p n} : \Do \times [0,+\infty) \vain \reali, \feq
$$ A_{p n}(c,\xi, u) :=
\dd{1 - c^2  \over [ (1 - 2 c \xi + \xi^2)^{p/2 +1/2} + u (1 - 2 c \xi + \xi^2)^{n/2+1/2}]^2}~. $$
It is easily checked that $A_{p n} \in C^\infty(\Do \times [0,+\infty), \reali)$.
Now, consider any $\aa \in \naturali$; by Taylor's
formula of order $\aa$ in the variables $\xi, u$ we can write
\beq A_{p n}(c, \xi, u) = \sum_{\ia,\ja \in \naturali, \ia+\ja \leqs \aa} A_{p n \ia \ja}(c) \, \xi^\ia u^\ja
+  \sum_{\ia,\ja \in \naturali, \ia+\ja= \aa+1} S_{p n \ia \ja}(c, \xi, u) \, \xi^\ia u^\ja  \label{espatk} \feq
$$ \qquad \mbox{for $(c, \xi, u) \in \Do \times [0,+\infty)$} $$
with suitable reminder functions  $S_{p n \ia \ja} \in C(\Do \times [0,+\infty), \reali)$.
The coefficients $A_{p n \ia \ja}(c)$ in the above expansions are related to the derivatives
of $A_{p n}$ at $\xi=0, u=0$; one finds by direct inspection of
the definitions of $A_{p n}$ that these coefficients are
polynomial functions of $c$. Inserting the expansions \rref{espatk}
into Eq.\,\rref{efabk} we get
\beq \label{wegek}
E_{p n}(c,\xi)
= \sum_{\ia,\ja \in \naturali, \ia+\ja \leqs \aa} A_{p n \ia \ja}(c) \xi^{\ia + (p-n)\ja} +
\sum_{\ia,\ja \in \naturali, \ia+\ja=\aa+1} S_{p n \ia \ja}(c,\xi, \xi^{p - n}) \xi^{\ia + (p-n)\ja}~.  \feq
All the exponents of $\xi$ in the above formula belong to the set
$\Gamma_{p n} = \{ 0 = \gamma_{p n 0} < \gamma_{p n 1} < ... \}$
defined by \rref{defgamma}. Now, after fixing $m \in \naturali$ we choose
$\aa \in \naturali$ so that
\beq \ia + (p-n) \ja \geqs \gamma_{p n,m+1}~\mbox{for
all $\ia,\ja \in \naturali$ such that $\ia+\ja = \aa+1$}~; \feq
then Eq.\,\rref{wegek} implies for $E_{p n}$ a representation of the form
\rref{onehask} for this value of $m$, where
\beq Q_{p n j}(c) = \sum_{\ia,\ja \in \naturali,\, \ia+\ja \leqs \aa, \,\ia + (p-n) \ja =
\gamma_{p n j}} \hspace{-0.7cm}
A_{p n \ia \ja}(c) \qquad (j=0,...,m), \feq
\beq S_{p n, m+1}(c, \xi) = \sum_{\ia,\ja,\in \naturali,\, \ia+\ja =\aa+1}
\!\! S_{p n \ia \ja}(c,\xi, \xi^{p - n}) \, \xi^{\ia + (p-n) \ja - \gamma_{p n, m +1}}
\feq
$$ + \sum_{\ia,\ja\, \in \naturali,\, \ia+\ja \leqs \aa,
\,\ia + (p-n) \ja \geqs \gamma_{p n, m+1}}  \hspace{-1cm}
A_{p n \ia \ja}(c) \, \xi^{\ia + (p-n) \ja  - \gamma_{p n, m +1}}~. $$
We note that the functions $Q_{p n j}$ ($j=0,...,m$) are polynomials in $c$
and $S_{p n, m+1}$ is continuous due to the previously mentioned features of $A_{p n \ia \ja}$
and $S_{p n \ia \ja}$. \fine
\begin{prop}
\label{deden}
\textbf{Lemma.} Let $p, n \in \reali$, $p \geqs n > 0$; then (i)(ii) hold. \parn
(i) Let $\Do$ be the domain in Eq.\,\rref{dedo}, and put
\parn
\vbox{
\beq F_{p n} : \Do \vain [0,+\infty)~, \label{defn} \feq
$$ F_{p n}(c,\xi) :=
\dd{1 - c^2  \over [ (1 - 2 c \xi + \xi^2)^{p/2} + \xi^{p-n} (1 - 2 c \xi + \xi^2)^{n/2}]^2}~$$
$$ \times
\Big\{ {(1 - \xi^p)^2 \over 1 - 2 c \xi + \xi^2} + \Big[{1- (1 - 2 c \xi + \xi^2)^{p/2} \over \xi}\Big]^2 \Big\}
\quad \mbox{if $\xi \neq 0$}~, $$
$$ F_{p n}(c, 0) := \dd{(1 - c^2)(1 + p^2 c^2) \over 1 + 3 \, \delta_{p n}}~. $$
}
Then the above function is well defined and continuous on $\Do$. \parn
(ii) Let $h, k \in \reali^d \setminus \{0\}$, $h \neq k$; then
$$
\sin^2 \te_{h, k-h} \left[{(|k|^p - |k - h|^p)^2 \over (|h|^p |k-h|^{n-1} + |h|^n |k-h|^{p-1})^2}
+  {(|k|^p - |h|^p)^2 \over (|k-h|^p |h|^{n-1} + |k-h|^n |h|^{p-1})^2}\right] $$
\beq = {1 \over |h|^{2 n -2}} F_{p n} \Big(\cos \te_{h k}\,, {|h| \over |k|} \Big)~.
\label{fung} \feq
\end{prop}
\textbf{Proof.} (i) Trivial (in particular, it is easy to check that
$F_n(c,0) = \lim_{\xi \vain 0} F_n(c,\xi)$). \parn
(ii) We consider the function in the left hand side of \rref{fung}, and
reexpress it using the identities \rref{iden1} \rref{iden2};
in this way, after some manipulations we obtain Eq.\,\rref{fung}. \fine
\begin{prop}
\label{deden2}
\textbf{Lemma.} Let $p,n \in \reali$, $p \geqs n > 0$, and consider the function
$F_{p n} : \Do \vain \reali$ of Lemma \ref{deden}. Introduce the set
\beq \Lambda_{p n} := \{ \ia + (p-n) \ja + p \,\ella \,|\, \ia,\ja \in \naturali, \ella \in \{0,1,2\} \} \label{deflanda} \feq
and represent it as an increasing sequence:
\beq \Lambda_{p n} = \{ 0 = \lambda_{p n 0} < \lambda_{p n 1} < \lambda_{p n 2} < ... \}~. \feq
There are two sequences of functions
\beq P_{p n j} \in C([-1,1], \reali), ~~ c \mapsto P_{p n j}(c) \qquad (j \in \naturali) \feq
\beq R_{p n j} \in C(\Do, \reali),~~(c, \xi) \mapsto R_{p n j}(c, \xi)
\qquad (j \in \naturali \setminus \{0 \}) \feq
uniquely determined by the following prescription: for each $m \in \naturali$
one has
\beq F_{p n}(c, \xi) = \sum_{j=0}^m P_{p n j}(c) \, \xi^{\lambda_{p n j}}
+ R_{p n, m+1}(c, \xi) \, \xi^{\lambda_{p n, m+1}} \quad \mbox{for all $(c, \xi) \in \Do$}~. \label{onehas} \feq
Moreover, each function $P_{p n j}$ is of polynomial type.
\end{prop}
\textbf{Proof.} It suffices to show the following:  \parn
(a) for each $m \in \naturali$, there is a unique family of functions $P_{p n 0},...,P_{p n m}
\in C([-1,1], \reali)$,
$R_{p n, m+1} \in C(\Do, \reali)$, such that \rref{onehas} holds.
Moreover, the functions $P_{p n j}$ ($j=0,...,m)$ are polynomials; \parn
(b) for $m < m' \in \naturali$, the family
$P_{p n 0},...,P_{p n m}, R_{p n, m+1}$ of item (a)
and the family $P'_{p n 0},...,P'_{p n m'}, R_{p n, m'+1}$ of item (a)
with $m$ replaced by $m'$ are such that $P_{p n 0} = P'_{p n 0}$,...,$P_{p n m} = P'_{p n m}$. \par
Let us first prove the uniqueness statement in (a), for a given $m \in \naturali$.
To this purpose we note that, given a family as in (a),
Eq.\,\rref{onehas} implies
\beq P_{p n 0}(c) = F_{p n}(c, 0)~, \label{carp0} \feq
$$ P_{p n j}(c) = \lim_{\xi \vain 0} {1 \over \xi^{\lambda_{p n j}}}
\Big( F_{p n}(c, \xi) - \sum_{\ell=0}^{j-1} P_{p n \ell} \xi^{\lambda_{p n \ell}} \Big) \qquad
\mbox{for $j =1,...,m$}~, $$
and this set of recursive relations determines uniquely the functions
$P_{p n j}$ for $j=0,...,m$. Once we have uniqueness for the sequence $(P_{p n j})_{j=0,...,m}$,
uniqueness of $R_{p n, m+1}$ follows noting that
\rref{onehas} implies
\beq R_{p n, m + 1}(c, \xi) = {1 \over \xi^{\lambda_{p n, m+1}}}
\Big( F_{p n}(c, \xi) - \sum_{j=0}^m P_{p n j}(c) \xi^{\lambda_{p n j}} \Big)
~\mbox{for $(c, \xi) \in \Do, \xi \neq 0$} \feq
and that, by the continuity requirement for $R_{p n, m+1}$,
$R_{p n, m+1}(c,0)$ is the $\xi \vain 0$ limit of the right hand
side in the above equation. \par
Now, let us prove statement (b) for given $m < m' \in \naturali$.
To this purpose we note that, besides the characterization
\rref{carp0} for $P_{p n 0}, ..., P_{p n m}$ we have a similar
characterization for $P'_{p n 0}, ..., P'_{p n m'}$; these imply
$P_{p n 0}(c) = F_{p n}(c,0) = P'_{p n 0}(c)$,
$P_{p n 1}(c) = \lim_{\xi \vain 0} \xi^{-\lambda_{p n 1}}(F_{p n}(c) - P_{p n 0}(c))$
$= \lim_{\xi \vain 0} \xi^{-\lambda_{p n 1}} (F_{p n}(c) - P'_{p n 0}(c))$
$= P'_{p n 1}(c)$ and so on, up to $P_{p n m}(c) = P'_{p n m}(c)$. \par
Let us pass to prove, for any $m \in \naturali$, the existence of functions $P_{p n 0},...,P_{p n m}$,
$R_{p n, m+1}$ fulfilling the conditions in (a) and the polynomial nature of the functions $P_{p n j}$;
for the sake of brevity we only discuss the case $p > n$, leaving
to the reader the case $p=n$ which is even simpler.
Let us note that Eq.\,\rref{defn} has the form
\beq F_{p n}(c, \xi) =
\Ax_{p n}(c, \xi, \xi^{p - n}) (1 - 2 \xi^p + \xi^{2 p}) + H_{p n}(c, \xi, \xi^{p-n})
\label{efab} \feq
where
\parn
\vbox{
\beq \Ax_{p n}, H_{p n} : \Do \times [0,+\infty) \vain \reali, \feq
$$ \Ax_{p n}(c,\xi, u) :=
\dd{1 - c^2  \over [ (1 - 2 c \xi + \xi^2)^{p/2} + u (1 - 2 c \xi + \xi^2)^{n/2}]^2} {1 \over 1 - 2 c \xi + \xi^2}~, $$
$$ H_{p n}(c,\xi, u) :=
\dd{1 - c^2  \over [ (1 - 2 c \xi + \xi^2)^{p/2} + u (1 - 2 c \xi + \xi^2)^{n/2}]^2}~$$
$$ \times \Big[{1- (1 - 2 c \xi + \xi^2)^{p/2} \over \xi}\Big]^2
\quad \mbox{if $\xi \neq 0$}~, $$
$$ H_{p n}(c,0, u) := \lim_{\xi \vain 0} H_{p n}(c, \xi, u) = p^2 \, {(1 - c^2) \, c^2 \over ( 1 + u)^2}~. $$
}
It is easily checked that $\Ax_{p n}, H_{p n}$ $\in C^\infty(\Do \times [0,+\infty), \reali)$.
Consider any $a \in \naturali$; by Taylor's
formula of order $a$ in the variables $\xi, u$ we can write
\beq \Ax_{p n}(c, \xi, u) = \sum_{\ia,\ja \in \naturali, \ia+\ja \leqs \aa} \Ax_{p n \ia \ja}(c) \xi^\ia u^\ja
+  \sum_{\ia,\ja \in \naturali, \ia+\ja= \aa+1} S_{p n \ia \ja}(c, \xi, u) \xi^\ia u^\ja~, \label{espat} \feq
$$ H_{p n}(c, \xi, u) = \sum_{\ia,\ja \in \naturali, \ia+\ja \leqs \aa} H_{p n \ia \ja}(c) \xi^\ia u^\ja
+  \sum_{\ia,\ja \in \naturali, \ia+\ja= \aa+1} T_{p n \ia \ja}(c, \xi, u) \xi^\ia u^\ja $$
$$ \qquad \mbox{for $(c, \xi, u) \in \Do \times [0,+\infty)$}~, $$
with suitable reminder functions  $S_{p n \ia \ja}, T_{p n \ia \ja} \in
C(\Do \times [0,+\infty), \reali)$.
The coefficients $\Ax_{p n \ia \ja}(c)$ and
$H_{p n \ia \ja}(c)$ in the above expansions are related to the derivatives
of $\Ax_{p n}, H_{p n}$ at $\xi=0, u=0$; one finds by direct inspection of
the definitions of $\Ax_{p n}, H_{p n}$ that these coefficients are
polynomial functions of $c$. Inserting the expansions \rref{espat}
into Eq.\,\rref{efab} we get
\parn
\vbox{
\beq F_{p n}(c,\xi)  \label{wege} \feq
$$ = \sum_{\ia,\ja,\ella \in \naturali, \ia+\ja \leqs \aa, \ella \leqs 2}
C_{p n \ia \ja \ella}(c) \, \xi^{\ia + (p-n)\ja + \ella p} +
\sum_{\ia,\ja,\ella \in \naturali, \ia+\ja=\aa+1, \ella \leqs 2} V_{p n \ia \ja \ella}(c,\xi) \,
\xi^{\ia + (p-n)\ja + \ella p}~, $$
$$ C_{p n \ia \ja 0}(c) := \Ax_{p n \ia \ja}(c) + H_{p n \ia \ja}(c),~~
C_{p n \ia \ja 1}(c) := - 2 \Ax_{p n \ia \ja}(c),~~C_{p n \ia \ja 2}(c) := \Ax_{p n \ia \ja}(c), $$
$$ V_{p n \ia \ja 0}(c,\xi) := S_{p n \ia \ja}(c,\xi,\xi^{p-n}) + T_{p n \ia \ja}(c,\xi,\xi^{p-n}),~~
V_{p n \ia \ja 1}(c) := - 2 S_{p n \ia \ja}(c,\xi,\xi^{p-n}),$$
$$ V_{p n \ia \ja 2}(c,\xi) := S_{p n \ia \ja}(c,\xi,\xi^{p-n})~. $$
}
\noindent
All the exponents of $\xi$ in the above formula belong to the set
$\Lambda_{p n} = \{ 0 = \lambda_{p n 0} < \lambda_{p n 1} < ... \}$
defined by \rref{deflanda}. Now, after fixing $m \in \naturali$ we choose
$a \in \naturali$ so that
\beq \ia + (p-n) \ja + \ella p \geqs \lambda_{p n,m+1}~\mbox{for
all $\ia,\ja,\ella \in \naturali$ such that $\ia+\ja = \aa+1$, $\ella \leqs 2$}~; \feq
then Eq.\,\rref{wege} implies for $F_{p n}$ a representation of the form
\rref{onehas} for this value of $m$, where
\beq P_{p n j}(c) = \sum_{\ia,\ja,\ella \in \naturali,\, \ia+\ja \leqs a, \ella \leqs 2, \,\ia + (p-n) \ja + \ella p =
\lambda_{p n j}}
C_{p n \ia \ja \ella}(c)~, \qquad (j=0,...,m) \feq
\beq R_{p n, m+1}(c, \xi) = \sum_{\ia,\ja,\ella \in \naturali,\, \ia+\ja =\aa+1, \ella \leqs 2}
V_{p n \ia \ja \ella}(c,\xi) \, \xi^{\ia + (p-n) \ja + \ella p - \lambda_{p n, m +1}}
\feq
$$ + \sum_{\ia,\ja,\ella \in \naturali,\, \ia+\ja \leqs \aa, \ella \leqs 2,
\,\ia + (p-n) \ja + \ella p \geqs \lambda_{p n, m+1}}
C_{p n \ia \ja \ella}(c) \, \xi^{\ia + (p-n) \ja + \ella p - \lambda_{p n, m +1}}~. $$
We note that the functions $P_{p n j}$ ($j=0,...,m$) are polynomials in $c$
and $R_{p n, m+1}$ is continuous due to the previously mentioned features of $\Ax_{p n \ia \ja}$, $H_{p n \ia \ja}$,
$S_{p n \ia \ja}$, $T_{p n \ia \ja}$. \fine
\begin{prop}
\label{lemz}
\textbf{Lemma.} Consider a real $\nu > d$. For any real $\ro >  2 \sqrt{d}$, one has
\beq \sum_{h \in \Zd, |h| \geqs \ro} {1 \over |h|^{\nu}}~\leqs
{2 \pi^{d/2} \over \Gamma(d/2)}
\sum_{i=0}^{d-1} \left( \barray{c} d - 1 \\ i \farray \right) {d^{d/2-1/2-i/2} \over (\nu - 1 - i)
(\ro - 2 \sqrt{d})^{\nu - 1 - i}}~. \label{desnu}\feq
\end{prop}
\textbf{Proof.} This is just Lemma C.2 of  \cite{accau} (with the variable $\lan$ of the cited reference
related to $\ro$ by $\lan = \ro - 2 \sqrt{d}$). \fine
The forthcoming statement uses the notation $\widehat{~}$ of
Definition \ref{versk} to indicate versors.
\begin{prop}
\label{lemcos4}
\textbf{Lemma.} Let $\ell \in \naturali$,   $\ro \in (1,+\infty)$, $\varphi: [1,\ro) \vain \reali$
and $k \in \reali^d \setminus \{0 \}$. Then
\beq
\sum_{h \in \Zd_0, |h| < \ro} \varphi(|h|) \cos^\ell \! \te_{h k} = P_{\varphi \ell}(\vk)~, \label{idew4} \feq
where $P_{\varphi \ell}$ is the following polynomial function on the spherical hypersurface $\Sd$:
\beq P_{\varphi \ell} : \Sd \vain \reali,~u \mapsto
P_{\varphi \ell}(u) := \hspace{-0.4cm} \sum_{i_1,...,i_d \in \naturali,i_1 + ... + i_d=\ell}
{\ell! \over i_1! ... i_d !} M_{\varphi_, i_1,...,i_d} {u_1}^{i_1}... {u_d}^{i_d}~,
\label{mmfi} \feq
\beq M_{\varphi, i_1,...,i_d} := \sum_{h \in \Zd_0, |h| < \rho}
\varphi(|h|) {\vh_1}^{\, i_1} ... {\vh_d}^{i_d}
\label{mfi} \feq
(in the above $u_r$ and $\vh_r$ stand for the $r$-th components of $u$ and $\vh$;
${u_r}^{i_r}$ and ${\vh_r}^{\,i_r}$ indicate their powers with exponent $i_r$). One has
\beq M_{\varphi, i_1,...,i_d} = 0 \qquad \mbox{if $i_r$ is odd for some $r \in \{1,...,d\}$}~, \label{sim1} \feq
\beq M_{\varphi, i_{\sigma(1)},...,i_{\sigma(d)}} = M_{\varphi, i_1,...,i_d} \qquad
\mbox{for each permutation $\sigma$ of $\{1,...,d\}$} \label{sim2} \feq
(so, the computation of the coefficients $M_{\varphi, i_1,...,i_d}$ can be
reduced to cases with $i_1 \leqs i_2 ... \leqs i_d$ and $i_r$ even for all $r$).
The previous facts imply
\beq P_{\varphi \ell}(u) = 0 \quad \mbox{for all $u \in \Sd$,} \quad \mbox{if $\ell$ is odd} \label{elod} \feq
and, in the case $\ell=2$,
\beq P_{\varphi 2}(u) = \mbox{constant} = {1 \over d} \sum_{h \in \Zd_0, |h| < \rho}
\varphi(|h|) \qquad \mbox{for $u \in \Sd$}~. \label{ell2} \feq
\end{prop}
\textbf{Proof.} We have $\cos \te_{h k} = \vh \sc \vk = \vh_1 \vk_1 + ... + \vh_d \vk_d$~,
which implies
\beq \cos^{\ell} \te_{h k} = (\vh_1 \vk_1 + ... + \vh_d \vk_d)^\ell
= \!\!\! \sum_{i_1,...,i_d \in \naturali,i_1 + ... + i_d=\ell} {\ell! \over i_1! ... i_d !}
{\vh_1}^{i_1} ...   {\vh_d}^{i_d}
{\vk_1}^{i_1} ...   {\vk_d}^{i_d}~.\feq
Multiplying this relation by $\varphi(|h|)$ and summing over $h$ we easily obtain Eqs.\,\rref{idew4}-\rref{mfi};
the definition \rref{mfi} of the coefficients $M_{\varphi, i_1,...,i_d}$
gives the relations \rref{sim1} \rref{sim2} by elementary considerations of symmetry.
Now, assume
$\ell$ is odd; then each one of the coefficients $M_{\varphi, i_1,...,i_d}$ appearing
in Eq.\,\rref{mmfi} is zero, because the list $(i_1,...,i_d)$ has some
odd element and \rref{sim1} can be applied, so we obtain Eq.\,\rref{elod}. \par
Let us pass to the case $\ell=2$; the only nonzero coefficients involved in \rref{mmfi} are
$$ M_{\varphi, 2,0,...,0} = M_{\varphi, 0,2,0,...,0} = ... = M_{\varphi,0,...,0,2} $$
and can be determined noting that
$$ M_{\varphi, 2,0,...,0} + ... + M_{\varphi,0,...,0,2}
= \!\! \sum_{h \in \Zd_0, |h| < \rho} \varphi(|h|) ({\vh_1}^{2} + ... +  {\vh_d}^{2}) =
\!\! \sum_{h \in \Zd_0, |h| < \rho} \varphi(|h|) ~.$$
In conclusion, the nonzero coefficients in \rref{mmfi} for $\ell=2$ are
\beq M_{\varphi, 2,0,...,0} = M_{\varphi, 0,2,0,...,0} = ... = M_{\varphi,0,...,0,2} =
{1 \over d} \sum_{h \in \Zd_0, |h| < \rho} \varphi(|h|)~,  \feq
and Eq.\,\rref{ell2} follows immediately.
\fine
We remark that statement \rref{ell2} in the above Lemma
is equivalent to Lemma A.5 in the arXiv version of \cite{cok}.
\section{The function $\boma{\KK_{p n}}$}
\label{sectk}
As in Section \ref{secbasic}, we consider $p, n \in \reali$ such that $p \geqs n > d/2$.
For $k \in \Zd_0$, we recall the definition \rref{kknd}
$$ \KK_{p n}(k) := 4 |k|^{2 p}
\sum_{h \in \Zd_{0 k}} {|P_{h, k - h}|^2 \over (|h|^{p} |k-h|^{n} + |h|^n |k-h|^{p})^2}~. $$
The general term of the above sum over $h$
is large when its denominator is small, which happens when $h$ is close to zero
or to $k$. Therefore, it is reasonable to approximate the infinite sum
in \rref{kknd} with a finite sum over the union of two balls of centers
$0$, $k$ and a suitable radius $\rho$. Such a finite
sum is the main character of the forthcoming proposition,
where it is denoted with $\Km_{p n}(k)$; the proposition
uses $\Km_{p n}(k)$ and other ingredients to estimate
$\KK_{p n}(k)$ and its $\sup$ for $k \in \Zd_0$.
\par
Let us remark that, for $p=n$, the results presented hereafter become
very similar to those appearing in Proposition B.1 of
\cite{cok} (extended arXiv version).
\begin{prop}
\label{prokknd}
\textbf{Proposition.}
Let us choose a ''cutoff''
\beq \ro \in (2 \sqrt{d},+\infty)~, \label{cutoffk} \feq
a ``factor''
\beq \M \in (1,+\infty) \label{mfacc} \feq
and an ``order''
\beq m \in \naturali~; \label{order} \feq
then the following holds
(with the functions and quantities
$\Km_{p n}$, $\dK_{p n}$, ..., $Y_{p n}$, $\Kk_{p n}$ mentioned in the sequel
depending parametrically on $\ro$, $\M$, $m$ and $d$: $\Km_{p n}(k) \equiv \Km_{p n \ro d}(k)$,
$\dK_{p n} \equiv \dK_{p n \ro d}$, ..., $Y_{p n} \equiv Y_{p n \rho m d},
\Kk_{p n} \equiv \Kk_{p n \ro \M m d}$). \parn
(i) The function $\KK_{p n}$ fulfills the inequalities
\beq \Km_{p n}(k) < \KK_{p n}(k) \leqs \Km_{p n}(k) + \dK_{p n}
~~\mbox{for all $k \in \Zd_0$}~. \label{dkmnd} \feq
Here
\beq \Km_{p n}(k) :=
4 |k|^{2 p} \hspace{-0.4cm}
\sum_{h \in \Zd_{0 k}, |h| < \ro \op |k-h| < \ro} {|P_{h, k - h}|^2 \over (|h|^{p} |k-h|^{n} + |h|^n |k-h|^{p})^2}~
~; \label{deakm} \feq
this function can be reexpressed as
\beq \Km_{p n}(k) = 4 |k|^{2 p} \sum_{h \in \Zd_{0 k}, |h | < \ro }\hspace{-0.3cm}
{| P_{h, k - h}|^2 + \teta(|k-h| - \ro) | P_{k-h, h}|^2
\over (|h|^p |k-h|^{n} + |h|^{n} |k-h|^{p})^2}  \label{kmnd} \feq
(with $\teta$ the Heaviside function, see Definition \ref{versk}; recall
that $| P_{k-h, h}| = | P_{h,k-h}|$ if $d \geqs 3$, due to \rref{recall}).
If $|k| \geqs 2 \ro$, in Eq.\,\rref{kmnd} one can replace $\Zd_{0 k}$ with $\Zd_0$
and $\teta(|k-h| - \ro)$ with $1$. Moreover
\beq \dK_{p n} :=
{2 \pi^{d/2} B_{p n} \over \Gamma(d/2)}
\sum_{i=0}^{d-1} \Big( \barray{c} d - 1 \\ i \farray \Big) {d^{d/2-1/2-i/2} \over (2 n - 1 - i)
(\ro - 2 \sqrt{d})^{2 n - 1 - i}}~, \label{dedek} \feq
with $B_{p n}$ as in \rref{debnk}. \parn
(ii) Consider the reflection and permutation operators
$R_r$, $P_{\si}$ defined by \rref{refper}. Then
\beq \Km_{p n}(R_r k) = \Km_{p n}(k)~, \quad \Km_{p n}(P_\sigma k) = \Km_{p n}(k)~
\qquad \mbox{for each $k \in \Zd_0$}~. \label{claib} \feq
\parn
(iii) Denoting with $E_{p n}$ the function in Lemma \ref{deken}, one has
\beq \Km_{p n}(k) \leqs 8 \sum_{h \in \Zd_0, |h| < \rho} {1 \over |h|^{2 n}}
E_{p n}\left(\cos \te_{h k}\,, {|h| \over |k|} \right) \qquad \mbox{for all $k \in \Zd_0$}. \label{kf} \feq
Now, consider the sequence of exponents
$0 = \gamma_{p n 0} < \gamma_{p n 1} < \gamma_{p n 2} ...$ and the sequences of
polynomials $Q_{p n j}$ ($j \in \naturali)$
and functions $S_{p n j}$ ($j \in \naturali \setminus \{ 0 \}$) involved
in the expansion of $E_{p n}$ according to Lemma \ref{deken2}; then
\beq \Km_{p n}(k) \leqs 8 \sum_{j=0}^{m} {\QQ_{p n j}(\vk) \over |k|^{\gamma_{p n j}}}  +
8 \,{ V_{p n} Y_{p n} \over |k|^{\gamma_{p n, m+1}}} \leqs \Kk_{p n}
\qquad \mbox{for $k \in \Zd,~ |k| \geqs \M \rho$}~. \label{gwk} \feq
Here, we recall, $\widehat{k}$ is the versor of $k$ (Definition \ref{versk}).
$\QQ_{p n j}$ are the functions defined as follows on the spherical hypersurface $\Sd$:
\beq \QQ_{p n j} : \Sd \vain \reali~, \qquad u \mapsto \QQ_{p n j}(u) :=
\sum_{h \in \Zd_0, |h| < \rho} {Q_{p n j}(\cos \te_{h u}) \over |h|^{2 n - \gamma_{p n j}}}
\label{qqpnj} \feq
(these are polynomials in the components of $u$, which can be computed
using Lemma \ref{lemcos4}). Moreover
\beq V_{p n} := \max_{c \in [-1,1], \xi \in [0,1/\M]} S_{p n, m+1}(c, \xi)\,, \label{wpnmk} \feq
\beq Y_{p n} := \sum_{h \in \Zd_0, |h| < \rho} {1 \over |h|^{2 n - \gamma_{p n, m+1}}}~.
\label{ypnmk} \feq
Finally,
\beq \Kk_{p n} := 8 \, \max_{u \in \Sd, \, \epsilon \in [0,1/(\M \rho)]}
\left( \sum_{j=0}^{m} \epsilon^{\gamma_{p n j}} \QQ_{p n j}(u) +
V_{p n} Y_{p n} \epsilon^{\gamma_{p n, m+1}} \right)~. \label{kpnm} \feq
(iv) Items (i) and (iii) imply
\beq {~} \hskip -0.4cm \sup_{k \in \Zd_0} \KK_{p n}(k)
\leqs \max \left( \max_{k \in \Zd_0, |k| < \M \rho} \Km_{p n}(k), \Kk_{p n}
\right) + \dK_{p n}~. \label{impifik} \feq
\end{prop}
The proof of the above statements will be given after the following comment.
\begin{rema}
\label{remkn}
\textbf{Remark}. From Theorem \ref{maink} we know that
the sharp constant $K_{p n}$ of
the inequality \rref{basineqa} fulfills
$K_{p n} \leqs \KP_{p n} := (2 \pi)^{-d/2} \sqrt{\sup_{k \in \Zd_0} \KK_{p n}(k)}$. Of course,
using for $\sup_{k \in \Zd_0} \KK_{p n}(k)$ the bound \rref{impifik},
we conclude
\beq K_{p n} \leqs \KP_{p n} \leqs \KPP_{p n} := {1 \over (2 \pi)^{d/2}}
\sqrt{\max \left( \max_{k \in \Zd_0, |k| < \M \rho} \Km_{p n}(k), \Kk_{p n}
\right) + \dK_{p n}}~~. \label{kpnbou} \feq
The bound $\KPP_{p n}$ is suitable for computer implementation,
a fact discussed in Section \ref{sec345}. From the statements in Proposition \ref{prokknd},
it is clear that $\KPP_{p n}$ gives a good approximation of $\KP_{p n}$
under two conditions: \parn
(a) $\Kk_{p n}$ (which is a theoretical upper bound on $\sup_{|k| \geqs \M \rho} \Km_{p n}(k)$),
must be smaller or not too larger than $\max_{|k| < \M \rho} \Km_{p n}(k)$.
If $\Kk_{p n}$ is smaller, one has
$\max \left(\max_{|k| < \M \rho} \Km_{p n}(k), \Kk_{p n} \right)$
$= \max_{|k| < \M \rho} \Km_{p n}(k)$ $= \sup_{k} \Km_{p n}(k)$;
if $\Kk_{p n}$ is a bit larger, $\max \left(\max_{|k| < \M \rho} \Km_{p n}(k), \Kk_{p n} \right)$
$= \Kk_{p n}$ is a good upper approximant for $\sup_{k} \Km_{p n}(k)$. \parn
(b) $\dK_{p n}$, which binds the distance between $\Km_{p n}(k)$ and $\KK_{p n}(k)$
for all $k$, must be small in comparison with $\max \left(\max_{|k| < \M \rho} \Km_{p n}(k), \Kk_{p n} \right)$;
under this condition and the one in item (a), $\max \left(\max_{|k| < \M \rho} \Km_{p n}(k), \Kk_{p n} \right)$
is a good approximation for $\sup_{k} \KK_{p n}(k)$.
\end{rema}
\textbf{Proof of Proposition \ref{prokknd}}. We fix $\ro, \M, m$ as in \rref{cutoffk}-\rref{order}.
Our argument is divided in several steps; more precisely,
Steps 1-4 prove the statements in (i) while
Steps 5, 6-7 and 8 are about the statements in (ii),(iii) and (iv), respectively.
The assumption \rref{cutoffk} $\ro > 2 \sqrt{d}$ is essential in Step 3.
\parn
\textsl{Step 1. One has
\beq \KK_{p n}(k)= \Km_{p n}(k) + \DK_{p n}(k)~\quad \mbox{for all $k \in \Zd_0$}~, \label{decompk} \feq
where, as in \rref{deakm},
$\Km_{p n}(k) :=\dd{
4 |k|^{2 p}
\sum_{h \in \Zd_{0 k}, |h| < \ro \op |k-h| < \ro}
{|P_{h, k - h}|^2 \over (|h|^{p} |k-h|^{n} + |h|^n |k-h|^{p})^2}~}$, while
\beq \DK_{p n}(k) := 4 |k|^{2 p} \hspace{-0.4cm} \sum_{h \in \Zd_0, |h| \geqs \ro, |k-h| \geqs \ro}
\hspace{-0.2cm}
{|P_{h, k - h}|^2 \over (|h|^{p} |k-h|^{n} + |h|^n |k-h|^{p})^2}
~ \in (0,+\infty)~. \label{deakp} \feq}
The above decomposition follows noting that $\Zd_{0 k}$ is the disjoint union of
the domains of the sums defining $\Km_{p n}(k)$ and $\DK_{p n}(k)$.
$\Km_{p n}(k)$ is finite, involving finitely many summands; $\DK_{p n}(k)$ is
finite as well, since we know that $\KK_{p n}(k) < +\infty$. \parn
\textsl{Step 2. For each $k \in \Zd_0$, one has the representation \rref{kmnd}
$$ \Km_{p n}(k) =
4 |k|^{2 p} \sum_{h \in \Zd_{0 k}, |h | < \ro }\hspace{-0.3cm}
{| P_{h, k - h}|^2 + \teta(|k-h| - \ro) | P_{k-h, h}|^2
\over (|h|^p |k-h|^{n} + |h|^{n} |k-h|^{p})^2}~.$$
If $|k| \geqs 2 \ro$, one can replace $\Zd_{0 k}$ with  $\Zd_{0}$
and $\teta(|k-h| - \ro)$ with $1$}. \parn
To prove \rref{kmnd} we reexpress the sum in Eq.\,\rref{deakm}, using
Eq.\,\rref{tesif} with $f(h) \equiv f_k(h)$ \hbox{$:=
|P_{h, k - h}|^2$ $/(|h|^{p} |k-h|^{n} + |h|^n |k-h|^{p})^2$}.
To go on, assume
$|k| \geqs 2 \ro$; then, for all $h \in \Zd_0$ with $|h| < \ro$ one has
$|k - h| \geqs |k| - |h| > \ro$ whence
$h \neq k$ (i.e., $h \in \Zd_{0 k}$) and $H(|k-h| - \ro) = 1$, two facts which
justify the replacements indicated above.
\parn
\textsl{Step 3. For each $k \in \Zd_0$ one has
\beq 0 <  \DK_{p n}(k) \leqs \dK_{p n}~, \label{boundk} \feq
with $\dK_{p n}$ as in Eq.\,\rref{dedek}.}
The obvious relation $0 < \DK_{p n}(k)$ has been already noted; in the sequel we prove that
$\DK_{p n}(k) \leqs \dK_{p n}$.
The definition \rref{deakp} of $\DK_{p n}(k)$ contains the term
\parn
\vbox{
\beq
{|k|^{2 p} |P_{h, k - h}|^2 \over (|h|^{p} |k-h|^{n} + |h|^n |k-h|^{p})^2}
\leqs  {|k|^{2 p} \sin^2 \te_{h, k - h} \over (|h|^{p} |k-h|^{n} + |h|^n |k-h|^{p})^2}
\label{onthesesk} \feq
$$ \leqs {B_{p n} \over 8} \left({1 \over |h|^{2 n}} + {1 \over |k-h|^{2 n}} \right)~. $$
}
\noindent
The first and the second inequality \rref{onthesesk} follow, respectively, from the bound \rref{bou}
$|P_{h \ell}| \leqs \sin \te_{h \ell}$ and from \rref{onthesik}, in both cases with $\ell=k-h$.
Inserting \rref{onthesesk} into \rref{deakp}, we obtain
\beq \DK_{p n}(k) \leqs
{B_{p n} \over 2} \, \Big( \sum_{h \in \Zd_0, |h| \geqs \ro, |k-h| \geqs \ro} {1 \over |h|^{2 n}} +
\sum_{h \in \Zd_0, |h| \geqs \ro, |k-h| \geqs \ro} {1 \over |k-h|^{2 n}} \Big)~. \feq
The domain of the above two sums is contained in each one of the sets
$\{ h \in \Zd~|~|h | \geqs \ro \}$ and $\{ h \in \Zd~|~|k - h | \geqs \ro \}$; so,
$$ \DK_{p n}(k) \leqs
{B_{p n} \over 2} \, \Big(\sum_{h \in \Zd_0, |h| \geqs \ro} {1 \over |h|^{2 n}} +
\sum_{h \in \Zd_0, |k-h | \geqs \ro} {1 \over |k-h|^{2 n}} \Big)~. $$
Now a change of variable $h \mapsto k - h$ in the second sum shows that the latter is equal to
the former, so
\beq \DK_{p n}(k) \leqs B_{p n}  \sum_{h \in \Zd_0, |h | \geqs \ro}
{1 \over |h|^{2 n}}~. \label{eefgk} \feq
Finally, Eq.\,\rref{eefgk} and Eq.\,\rref{desnu} with $\nu = 2 n$ give
$$ \DK_{p n}(k)
\leqs
{2 \pi^{d/2} B_{p n} \over \Gamma(d/2)}
\sum_{i=0}^{d-1} \left( \barray{c} d - 1 \\ i \farray \right) {d^{d/2-1/2-i/2} \over (2 n - 1 -i)
(\ro - 2 \sqrt{d})^{2 n - 1 -i}}~=\dK_{p n}~\mbox{as in \rref{dedek}}~. $$
\parn
\textsl{Step 4. One has the inequalities \rref{dkmnd}
$\Km_{p n}(k) < \KK_{p n}(k) \leqs \Km_{p n}(k) + \dK_{p n}$}.
These relations follow immediately from the decomposition
\rref{decompk}
$\KK_{p n}(k)= \Km_{p n}(k) + \DK_{p n}(k)$
and from the bounds \rref{boundk} on $\DK_{p n}(k)$. \parn
\textsl{Step 5. One has the equalities \rref{claib}
$\Km_{n}(R_r k) = \Km_{n}(k)$, $\Km_{n}(P_\sigma k) = \Km_{n}(k)$,
involving the reflection and permutation operators  $R_r, P_\sigma$.}
The verification is based on considerations very similar
to the ones that follow Eq.\,\rref{claim}.
\parn
\textsl{Step 6.  For all $k \in \Zd_0$ we have the inequality
\rref{kf}}
$$ \Km_{p n}(k) \leqs 8 \sum_{h \in \Zd_0, |h| < \rho} {1 \over |h|^{2 n}}
E_{p n}\left(\cos \te_{h k}\,, {|h| \over |k|} \right)~. $$
To prove this we start from the expression \rref{kmnd} for $\Km_{p n}(k)$;
we substitute therein the obvious inequality $\teta(|k-h| - \rho)
\leqs 1$, and the relations $|P_{h, k - h}|, |P_{k-h,h}| \leqs \sin \te_{h, k -h}$
following from \rref{bou}. This gives
\beq \Km_{p n}(k) \leqs 8 \sum_{h \in \Zd_{0 k}, |h | < \ro }\hspace{-0.3cm}
{|k|^{2 p} \sin^2 \te_{h, k - h}
\over (|h|^{p} |k-h|^{n} + |h|^{n} |k-h|^{p})^2}~.  \feq
On the other hand, the $h$-th term in the above sum equals
$\dd{{1 \over |h|^{2 n}}
E_{p n}\left(\cos \te_{h k}, {|h| \over |k|} \right)}$ due to
Eq.\,\rref{funk}, so we have the thesis \rref{kf}.
\parn
\textsl{Step 7. Let $k \in \Zd$, $|k| \geqs \M \rho$;
we have the inequalities \rref{gwk}
$$ \Km_{p n}(k) \leqs 8 \sum_{j=0}^{m} {\QQ_{p n j}(\vk) \over |k|^{\gamma_{p n j}}}  +
8 \,{ V_{p n} Y_{p n} \over |k|^{\gamma_{p n, m+1}}} \leqs \Kk_{p n}~, $$
where all the objects in the right hand side are defined as indicated in
item (iii).}
In order to prove this we start from the inequality \rref{kf};
for each $h \in \Zd_0$ with $|h| < \rho$, on account of Eq.\,\rref{onehask} the general term of the sum in
\rref{kf} fulfills the following:
\beq {1 \over |h|^{2 n}}
E_{p n}\left(\cos \te_{h k}, {|h| \over |k|} \right) \label{ineqqk} \feq
$$ = {1 \over |h|^{2 n}} \left( \sum_{j=0}^m Q_{p n j}(\cos \te_{h k})
{\left(|h| \over |k|\right)}^{\gamma_{p n j}} + S_{p n, m+1}(\cos \te_{h k}\,, {|h| \over |k|})
{\left(|h| \over |k|\right)}^{\gamma_{p n, m+1}} \right) $$
$$ \leqs {1 \over |h|^{2 n}} \left( \sum_{j=0}^m Q_{p n j}(\cos \te_{h k})
{\left(|h| \over |k|\right)}^{\gamma_{p n j}} + V_{p n}
{\left(|h| \over |k|\right)}^{\gamma_{p n, m+1}} \right)~.
$$
The last inequality depends on the remark that $|h|/|k| < 1/\M$
and from the definition \rref{wpnmk} of $V_{p n}$. The relation
\rref{kf}, the inequality in \rref{ineqqk} and the definitions
\rref{qqpnj} of $\QQ_{p n j}$, \rref{ypnmk} of $Y_{p n}$ give the first inequality \rref{gwk};
the second inequality \rref{gwk} is an obvious consequence
of the first one and of the definition \rref{kpnm} of $\Kk_{p n}$.
\parn
\textsl{Step 8. The previous results imply the inequality \rref{impifik}.}
This statement is obvious.~ \fine
\vfill \eject \noindent
{~}
\vskip -1.5cm \noindent
\section{The function $\boma{\GG_{p n}}$}
\label{sectg}
As in Section \ref{seckato}, we consider $p, n \in \reali$ such that $p \geqs n > d/2 +1$.
For $k \in \Zd_0$, we recall the definition \rref{ggnd}
$$ \GG_{p n}(k) :=
4 \sum_{h \in \Zd_{0 k}}
{(|k|^p - |k - h|^p)^2 | P_{h, k - h}|^2  \over (|h|^p |k-h|^{n-1} + |h|^n |k-h|^{p-1})^2}~. $$
The forthcoming proposition presents estimates
about $\GG_{p n}(k)$ and its $\sup$ for $k \in \Zd_0$;
its structure is very similar to the one of Proposition \ref{prokknd}
about $\KK_{p n}$ and its sup.
The proof is given only for completeness, since it
is just a rephrasing of arguments employed in Section \ref{sectk}.
Let us also remark that, for $p=n$, the results presented hereafter become
very similar to those appearing in Proposition B.1 of
\cite{cog}.
\begin{prop}
\label{proggnd}
\textbf{Proposition.} As in Eqs.\,\rref{cutoffk}-\rref{order},
let us choose a cutoff, a factor and an order
$$ \ro \in (2 \sqrt{d},+\infty)~, \qquad \M \in (1,+\infty)~, \qquad m \in \naturali~; $$
then the following holds
(with the functions and quantities
$\Gg_{p n}$, $\dG_{p n}$, ..., $Z_{p n}$, $\Gk_{p n}$ mentioned in the sequel
depending parametrically on $\ro$, $\M$, $m$ and $d$: $\Gg_{p n}(k) \equiv \Gg_{p n \ro d}(k)$,
$\dG_{p n} \equiv \dG_{p n \ro d}$, ..., $Z_{p n} \equiv Z_{p n \rho m d},
\Gk_{p n} \equiv \Gk_{p n \ro \M m d}$). \parn
(i) The function $\GG_{p n}$ fulfills the inequalities
\beq \Gg_{p n}(k) < \GG_{p n}(k) \leqs \Gg_{p n}(k) + \dG_{p n}
~~\mbox{for all $k \in \Zd_0$}~. \label{dgmnd} \feq
Here
\beq \Gg_{p n}(k) :=
4 \sum_{h \in \Zd_{0 k}, |h| < \ro \op |k-h| < \ro}
{(|k|^p - |k - h|^p)^2 | P_{h, k - h}|^2  \over (|h|^p |k-h|^{n-1} + |h|^n |k-h|^{p-1})^2}~; \label{deagm} \feq
this function can be reexpressed as
\beq \Gg_{p n}(k) = 4 \sum_{h \in \Zd_{0 k}, |h | < \ro }
\Big(
{(|k|^p - |k - h|^p)^2 | P_{h, k - h}|^2 \over (|h|^p |k-h|^{n-1} + |h|^n |k-h|^{p-1})^2}
\label{gmnd} \feq
$$ +  {\teta(|k-h| - \ro) (|k|^p - |h|^p)^2 | P_{k-h, h}|^2
\over (|k-h|^p |h|^{n-1} + |k-h|^n |h|^{p-1})^2} \Big) $$
(recall that $| P_{k-h, h}| = | P_{h,k-h}|$ if $d \geqs 3$, due to \rref{recall}).
If $|k| \geqs 2 \ro$, in Eq.\,\rref{gmnd} one can replace $\Zd_{0 k}$ with $\Zd_0$
and $\teta(|k-h| - \ro)$ with $1$.
Moreover
\beq \dG_{p n} :=
{2 \pi^{d/2} C_{p n} \over \Gamma(d/2)}
\sum_{i=0}^{d-1} \Big( \barray{c} d - 1 \\ i \farray \Big) {d^{d/2-1/2-i/2} \over (2 n - 3 - i)
(\ro - 2 \sqrt{d})^{2 n - 3 - i}}~, \label{dedeg} \feq
with $C_{p n}$ as in \rref{decin}. \parn
(ii) Consider the reflection and permutation operators
$R_r$, $P_{\si}$ defined by \rref{refper}. Then
\beq \Gg_{p n}(R_r k) = \Gg_{p n}(k)~, \quad \Gg_{p n}(P_\sigma k) = \Gg_{p n}(k)~
\qquad \mbox{for each $k \in \Zd_0$.} \label{claig} \feq
(iii) Denoting with $F_{p n}$ the function in Lemma \ref{deden}, one has
\beq \Gg_{p n}(k) \leqs 4 \sum_{h \in \Zd_0, |h| < \rho} {1 \over |h|^{2 n - 2}}
\, F_{p n}\left(\cos \te_{h k}\,, {|h| \over |k|} \right) \qquad \mbox{for all $k \in \Zd_0$.} \label{gf} \feq
Now, consider the sequence of exponents
$0 = \lambda_{p n 0} < \lambda_{p n 1} < \lambda_{p n 2} ...$ and the sequences of
polynomials $P_{p n j}$ ($j \in \naturali)$
and functions $R_{p n j}$ ($j \in \naturali \setminus \{ 0 \}$) involved
in the expansion of $F_{p n}$ according to Lemma \ref{deden2}; then
\beq \Gg_{p n}(k) \leqs 4 \sum_{j=0}^{m} {\PP_{p n j}(\vk) \over |k|^{\lambda_{p n j}}}  +
4 \, {W_{p n} Z_{p n} \over |k|^{\lambda_{p n, m+1}}} \leqs \Gk_{p n}
\qquad \mbox{for $k \in \Zd,~ |k| \geqs \M \rho$}~. \label{gwg} \feq
Here $\PP_{p n j}$ are the functions defined as follows on the spherical hypersurface $\Sd$:
\beq \PP_{p n j} : \Sd \vain \reali~, \qquad u \mapsto \PP_{p n j}(u) :=
\sum_{h \in \Zd_0, |h| < \rho} {P_{p n j}(\cos \te_{h u}) \over |h|^{2 n - 2 - \lambda_{p n j}}}
\label{pppnj} \feq
(these are polynomials in the components of $u$, which can be computed
using Lemma \ref{lemcos4}). Moreover
\beq W_{p n} := \max_{c \in [-1,1], \xi \in [0,1/\M]} R_{p n, m+1}(c, \xi), \label{wpnm} \feq
\beq Z_{p n} := \hspace{-0.4cm} \sum_{h \in \Zd_0, |h| < \rho} {1 \over |h|^{2 n - 2 - \lambda_{p n, m+1}}}~.
\label{zpn} \feq
Finally,
\beq \Gk_{p n} := 4 \, \max_{u \in \Sd, \, \epsilon \in [0,1/(\M \rho)]}
\left( \sum_{j=0}^{m} \epsilon^{\lambda_{p n j}} \PP_{p n j}(u)
 + \epsilon^{\lambda_{p n, m+1}} W_{p n} Z_{p n} \right)~. \label{gpnm} \feq
(iv) Items (i) and (iii) imply
\beq {~} \hskip -0.4cm \sup_{k \in \Zd_0} \GG_{p n}(k)
\leqs \max \left( \max_{k \in \Zd_0, |k| < \M \rho} \Gg_{p n}(k), \Gk_{p n}
\right) + \dG_{p n}~. \label{impifi} \feq
\end{prop}
The proof of the above statements will be given after the following comment.
\begin{rema}
\label{remgn}
\textbf{Remark} (very similar to Remark \ref{remkn}). From Theorem \ref{maing} we know that
the sharp constant $G_{p n}$ of \rref{katineqa} fulfills
$G_{p n} \leqs \GP_{p n} := (2 \pi)^{-d/2} \sqrt{\sup_{k \in \Zd_0} \GG_{p n}(k)}$.
So, using for $\sup_{k \in \Zd_0} \GG_{p n}(k)$ the bound \rref{impifi},
we conclude
\beq G_{p n} \leqs \GP_{p n} \leqs \GPP_{p n} := {1 \over (2 \pi)^{d/2}}
\sqrt{\max \left( \max_{k \in \Zd_0, |k| < \M \rho} \Gg_{p n}(k), \Gk_{p n}
\right) + \dG_{p n}}~~. \label{gpnbou} \feq
The bound $\GPP_{p n}$ is suitable for computer implementation,
an aspect treated in Section \ref{sec345}.
From the statements in Proposition \ref{proggnd},
it is clear that $\GPP_{p n}$ gives a good approximation of
$\GP_{p n}$ under two conditions: \parn
(a) $\Gk_{p n}$ must be smaller or not too larger than
$\max_{|k| < \M \rho} \Gg_{p n}(k)$. \parn
(b) $\dG_{p n}$ must be small in comparison with $\max \left(\max_{|k| < \M \rho} \Gg_{p n}(k), \Gk_{p n} \right)$.
\end{rema}
\textbf{Proof of Proposition \ref{proggnd}.} We fix $\ro, \M, m$ as in Eqs.\,\rref{cutoffk}-\rref{order}.
Our argument is divided in several steps; more precisely,
Steps 1-4 prove the statements in (i) while
Steps 5, 6-7 and 8 are about the statements in (ii),(iii) and (iv), respectively.
The assumption \rref{cutoffk} $\ro > 2 \sqrt{d}$ is essential in Step 3.
\parn
\textsl{Step 1. One has
\beq \GG_{p n}(k)= \Gg_{p n}(k) + \DG_{p n}(k)~\quad \mbox{for all $k \in \Zd_0$}~, \label{decompg} \feq
where, as in \rref{deagm},
$\Gg_{p n}(k) :=\dd{
4 \sum_{h \in \Zd_{0 k}, |h| < \ro \op |k-h| < \ro}
{(|k|^p - |k - h|^p)^2 | P_{h, k - h}|^2  \over (|h|^p |k-h|^{n-1} + |h|^n |k-h|^{p-1})^2}}~
$,
while
\beq \DG_{p n}(k) := 4 \sum_{h \in \Zd_0, |h| \geqs \ro, |k-h| \geqs \ro}
{(|k|^p - |k - h|^p)^2 | P_{h, k - h}|^2  \over (|h|^p |k-h|^{n-1} + |h|^n |k-h|^{p-1})^2}
~ \in (0,+\infty)~. \label{deagp} \feq}
The above decomposition follows noting that $\Zd_{0 k}$ is the disjoint union of
the domains of the sums defining $\Gg_{p n}(k)$ and $\DG_{p n}(k)$.
$\Gg_{p n}(k)$ is finite, involving finitely many summands; $\DG_{p n}(k)$ is
finite as well, since we know that $\GG_{p n}(k) < +\infty$. \parn
\textsl{Step 2. For each $k \in \Zd_0$, one has the representation \rref{gmnd}
$$ \Gg_{p n}(k)
= 4 \sum_{h \in \Zd_{0 k}, |h | < \ro } \Big(
{(|k|^p - |k - h|^p)^2 | P_{h, k - h}|^2 \over (|h|^p |k-h|^{n-1} + |h|^n |k-h|^{p-1})^2} $$
$$ +  {\teta(|k-h| - \ro) (|k|^p - |h|^p)^2 | P_{k-h, h}|^2
\over (|k-h|^p |h|^{n-1} + |k-h|^n |h|^{p-1})^2} \Big). $$
If $|k| \geqs 2 \ro$, one can replace $\Zd_{0 k}$ with  $\Zd_{0}$
and $H(|k-h| - \ro)$ with $1$}.
\parn
To prove \rref{gmnd} we reexpress the sum in Eq.\,\rref{deagm}, using
Eq.\,\rref{tesif} with $f(h) \equiv f_k(h)$
\hbox{$:=(|k|^p - |k - h|^p)^2 $ $| P_{h, k - h}|^2 / (|h|^p |k-h|^{n-1} + |h|^n |k-h|^{p-1})^2$}.
To go on, assume
$|k| \geqs 2 \ro$; then, for all $h \in \Zd_0$ with $|h| < \ro$ one has
$|k - h| \geqs |k| - |h| > \ro$ whence
$h \neq k$ (i.e., $h \in \Zd_{0 k}$) and $H(|k-h| - \ro) = 1$, two facts which
justify the replacements indicated above.
\parn
\textsl{Step 3. For each $k \in \Zd_0$ one has
\beq 0 <  \DG_{p n}(k) \leqs \dG_{p n}~, \label{boundg} \feq
with $\dG_{p n}$ as in Eq.\,\rref{dedeg}.}
The obvious relation $0 < \DG_{p n}(k)$ was already noted; in the sequel we prove that
$\DG_{p n}(k) \leqs \dG_{p n}$.
The definition \rref{deagp} of $\DG_{p n}(k)$ contains the term
\parn
\vbox{
\beq
{(|k|^p - |k - h|^p)^2 | P_{h, k - h}|^2 \over (|h|^p |k-h|^{n-1} + |h|^n |k-h|^{p-1})^2}
\leqs  {(|k|^p - |k - h|^p)^2 \sin^2 \te_{h, k-h} \over (|h|^{p} |k-h|^{n-1} + |h|^{n} |k-h|^{p-1})^2}
\label{ontheses} \feq
$$ \leqs {C_{p n} \over 8} \left({1 \over |h|^{2 n - 2}} + {1 \over |k-h|^{2 n - 2}} \right) $$
}
\noindent
The first and the second inequality \rref{ontheses} follow, respectively, from the bound \rref{bou}
$|P_{h \ell}| \leqs \sin \te_{h \ell}$ and from \rref{onthesi}, in both cases with $\ell=k-h$.
Inserting \rref{ontheses} into \rref{deagp} we obtain
\beq \DG_{p n}(k) \leqs
{C_{p n} \over 2} \, \Big( \sum_{h \in \Zd_0, |h| \geqs \ro, |k-h| \geqs \ro} {1 \over |h|^{2 n - 2}} +
\sum_{h \in \Zd_0, |h| \geqs \ro, |k-h| \geqs \ro} {1 \over |k-h|^{2 n - 2}} \Big)~. \feq
The domain of the above two sums is contained in each one of the sets
$\{ h \in \Zd~|~|h | \geqs \ro \}$ and $\{ h \in \Zd~|~|k - h | \geqs \ro \}$; so,
$$ \DG_{p n}(k) \leqs
{C_{p n} \over 2} \, \Big(\sum_{h \in \Zd_0, |h| \geqs \ro} {1 \over |h|^{2 n - 2}} +
\sum_{h \in \Zd_0, |k-h | \geqs \ro} {1 \over |k-h|^{2 n - 2}} \Big)~. $$
Now a change of variable $h \mapsto k - h$ in the second sum shows that the latter is equal to
the former, so
\beq \DG_{p n}(k) \leqs C_{p n}  \sum_{h \in \Zd_0, |h | \geqs \ro}
{1 \over |h|^{2 n - 2}}~. \label{eefg} \feq
Finally, Eq.\,\rref{eefg} and Eq.\,\rref{desnu} with $\nu = 2 n - 2$ give
$$ \DG_{p n}(k)
\leqs
{2 \pi^{d/2} C_{p n} \over \Gamma(d/2)}
\sum_{i=0}^{d-1} \left( \barray{c} d - 1 \\ i \farray \right) {d^{d/2-1/2-i/2} \over (2 n - 3 -i)
(\ro - 2 \sqrt{d})^{2 n - 3 -i}}~=\dG_{p n}~\mbox{as in \rref{dedeg}}~. $$
\parn
\textsl{Step 4. One has the inequalities \rref{dgmnd}
$\Gg_{p n}(k) < \GG_{p n}(k) \leqs \Gg_{p n}(k) + \dG_{p n}$}.
These relations follow immediately from the decomposition
\rref{decompg}
$\GG_{p n}(k)= \Gg_{p n}(k) + \DG_{p n}(k)$
and from the bounds \rref{boundg} on $\DG_{p n}(k)$. \parn
\textsl{Step 5. One has the equalities \rref{claig}
$\Gg_{p n}(R_r k) = \Gg_{p n}(k)$, $\Gg_{p n}(P_\sigma k) = \Gg_{p n}(k)$,
involving the reflection and permutation operators  $R_r, P_\sigma$.}
The verification is based on considerations very similar
to the ones that follow Eq.\,\rref{claim}.
\parn
\textsl{Step 6. For all $k \in \Zd_0$ we have the inequality
\rref{gf}}
$$ \Gg_{p n}(k) \leqs 4 \sum_{h \in \Zd_0, |h| < \rho} {1 \over |h|^{2 n - 2}}
F_{p n}\left(\cos \te_{h k}, {|h| \over |k|} \right)~. $$
To prove this we start from
the expression \rref{gmnd} for $\Gg_{p n}(k)$; we substitute therein
the obvious inequality $\teta(|k-h| - \rho) \leqs 1$,
and the relations $|P_{k - h,h}|, |P_{h, k - h}| \leqs \sin \te_{h, k-h}$
following from \rref{bou}.
This gives
\beq \Gg_n(k) \leqs 4 \sum_{h \in \Zd_{0 k}, |h | < \ro }
\sin^2 \te_{h, k-h} \left({(|k|^p - |k - h|^p)^2 \over (|h|^p |k-h|^{n-1} + |h|^n |k-h|^{p-1})^2} \right.
\feq
$$ +  \left. {(|k|^p - |h|^p)^2 \over (|k-h|^p |h|^{n-1} + |k-h|^n |h|^{p-1})^2} \right); $$
on the other hand, the $h$-th term in the above sum equals
$\dd{{1 \over |h|^{2 n - 2}}
F_{p n}\left(\cos \te_{h k}, {|h| \over |k|} \right)}$ due to
Eq.\,\rref{fung}, so we have the thesis \rref{gf}.
\parn
\textsl{Step 7. Let $k \in \Zd$, $|k| \geqs \M \rho$;
we have the inequalities \rref{gwg}
$$ \Gg_{p n}(k) \leqs 4 \sum_{j=0}^{m} {\PP_{p n j}(\vk) \over |k|^{\lambda_{p n j}}}  +
4 \, {W_{p n} Z_{p n} \over |k|^{\lambda_{p n, m+1}}} \leqs \Gk_{p n} ~,$$
where all the objects in the right hand side are defined as indicated in
item (iii).}
In order to prove this we start from the inequality \rref{gf};
for each $h \in \Zd_0$ with $|h| < \rho$, the general term of the sum in
\rref{gf} fulfills, due to \rref{onehas}:
\beq {1 \over |h|^{2 n - 2}}
F_{p n}\left(\cos \te_{h k}\,, {|h| \over |k|} \right) \label{ineqq} \feq
$$ = {1 \over |h|^{2 n - 2}} \left( \sum_{j=0}^m P_{p n j}(\cos \te_{h k})
{\left(|h| \over |k|\right)}^{\lambda_{p n j}} + R_{p n, m+1}(\cos \te_{h k}\,, {|h| \over |k|})
{\left(|h| \over |k|\right)}^{\lambda_{p n, m+1}} \right) $$
$$ \leqs {1 \over |h|^{2 n - 2}} \left( \sum_{j=0}^m P_{p n j}(\cos \te_{h k})
{\left(|h| \over |k|\right)}^{\lambda_{p n j}} + W_{p n}
{\left(|h| \over |k|\right)}^{\lambda_{p n, m+1}} \right)~.
$$
The last inequality depends on the remark that $|h|/|k| < 1/\M$
and from the definition \rref{wpnm} of $W_{p n}$.
The relation
\rref{gf}, the inequality in \rref{ineqq} and the definitions
\rref{pppnj} of $\PP_{p n j}$, \rref{zpn} of $Z_{p n}$ give the first inequality \rref{gwg};
the second inequality \rref{gwg} is an obvious consequence
of the first one and of the definition \rref{gpnm} of $\Gk_{p n}$.
\parn
\textsl{Step 8. The previous results imply the inequality \rref{impifi}.}
This statement is obvious. \par
\fine
\section{Computation of $\boma{\KPP_{p n}}$, $\boma{\GPP_{p n}}$
in the cases \rref{casesdue} or \rref{casestre}}
\label{sec345}
Throughout this section the space dimension is
\beq d = 3~. \feq
We make extensive reference to Propositions \ref{prokknd} and \ref{proggnd}
and to the definitions of the upper bounds  $\KPP_{p n}$, $\GPP_{p n}$
in Eqs.\,\rref{kpnbou}\,\rref{gpnbou}; our purpose is to describe
the computer implementation of these definitions and,
in particular, the calculations performed for
the cases \rref{casesdue} or \rref{casestre}.
\salto
\textbf{Computation of $\boma{\KPP_{p n}}$.} Eq.\,\rref{kpnbou} reads
\beq \KPP_{p n} = {1 \over (2 \pi)^{3/2}}
\sqrt{\max \left( \max_{k \in \Zt_0, |k| < \M \rho} \Km_{p n}(k), \Kk_{p n}
\right) + \dK_{p n}}~~. \label{eqqk} \feq
We recall that the right hand side of the above equation depends on
a cutoff $\rho > 2 \sqrt{3}$, on a factor $\mu > 1$ and on the order $m \in \naturali$
involved in the expansion determining $\Kk_{p n}$. \par
Hereafter we describe the choices
of $\rho, \mu, m$ and give details on the calculation
of $\Km_{p n}(k)$, $\Kk_{p n}$ and $\dK_{p n}$, for
the cases \rref{casesdue} \rref{casestre}.
The results obtained in these cases are summarized
in Table \TC; here the last column contains the final values of $\KPP_{p n}$ given
by the overall procedure, which coincide with the ones anticipated in Table \TA.
Any one of the necessary computations has been performed on a PC with an 8 Gb RAM,
with the software utilities mentioned below.
\salto
\textsl{Choosing $\rho,\mu,m$.} As explained in Remark \ref{remkn},
in the implementation of Eq.\,\rref{eqqk}
it is convenient to choose the cutoff $\rho$, the factor $\M$ and the order $m$ of the
expansion determining $\Kk_{p n}$ so that
$\Kk_{p n}$ is smaller (or not too larger) than
$\max_{|k| < \M \rho} \Km_{p n}(k)$, and
$\dK_{p n}$ is small with respect to
$\max \big( \max_{|k| < \M \rho} \Km_{p n}(k),$
$\Kk_{p n} \big)$. In all our computations we have taken
\beq \M = 2~, \qquad m = 6 \feq
and we have used cutoffs $\rho$ between $20$ and $100$, chosen
empirically so as to fulfill the previous requirements
(starting from the lower value $\rho=20$, and increasing this
if necessary. Of course, the computational
costs increase with $\rho$).
\salto
\textsl{The error bound $\dK_{p n}$.} This is readily computed
in terms of the cutoff $\rho$ via Eq.\,\rref{dedek}. The constants
$B_{p n}$ in the cited equation are the maxima
of certain functions, see Eqs.\,\rref{debn}\,\rref{debnk}
and Remark \ref{rembpn}.
Our actual computation of the $\dK_{p n}$'s
has been made using {\tt{Mathematica}}. \par
Here and in the sequel, to give an example
of our calculations we consider the case $(p,n) = (5,2)$.
As indicated in Table \TC, we have chosen
\beq \rho=50 \qquad \mbox{for $p=5,~ n=2$}~, \feq
that gives
\beq \delta \Km_{5 2} = 65.0229... \,. \feq
\salto
\textsl{Computation of $\Km_{p n}(k)$ for $k \in \Zt_0$, $|k| < 2 \rho$
and of its maximum over this ball.}
The computation of $\Km_{p n}(k)$ at all
nonzero points in the ball $\{ |k| < 2 \rho \}$
have been performed using an ad hoc {\tt{C}} program,
for each one of the cases \rref{casesdue} \rref{casestre}
These computations allow to determine
the maximum point of $\Km_{p n}$ over the ball; the value of
$\Km_{p n}$ at this point has been subsequently validated using {\tt{Arb}} \cite{arb}, a {\tt{C}}-library that
produces certified estimates on the roundoff errors. For example
in the case $(p,n)=(5,2)$, treated with a cutoff $\rho=50$,
we have found
\beq \max_{k \in \Zt_0, |k| < 100} \Km_{5 2}(k) = \Km_{5 2}(2,1,0) = 263.364...~; \label{k52} \feq
for the other cases, see Table \TC. Computations for all values of $k$ involved in
Eq.\,\rref{k52} and the subsequent determination of the maximum
have required a CPU time of about 3 hours on our PC; the validation
of the value at the maximum point $k=(2,1,0)$ has
been been performed very quickly using {\tt{Arb}}
({\footnote{The {\tt{Arb}} result is
$\Km_{5 2}(2,1,0) =  263.36493191766936106 \pm 9.6212 \times 10^{-14}$.}}).
 The analogous computations
for $n=2$ and each one of the cases $p=7,8,9,10$, based on the larger cutoff $\rho=100$, have required a CPU time of
about 4 days. \par
\salto
\textsl{Computation of $\Kk_{p n}$.} This requires
a rather long procedure, described by item (iii) of
Proposition \ref{prokknd}. First of
all one should consider the function
$E_{p n}(c, \xi)$ of Lemma \ref{deken} and build
the expansion described by Lemma \ref{deken2}.
We have already mentioned the choice $m=6$ for the order
of this expansion; for the cases that we have considered,
where $p,n$ are integers,
the expansion of $E_{p n}(c, \xi)$ involves
only integer powers of $\xi$ and takes the form
\beq E_{p n}(c, \xi) = \sum_{j=0}^6 Q_{p n j}(c) \xi^{j} + S_{p n 7}(c, \xi)
\xi^{7}~. \label{onehask6} \feq
The coefficient $S_{p n 7}(c, \xi)$ in the reminder term
is important, since its maximization (for $c \in [-1,1], \xi \in [0,1/2]$)
gives a constant $V_{p n}$ to be used later (see Eq.\,\rref{wpnmk}). The coefficients $Q_{p n j}(c)$ must be used
to build certain functions $\QQ_{p n j} : \St \vain \reali$, see Eq.\,\rref{qqpnj};
next one should compute a finite zeta-type sum $Y_{p n}$
(see again Eq.\,\rref{wpnmk}). Finally, one determines
$\Kk_{p n}$ via maximization for $u \in \St, \epsilon \in [0,1/2 \rho]$
of a certain function of these variables, built from the previous
ingredients (see Eq.\,\rref{kpnm}). \par
For all values of $p, n$ in \rref{casesdue} \rref{casestre}
the coefficients $Q_{p n j}$, the related functions $\QQ_{p n j}:
\St \vain \reali$ ($j=0,...,6$) and the sums $Y_{p n}$ have been
computed symbolically using {\tt{Mathematica}}.
The maxima (for $c \in [-1,1], \xi \in [0,1/2]$
or for $u \in \St, \ep \in [0, 1/2 \rho]$) defining
$V_{p n}$ and $\Kk_{p n}$ have been computed
numerically, using the internal routines of {\tt{Mathematica}}.
As an example, for $(p,n)=(5,2)$ we have found
\parn
\vbox{
\beq Q_{5 2 0}(c) = 1 - c^2~,~~Q_{5 2 1}(c) = 12 (c - c^3)~,~~Q_{5 2 2}(c) =
-6 + 90 c^2 - 84 c^4~,... \feq
$$ Q_{5 2 6}(c) = -53 + 255 c + 3077 c^2 - 1870 c^3 - 23184 c^4 + 1615 c^5 + 49728 c^6 -
29568 c^8; $$
}
\noindent
\beq V_{5 2} = 2211.24...~.\feq
The other ingredients in these calculations depend
on the cutoff, so we recall the choice $\rho=50$
for the present case.
We have recalled that one can compute from the coefficients $Q_{5 2 j}$ the
functions $\QQ_{5 2 j}$; we report only one of the functions, namely
$$ \QQ_{5 2 2}(u) = 14861.4... - 10448.7...({u_1}^4 + {u_2}^4 + {u_3}^4)
-20668.7...({u_1}^2 {u_2}^2 + {u_1}^2 {u_3}^2 + {u_2}^2 {u_3}^2) $$
\beq \mbox{for $u = (u_1, u_2, u_3) \in \St$}~. \feq
Moreover,
\beq Y_{5 2} = 3.26693...\times 10^{10}~. \feq
Finally, maximizing for $u \in \St$ and $\epsilon \in [0,1/2 \rho] = [0, 1/100]$
the function in Eq.\,\rref{kpnm} we obtain
\beq \Kk_{5 2} = 92.5195... ~,\feq
as indicated in Table \TC.
\salto
\textsl{The final step: determination of $\KPP_{p n}$.}
After computing $\dK_{p n}$,  $\max_{k \in \Zt_0, |k| < 2 \rho} \Km_{p n}(k)$
and $\Kk_{p n}$ one returns to Eq.\,\rref{eqqk} and obtains
$\KPP_{p n}$.
The values indicated with $\KPP_{p n}$ in Table \TC are in fact
the roundups to 3 digits of
the right hand side of Eq.\,\rref{eqqk}. These roundups
are our final upper bounds for the sharp constants $K_{p n}$ in the inequality
\rref{basineqa}, in the cases under consideration.
\spazio
\vbox{
\noindent
\textbf{Table \TC. Details on the computation of $\boma{\KPP_{p n}}$
in the cases \rref{casesdue} \rref{casestre}
(using Proposition \ref{prokknd} with $\boma{\M=2}$, $\boma{m=6}$)}
\vskip 0.2cm \noindent
\vbox{
$$ \begin{tabular}{|c|| c|c|c|c|c||c|}
\hline
$(p,n)$ & ~~~$\ro$~~~~  & $\max \Km_{p n}$ (*) & $k_{max}$ (**) & $\Kk_{p n}$ & $\dK_{p n}$ & $\KPP_{p n}$  \\[0.1cm]
\hline \hline
$(2,2)$ & 20  & 22.0223... & (9,9,9) & 21.6447... & 5.68568... & 0.335 \\ \hline
$(3,2)$ & 20  & 77.8597... & (25, 23, 21) & 84.8166...  & 17.6648... & 0.643   \\ \hline
$(4,2)$ & 20 & 113.227... & (2,2,2) & 89.5797...  & 57.7725...  & 0.831   \\ \hline
$(5,2)$ & 50 & 263.364... & (2,1,0) & 92.5195... & 65.0229... & 1.16 \\ \hline
$(6,2)$ & 50 & 702.295... & (2,1,0) & 97.2697... & 225.357... &  1.94 \\ \hline
$(7,2)$ & 100 & 1884.65... & (2,1,0) & 96.4078... & 376.103...  & 3.02    \\ \hline
$(8,2)$ & 100 & 5018.97... & (2,1,0) & 101.904... & 1345.89... & 5.07   \\ \hline
$(9,2)$ & 100 & 13205.2... & (2,1,0) & 110.191,... & 4870.46... & 8.54    \\ \hline
$(10,2)$ & 100 & 34334.1... & (2,1,0) & 122.389... & 17786.6... & 14.5   \\ \hline
$(3,3)$ & 20 & 25.3013... & (2,1,1) & 11.2784...  & 0.0226087...  & 0.320 \\ \hline
$(4,3)$ & 20 & 71.8198... & (2,1,0) & 44.8074... & 0.0739415...  & 0.539  \\ \hline
$(5,3)$ & 20 & 204.342... & (2,1,0) & 45.4808...  & 0.250165... & 0.909    \\ \hline
$(6,3)$ & 20 & 581.166... & (2,1,0) & 45.9450... & 0.867027... & 1.54    \\ \hline
$(7,3)$ & 20 & 1636.38... & (2,1,0) & 46.3859... & 3.05959... & 2.58    \\ \hline
$(8,3)$ & 20 & 4521.94... & (2,1,0) & 46.9192... & 10.9488... & 4.28    \\ \hline
$(9,3)$ & 20 & 12237.3... & (2,1,0) & 47.5671... & 39.6211... & 7.04    \\ \hline
$(10,3)$ & 20 & 32495.3... & (2,1,0) & 48.3602... & 144.694... & 11.5  \\ \hline
\end{tabular} $$
\vskip 0.2cm \noindent
{\footnotesize{
(*) $\max \Km_{p n}$ stands for $\max_{k \in \Zt_0, |k| < 2 \ro} \Km_{p n}(k)$. \parn
(**) $k_{max}$ is a maximum point of the function $\Km_{p n}$
on the set $\{k \in \Zt_0~|~ |k| < 2 \ro\}$.}}
}}
\vskip 0.5cm \noindent
\textbf{Computation of $\boma{\GPP_{p n}}$.} We follow a general scheme very similar
to the one employed to determine $\KPP_{p n}$; however, in this case we refer
to Proposition \ref{proggnd} and to
Eq.\,\rref{gpnbou}; this reads
\beq \GPP_{p n} = {1 \over (2 \pi)^{3/2}}
\sqrt{\max \left( \max_{k \in \Zt_0, |k| < \M \rho} \Gg_{p n}(k), \Gk_{p n}
\right) + \dG_{p n}}~~, \label{eqqg} \feq
with the right hand side depending on
a cutoff $\rho > 2 \sqrt{3}$, on a factor $\mu > 1$
and on the order $m \in \naturali$ in the expansion
giving $\Gk_{p n}$. \par
Hereafter we describe the choices
of $\rho, \mu, m$ and give details on the calculation
of $\Gg_{p n}(k)$, $\Gk_{p n}$ and $\dG_{p n}$
in the cases \rref{casestre}. The results are summarized
in Table \TD; here the last column contains the final values of $\GPP_{p n}$
anticipated in Table \TB.
For the necessary computations, we have used the hardware and software utilities
already mentioned in relation to the constants $K_{p n}$.
\vfill \eject \noindent
{~}
\vskip -1cm \noindent
\textsl{Choosing $\rho,\mu,m$.} According to Remark \ref{remgn},
in the implementation of Eq.\,\rref{eqqg}
it is convenient to choose the cutoff $\rho$, the factor $\M$ and the order $m$ of the
expansion determining $\Kk_{p n}$ so that
$\Gk_{p n}$ be smaller (or not too larger) than
$\max_{|k| < \M \rho} \Gg_{p n}(k)$, and
$\dG_{p n}$ be small with respect to
$\max \big( \max_{|k| < \M \rho} \Gg_{p n}(k),$
$\Gk_{p n} \big)$. In all our computations we have taken
\beq \M = 2~, \qquad m = 6 \feq
and we have used cutoffs $\rho$ between $20$ and $100$, chosen
empirically so as to fulfill the previous requirements.
\salto
\textsl{The error bound $\dG_{p n}$.} This is readily computed
in terms of the cutoff $\rho$ via Eq.\,\rref{dedeg},
which has been implemented via {\tt{Mathematica}}. The constants
$C_{p n}$ in the cited equation are the maxima
of certain functions, see Eqs.\,\rref{decn}\,\rref{decin};
for the values of $(p,n)$ in \rref{casestre},
these constants have been determined by numerical maximization
(see Eq.\,\rref{valc2} for some examples). \par
\salto
\textsl{Computation of $\Gg_{p n}(k)$ for $k \in \Zt_0$, $|k| < 2 \rho$
and of its maximum over this ball.} The necessary computations
have been performed using a {\tt{C}} program; for each
case $(p,n)$ in \rref{casestre} involving the largest cutoff
$\rho = 100$, the CPU time has been of 4 days approximately.
For all cases in the table, the value of $\Gg_{p m}$ at the maximum point
in the ball $\{ |k| < 2 \rho \}$ has been validated using {\tt{Arb}}.
\salto
\textsl{Computation of $\Gk_{p n}$.} This is based
on the procedure in item (iii) of
Proposition \ref{proggnd}. First of
all one should consider the function
$F_{p n}(c, \xi)$ of Lemma \ref{deden} and build
the expansion described by Lemma \ref{deden2}.
We have already mentioned the choice $m=6$ for the order
of this expansion; for the cases that we have considered,
where $p,n$ are integers,
the expansion of $F_{p n}(c, \xi)$ involves
only integer powers of $\xi$ and takes the form
\beq F_{p n}(c, \xi) = \sum_{j=0}^6 P_{p n j}(c) \,\xi^{j} + R_{p n 7}(c, \xi)
\,\xi^{7}~. \label{onehasg6} \feq
The coefficient $R_{p n 7}$ in the reminder of this expansion determines, after a maximization for
$c \in [-1,1], \xi \in [0,1/2]$, a constant $W_{p n}$ to be used later (see Eq.\,\rref{wpnm}).
The coefficients $P_{p n j}$ are used
to build certain functions $\PP_{p n j} : \St \vain \reali$
(see Eq.\,\rref{pppnj}).
After computing a finite zeta-type sum $Z_{p n}$
(see again Eq.\,\rref{wpnm}), one finally determines
$\Gk_{p n}$ via maximization for $u \in \St, \epsilon \in [0,1/2 \rho] $
of a certain function of these variables, built from the previous
ingredients (see Eq.\,\rref{gpnm}). For the values
of $(p,n)$ in \rref{casestre}, all the above mentioned computations
have been performed using {\tt{Mathematica}}.
\salto
\textsl{The final step: determination of $\GPP_{p n}$.}
After computing $\dG_{p n}$,  $\max_{k \in \Zt_0, |k| < 2 \rho} \Gg_{p n}(k)$
and $\Gk_{p n}$ one returns to Eq.\,\rref{eqqg} and obtains
$\GPP_{p n}$.
The values indicated in Table \TD with $\GPP_{p n}$ are in fact
the roundups to 3 digits of the right hand side of Eq.\,\rref{eqqg}. These roundups
are our final upper bounds for the sharp constants $G_{p n}$ in the inequality
\rref{katineqa}, in the cases under consideration.
\salto
\vbox{
\noindent
\textbf{Table \TD. Details on the computation of $\boma{\GPP_{p n}}$
in the cases \rref{casestre} (using Proposition \ref{proggnd} with $\boma{\M=2}$, $\boma{m=6}$)}
\spazio
$$ \begin{tabular}{|c|| c|c|c|c|c||c|}
\hline
$(p,n)$ & ~~~$\ro$~~~~  & $\max \Gg_{p n}$ (*) & $k_{max}$ (**) & $\Gk_{p n}$ &  $\dG_{p n}$ &
$\GPP_{p n}$  \\[0.1cm]
\hline \hline
$(3,3)$ & 20  & 34.9016... & (9,9,9) & 34.4741... & 12.4785... & 0.438    \\ \hline
$(4,3)$ & 20  & 190.684... & (23,23,23) &  206.799... & 51.5254... & 1.03   \\ \hline
$(5,3)$ & 50  & 325.352... & (4,4,4) & 309.674... & 64.3690... & 1.26    \\ \hline
$(6,3)$ & 50  & 816.449... & (2,1,0) & 437.386... & 230.273... & 2.06    \\ \hline
$(7,3)$ & 50  & 2356.09... & (2,1,0) & 593.730... & 817.263... & 3.58    \\ \hline
$(8,3)$ & 100  & 6611.94... & (2,1,0) & 755.564.... & 1380.53... & 5.68    \\ \hline
$(9,3)$ & 100  & 18068.8... & (2,1,0) & 965.898... & 4977.42... & 9.64    \\ \hline
$(10,3)$ & 100  & 48275.0... & (2,1,0) & 1218.84... & 18110.5... & 16.4    \\ \hline
\end{tabular} $$
\vskip 0.2cm \noindent
{\footnotesize{
(*) $\max \Gg_{p n}$ stands for $\max_{k \in \Zt_0, |k| < 2 \ro} \Gg_{p n}(k)$. \parn
(**) $k_{max}$ is a maximum point of the function $\Gg_{p n}$
on the set $\{k \in \Zt_0~|~ |k| < 2 \ro\}$.}}
}
\salto
\textbf{Comparison with previous works.}
We have already mentioned that
the constants $K_{p n}$, $G_{p n}$ have been discussed
in the previous works \cite{cog} \cite{cok} in the special
case $p=n$ (using the notations $K_{n}, G_n$ for
$K_{n n}$ and $G_{n n}$).
The values of $K^{(+)}_{2 2}$ and $G^{(+)}_{3 3}$ in
Tables \TC and \TD agree with the upper bounds on $K_2$ and $G_3$ computed in
the cited works (for $d=3$).
In \cite{cok} it has been shown that (for $d=3$) $K_{3}$
has an upper bound equal to $0.323$;
the result in Table \TC (namely, the upper bound $K^{(+)}_{3 3} = 0.320$) is a slight improvement,
due to the use of balls with a larger value of the radius $\rho$.
\section{Lower bounds for the sharp constant $\boma{K_{p n}}$
of the inequality \rref{basineqa}}
\label{seclowk}
Consider a pair $p,n$ as in \rref{31}.
To obtain the above mentioned lower bounds it is sufficient to use
the tautological inequality
\beq K_{p n} \geqs {2 \, \| \PPP(v, w) \|_p \over
\| v \|_p \| w \|_{n+1} + \| v \|_n \| w \|_{p+1}} \quad
\mbox{for $v \in \HM{p} \setminus \{0 \}$, $w \in \HM{p+1} \setminus \{0 \}$}~,
\label{kneq} \feq
choosing for $v$ and $w$ two suitable
non zero ``trial vector fields''; rather simple choices
yield the results presented hereafter. Throughout
this section we denote with $a,b$ the first two elements in
the canonical basis of $\reali^d$ and, if $d \geqs 3$, we
write $c$ for the third element. Thus
\beq a := (1,0,...,0),~~ b := (0,1,0,...0)~; \qquad c := (0,0,1,0,...,0)~\,
\mbox{if $d \geqs 3$}. \label{basnat} \feq
Let us also recall the notation $e_k$ for the Fourier basis, see Eq.\,\rref{furba}.
\begin{prop}
\label{propkmenoa}
\textbf{Proposition.} For all real $p,n$ with $p \geqs n > d/2$ one has
\beq K_{p n} \geqs K^{[-]}_p~,
\qquad K^{[-]}_p := {U_d \over (2 \pi)^{d/2}} \, 2^{p/2}~,\label{kpmenoa} \feq
\beq
U_d := \left\{ \barray{cc}  1 & \mbox{if $d \geqs 3$,} \\  {1/\sqrt{2}} & \mbox{if $d=2$}.
\farray \right.
\label{kpmeno1a}
\feq
\end{prop}
\textbf{Proof.}
\textsl{Step 1. For $d \geqs 3$ one has}
\beq K_{p n} \geqs {2^{p/2} \over (2 \pi)^{d/2}}~; \label{step1} \feq
\textsl{so, the thesis \rref{kpmenoa} \rref{kpmeno1a} holds if $d \geqs 3$.}
To prove this we set
\beq v := i b (e_{a} - e_{-a}), \quad w := i c (e_{b} - e_{-b})~. \label{44} \feq
The above vector fields have vanishing mean and divergence
(since $k \sc v_k = k \sc w_k=0$ for each $k$) and clearly belong to
$\HM{m}$ for any real $m$, with
\beq \| v \|_m = \| w \|_m = \sqrt{2}~. \label{step1a} \feq
Using Eqs.\,\rref{deflpk}\,\rref{infert} one finds
\beq \PPP(v,w) =  - {i c \over (2 \pi)^{d/2}} \big( e_{a+b} - e_{-a - b} + e_{a-b} - e_{-a + b} \big) \feq
(note that $\LP_{k} c = c$ for $k=\pm (a+b), \pm(a-b)$); this implies
\beq \| \PPP(v,w) \|_p = {2^{p/2 + 1} \over (2 \pi)^{d/2}}~. \label{step1b} \feq
Now, using Eqs.\,\rref{kneq} and \rref{step1a}\,\rref{step1b} one readily infers
Eq.\,\rref{step1}. \parn
\textsl{Step 2. For $d =2$ one has}
\beq K_{p n} \geqs {2^{p/2} \over 2 \sqrt{2} \pi}~; \label{step10} \feq
\textsl{so, the thesis \rref{kpmenoa} \rref{kpmeno1a} holds if $d =2$.}
To prove this we define $v$ as in \rref{44}, and put
\beq w := i a (e_{b} - e_{-b})~.
\label{440} \feq
Again, $v,  w \in \HM{m}$ for any real $m$, with $\| v \|_m$, $ \| w \|_m $ as in \rref{step1a}.
Using Eqs.\,\rref{lpk}\,\rref{infert} one finds
\beq \PPP(v,w) =
- {i  \over 4 \pi} \big( (a-b) (e_{a+b} - e_{-a - b}) +
(a + b) (e_{a-b} - e_{-a +b}) \big) ; \feq
this implies
\beq \| \PPP(v,w) \|_p = {2^{p/2 - 1/2} \over \pi}~. \label{step1b0} \feq
Now, using Eqs.\,\rref{kneq} and \rref{step1a}\,\rref{step1b0} one infers
Eq.\,\rref{step10}. \fine
\begin{rema}
\textbf{Remark.}
For $p=n$ Eq.\,\rref{kpmenoa} becomes
(writing $K_{n}, K^{[-]}_n$ for $K_{n n}$, $K^{[-]}_{n n}$)
\beq K_{n} \geqs K^{[-]}_{n}~, \qquad K^{[-]}_{n} := {U_d \over (2 \pi)^{d/2}} \cdot 2^{n/2}~;
\label{kpmenonn} \feq
this bound has been already proposed in \cite{cok}
({\footnote{See Eq.\,(3.9) of \cite{cok}, corresponding to Eq.\,(3.22) in
the extended arXiv version of the same paper. The cited work states that
$U_d := 1$ for $d \geqs 3$ and $U_d := (2 - \sqrt{2})^{1/2}$
for $d=2$; the value of $U_2$ is wrong, and reflects an error
in calculations with the trial vector fields presented therein.
The correct formulation of the estimate
of \cite{cok} for $d=2$ requires $U_2$ to be defined as in \rref{kpmeno1a}.}}).
\end{rema}
\begin{prop}
\label{propkmeno}
\textbf{Proposition.} Let $p, n$ be real, with $p \geqs n > d/2$.
For each $\r \in \naturali_0 := \{1,2,3,...\}$ one has
\beq K_{p n} \geqs K^{\{-\}}_{p n}(\r)~,
\label{kpmeno} \feq
where
\beq
K^{\{-\}}_{p n}(\r) := \left\{ \barray{ll} \dd{
{\sqrt{2} \over (2 \pi)^{d/2}}\,
{\sqrt{1 + (1 + 4 \r^2)^p} \over
\r^p (1+\r^2)^{n/2 + 1/2} + \r^n (1 + \r^2)^{p/2 + 1/2}} }
 & \mbox{if $d \geqs 3$,} \vspace{0.2cm} \\
\dd{{\sqrt{2} \over 2 \pi} \, {\sqrt{1 + (1 + 2 \r^2)^2(1 + 4 \r^2)^{p-1}} \over
\r^p (1 + \r^2)^{n/2 + 1} + \r^n (1 + \r^2)^{p/2 + 1}} } & \mbox{if $d=2$}.
\farray \right.
\label{kpmeno1}
\feq
Thus
\beq K_{p n} \geqs K^{\{-\}}_{p n} := \sup_{\r \in \naturali_0} K^{\{-\}}_{p n}(\r)~. \label{kpmenoz}
\feq
\end{prop}
\textbf{Proof.} Of course, Eq.\,\rref{kpmenoz} is an obvious consequence
of \rref{kpmeno}; in the sequel we derive Eqs.\,\rref{kpmeno}\,\rref{kpmeno1}
for any $\r \in \naturali_0$, proceeding in two steps. \parn
\textsl{Step 1. For $d \geqs 3$ one has}
\beq K_{p n} \geqs  {\sqrt{2} \over (2 \pi)^{d/2}}\,
{\sqrt{1 + (1 + 4 \r^2)^p} \over
\r^p (1+\r^2)^{n/2 + 1/2} + \r^n (1 + \r^2)^{p/2 + 1/2}}~; \label{step2} \feq
\textsl{so, the thesis \rref{kpmeno} \rref{kpmeno1} holds if $d \geqs 3$.}
To prove this we use Eq.\,\rref{kneq} with
\beq v \equiv v(\ell) := i b(e_{\r a} - e_{- \r a})~, \quad
w \equiv w(\r) := i c (e_{\r a + b} - e_{-\r a-b})~. \label{defvi} \feq
We have $v,w \in\HM{m}$ for any real $m$, with
\beq \| v \|_m = \sqrt{2} \,\r^m~, \quad
\| w \|_m = \sqrt{2} \,(1 + \r^2)^{m/2}~. \label{step2a} \feq
Using Eqs.\,\rref{deflpk}\,\rref{infert} we find
\beq \PPP(v,w) = - {i c \over (2 \pi)^{d/2}} \big( e_{2 \r a+b} - e_{- 2 \r a - b} - e_{b} + e_{-b} \big) \feq
(note that $\LP_{k} c = c$ for $k=\pm (2 \r a+b), \pm b$); this implies
\beq \| \PPP(v,w) \|_p = {\sqrt{2} \over (2 \pi)^{d/2}} \, \sqrt{1 + (1 + 4 \r^2)^p} ~. \label{step2b} \feq
Now, using Eqs.\,\rref{kneq} and \rref{step2a}\,\rref{step2b} we infer
Eq.\,\rref{step2}.
\parn
\textsl{Step 2. For $d =2$ one has}
\beq K_{p n} \geqs {\sqrt{2} \over 2 \pi} \, {\sqrt{1 + (1 + 2 \r^2)^2(1 + 4 \r^2)^{p-1}} \over
\r^p (1 + \r^2)^{n/2 + 1} + \r^n (1 + \r^2)^{p/2 + 1}} ~~; \label{step20} \feq
\textsl{so, the thesis \rref{kpmeno} \rref{kpmeno1} holds if $d=2$.}
To prove this we use Eq.\,\rref{kneq} with $v \equiv v(\r)$ as in \rref{defvi}, and
\beq w \equiv w(\r) := i \cc (e_{\r a + b} - e_{-\r a - b})~, \qquad \cc := {a - \r b \over \sqrt{1 + \r^2}}~. \feq
Again $v,w \in \HM{m}$ for any real $m$, with $\| v \|_m$, $\| w \|_m$ as in \rref{step2a}.
Using Eqs.\,\rref{lpk}\,\rref{infert} we find
\beq \PPP(v,w) =
- {i \over 2 \pi} {(1 + 2 \r^2) (a - 2 \r b)
\over \sqrt{1 + \r^2} (1 + 4 \r^2)} (e_{2 \r a + b} - e_{-2 \r a-b}) +
{i \over 2 \pi} {a \over \sqrt{1 + \r^2}} (e_{b}- e_{-b})~;
\feq
this implies
\beq \| \PPP(v,w) \|_p = {\sqrt{2} \over 2 \pi} \, {\sqrt{1 + (1 + 2 \r^2)^2(1 + 4 \r^2)^{p-1}} \over \sqrt{1 + \r^2}} ~.
\label{step2b0} \feq
Now, using Eqs.\,\rref{kneq} and \rref{step2a}\,\rref{step2b0} we infer
Eq.\,\rref{step20}. \fine
\begin{rema}
\label{remtoprove}
\textbf{Remark.} Obviously enough, Propositions \ref{propkmenoa} and \ref{propkmeno} imply
\beq K_{p n} \geqs K^{(-)}_{p n} := \max( K^{[-]}_p, K^{\{ - \}}_{p n}) \label{kpmenoimp} \feq
for all $p \geqs n > d/2$. We anticipate that
the above maximum equals $K^{\{ - \}}_{p n}$ for fixed $d, n$ and $p$ sufficiently large;
this fact is suggested by the subsequent numerical examples, and is proved
in the forthcoming Remark \ref{remadd}.
\end{rema}
\textbf{Some numerical examples.} Let us consider the ($d=3$) cases
\rref{casesdue} \rref{casestre}.
In all these cases, the function $K^{\{-\}}_{p n}(~) : \r \in \naturali_0
\mapsto K^{\{-\}}_{p n}(\r)$ attains its sup
at $\r=1$. In the forthcoming Table \TE we report for each one of the above cases
the values of $K^{[-]}_p$ (defined by \rref{kpmenoa}) and of $K^{\{ - \}}_{p n}$ (the sup of the previous function),
which immediately determine $K^{(-)}_{p n}$ according to \rref{kpmenoimp}
({\footnote{To be precise, the table reports the rounddown to
three digits of all the above mentioned quantities.}}).
The values of $K^{(-)}_{p n}$ have been anticipated in Table \TA of the Introduction.
\par
It should be noted that for larger values of $p$, not considered in these examples,
the $\sup$ of the function $K^{\{-\}}_{p n}(~)$ is attained at a point
$\ell_{p n} > 1$: see the forthcoming subsection. \parn
\spazio
\vbox{
\noindent
\textbf{Table \TE. Lower bounds on $\boma{K_{p n}}$ in the cases \rref{casesdue}
\rref{casestre}}
\vskip 0.2cm \noindent
\vbox{
$$ \begin{tabular}{|c|| c|c||c|}
\hline
$(p,n)$ & $K^{[-]}_p$ & $K^{\{ - \}}_{p n}$ & $K^{(-)}_{p n}$ \\ [0.1cm]
\hline \hline
$(2,2)$ & 0.126 & 0.0809 & 0.126 \\ \hline
$(3,2)$ & 0.179 & 0.147 & 0.179   \\ \hline
$(4,2)$ & 0.253 & 0.264 & 0.264   \\ \hline
$(5,2)$ & 0.359 & 0.463 & 0.463 \\ \hline
$(6,2)$ & 0.507 & 0.793 & 0.793 \\ \hline
$(7,2)$ & 0.718 & 1.33 & 1.33   \\ \hline
$(8,2)$ & 1.01 & 2.20 & 2.20   \\ \hline
$(9,2)$ & 1.43 & 3.60 & 3.60    \\ \hline
$(10,2)$ & 2.03 & 5.83 & 5.83   \\ \hline
$(3,3)$ & 0.179 & 0.125 & 0.179 \\ \hline
$(4,3)$ & 0.253 & 0.232 &  0.253  \\ \hline
$(5,3)$ & 0.359 & 0.418 & 0.418    \\ \hline
$(6,3)$ & 0.507 & 0.732 & 0.732   \\ \hline
$(7,3)$ & 0.718 & 1.25 & 1.25    \\ \hline
$(8,3)$ & 1.01 & 2.10 & 2.10    \\ \hline
$(9,3)$ & 1.43 & 3.48 & 3.48    \\ \hline
$(10,3)$ & 2.03 &  5.69 & 5.69  \\ \hline
\end{tabular} $$
}}
\spazio
\textbf{More on the lower bounds}.
Some numerical experiments performed on the function $K^{\{-\}}_{p n}(~) : \r \in \naturali_0
\mapsto K^{\{-\}}_{p n}(\r)$ for $n \leqs 10$ and
$n \leqs p \leqs 1000$ indicate that, within this range,
$K^{\{-\}}_{p n}(~)$ attains its maximum at a point $\r_{p n}$, which is
approximated with good accuracy by $\rr_{p n} \in \naturali_0$ defined hereafter.
To define this approximant we first introduce the nonnegative real number
\beq \lambda_{p n} := \left\{ \barray{ll}
\sqrt{{p-n \over 2 (n+1)}} & \mbox{for $d \geqs 3$,} \\
\sqrt{{p-n + 2 \over 2 (n+1)}} & \mbox{for $d = 2$} \farray \right.
\label{defla} \feq
and then we put
\beq
\rr_{p n} := \left\{ \barray{ll} \llfloor \lambda_{p n} \rrfloor
& \mbox{if $K^{\{-\}}_{p n}(\llfloor \lambda_{p n}\rrfloor) >
K^{\{-\}}_{p n}(\llceil \lambda_{p n}\rrceil)$} \\
\llceil \lambda_{p n} \rrceil
& \mbox{if $K^{\{-\}}_{p n}(\llceil \lambda_{p n}\rrceil) \geqs
K^{\{-\}}_{p n}(\llfloor \lambda_{p n}\rrfloor)$.}
\farray \right. \label{defer} \feq
In the above, for each real number $x \geqs 0$ we intend
\beq \llfloor x \rrfloor := \max(\lfloor x \rfloor, 1)~, \quad \llceil x \rrceil := \max(\lceil x \rceil, 1)~,
\label{defll} \feq
where $\lfloor x \rfloor$ and $\lceil x \rceil$ are the usual lower and
upper integer parts of $x$ (i.e., the unique integers such that
$\lfloor x \rfloor \leqs x < \lfloor x \rfloor + 1$ and $\lceil x \rceil -1 < x \leqs \lceil x \rceil$;
this implies, for example, $x \leqs \lceil x \rceil < x + 1$ and $x \leqs \llceil x \rrceil \leqs x + 1$).
\par
Let us exemplify the accuracy of the approximation of $\r_{p n}$ via
$\rr_{p n}$ in a number of cases. As a matter of fact, $\r_{p n}$ coincides \textsl{exactly}
with $\rr_{p n}$ for $d=n=3$, $p=3,4,...,69,70$
and for $d=n=2$, $p=2,3,...,89,90$; moreover, $\r_{p n}$ equals $\rr_{p n}$ or $\rr_{p n}+1$ for
$d=n=3$, $p=80,90,...,990,1000$ and for $d=n=2$, $p=100,110,....,790,800$.
\vskip 0.2cm \noindent
\textbf{The large $\boma{p}$ limit of the previous lower bounds.}
We know that $K_{p n} \geqs K^{\{-\}}_{p n}(\r)$ for all $\r \in \naturali_0$;
therefore (independently of the previous considerations and
numerical experiments on the maximum over $\r$ of $K^{\{-\}}_{p n}(\r)$)
we are granted that
\beq K_{p n} \geqs K^{\{-\}}_{p n}(\llceil \lambda_{p n} \rrceil)~, \label{wehaves} \feq
where $\lambda_{p n}$ and $\llceil~\rrceil$ are defined via
Eqs.\,\rref{defla}\,\rref{defll}. This obvious remark yield the
rigorous statement presented hereafter.
\begin{prop}
\label{lowlim}
\textbf{Proposition.}
(i) Let $d \geqs 3$. For all real $p, n$ with $p \geqs n > d/2$ one has
\beq K_{p n} \geqs K^{\la - \ra}_{p n}~, \feq
where
\beq K^{\la - \ra}_{p n} := {1 \over (2 \pi)^{d/2}} \, \big({n+1 \over p}\big)^{(n + 1)/2} \,
{ \Xa_{p n}^{1/2} \cdot 2^{p+n/2 + 1} \over
\Xb_{p n}^{p/2} \Xc_{p n}^{n/2 + 1/2} + \Xb_{p n}^{n/2} \Xc_{p n}^{p/2 + 1/2} }~,
 \label{935} \feq
\parn
\vbox{
\beq \Xa_{p n} := \big(1 - {n-1 \over 2 p}\big)^p + \big( {n+1 \over 2 p} \big)^p ~,
\qquad \label{xa} \feq
$$ \Xb_{p n} := 1  + 2 \sqrt{1 - {n \over p}} \, \, \sqrt{{2 (n+1) \over p}}
+ {n+2 \over p}~, $$
$$ \Xc_{p n} := 1  + 2 \sqrt{1 - {n \over p}} \, \, \sqrt{{2 (n+1) \over p}}
+ {3 n+4 \over p}~. $$
}
(ii) Let $d =2$. For all real $p, n$ with $p \geqs n > 1$ one has
\beq K_{p n} \geqs K^{\la - \ra}_{p n}~, \feq
where
\beq K^{\la - \ra}_{p n} := {1 \over \pi} \, \big({n+1 \over p}\big)^{(n+1)/2} \,
{\Xa_{p n}^{1/2} \cdot 2^{p+n/2} \over
\Xb_{p n}^{p/2} \Xc_{p n}^{n/2 + 1} + \Xb_{p n}^{n/2} \Xc_{p n}^{p/2 + 1} }
~, \label{835} \feq
\parn
\vbox{
\beq \Xa_{p n} := \big(1 + {3 \over p}\big)^2 \big(1 - {n-5 \over 2 p}\big)^{p-1} +
4 \big( {n+1 \over 2 p} \big)^{p+1} ~, \label{xba} \feq
$$ \Xb_{p n} := 1  + 2 \sqrt{1 - {n - 2\over p}} \, \, \sqrt{{2 (n+1) \over p}}
+ {n+4 \over p}~, $$
$$ \Xc_{p n} := 1  + 2 \sqrt{1 - {n - 2\over p}} \, \, \sqrt{{2 (n+1) \over p}}
+ {3 (n+ 2) \over p}~. $$
}
(iii) Let $d \geqs 2$. For fixed $n$ and $p \vain + \infty$, one has
\beq K^{\la - \ra}_{p n} \sim~{2 \over (2 \pi)^{d/2}} \, {e^{5 (n+1)/4} \over e^{n+1} + 1} \,
\big({n+1 \over p}\big)^{(n + 1)/2} \,
{2^{p + n/2} \over e^{\sqrt{2 (n+1) p}}} \label{936} \feq
and
\beq (K^{\la - \ra}_{p n})^{1/p} \vain 2~. \label{937} \feq
\end{prop}
\textbf{Proof.} (i) Let $d \geqs 3$. We have
\parn
\vbox{
\beq K_{p n} \geqs_{(1)}
{\sqrt{2} \over (2 \pi)^{d/2}}\,
{\sqrt{1 + (1 + 4 \llceil \lambda_{p n}\rrceil^2)^p} \over
\llceil \lambda_{p n}\rrceil^p (1+ \llceil \lambda_{p n}\rrceil^2)^{n/2 + 1/2} +
\llceil \lambda_{p n}\rrceil^n (1 + \llceil \lambda_{p n}\rrceil^2)^{p/2 + 1/2}}
\label{same} \feq
$$ \geqs_{(2)} {\sqrt{2} \over (2 \pi)^{d/2}} \,
{\sqrt{1 + (1 + 4 \lambda^2_{p n})^p} \over
(1 + \lambda_{p n})^p (1+ (1 + \lambda_{p n})^2)^{n/2 + 1/2} +
(1 + \lambda_{p n})^n (1+ (1 + \lambda_{p n})^2)^{p/2 + 1/2}} $$
$$ =_{(3)} K^{\la - \ra}_{p n}~\mbox{as in \rref{935}\rref{xa}}. $$
}
In the above: the inequality $\geqs_{(1)}$ is just the relation
\rref{wehaves} with the explicit expression coming from \rref{kpmeno1} for
$K^{\{-\}}_{p n}(\llceil \lambda_{p n}\rrceil)$;
the inequality $\geqs_{(2)}$ follows from $\lambda_{p n} \leqs \llceil \lambda_{p n}\rrceil
\leqs 1 + \lambda_{p n}$; the equality $=_{(3)}$ follows using
the explicit expression \rref{defla} for $\lambda_{p n}$ and
performing some elementary manipulations. \parn
(ii) Let $d=2$. A chain of relations with the same
structure as \rref{same}, justified by the same arguments
employed in the proof of (i), yields the conclusion
$K_{p n} \geqs K^{\la - \ra}_{p n}~\mbox{as in \rref{835}\rref{xba}}$. \parn
(iii) To prove Eq.\,\rref{936} one starts from the explicit expression
of $K^{\la - \ra}_{p n}$, given by Eq.\,\rref{935} for $d \geqs 3$ and by Eq.\,\rref{835} for $d=2$;
both for $d \geqs 3$ and for $d=2$, the $p \vain +\infty$ limits are performed in
an elementary way. The derivation of Eq. \rref{936} for $d \geqs 3$ uses, e.g., the relations
\beq
\Big(1  - {n-1 \over 2 p}\Big)^p =
e^{- (n-1)/2} \big(1 + O({1 \over p})\big)~, \feq
$$
\Big(1  + 2 \sqrt{1 - {n \over p}} \, \, \sqrt{{2 (n+1) \over p}}
+ {n+2 \over p}\Big)^{p/2} = e^{- 3 n/2 - 1} e^{\sqrt{2 (n+1) p}} \big(1 + O({1 \over \sqrt{p}})\big)~;
$$
for $d=2$ one uses, e.g., the relations
\beq
\Big(1  - {n-5 \over 2 p}\Big)^{p-1} =
e^{- (n-5)/2} \big(1 + O({1 \over p})\big)~, \feq
$$
\Big(1  + 2 \sqrt{1 - {n - 2\over p}} \, \, \sqrt{{2 (n+1) \over p}}
+ {3(n+2) \over p}\Big)^{p/2 + 1} = e^{- n/2 + 1} e^{\sqrt{2 (n+1) p}} \big(1 + O({1 \over \sqrt{p}})\big)\,.
$$
Eq.\,\rref{937} is a consequence of \rref{936}. \fine
\begin{rema} \label{remadd}
\textbf{Remark.} Let us compare the $p \vain + \infty$ asymptotics
\rref{936} of $K^{\la - \ra}_{p n}$ with the explicit expression
\rref{kpmenoa} of the alternative lower bound $K^{[-]}_{p}$, which
has the form $2^{p/2}$ $\times$ a constant depending only on $d$.
From here it is evident that $K^{[-]}_p < K^{\la - \ra}_{p n}$ for
fixed $n > d/2$ and all sufficiently large $p$. On the other hand,
from the previous considerations it is evident that
$K^{\la - \ra}_{p n} \leqs K^{\{-\}}_{p n}(\llceil \lambda_{p n} \rrceil)
\leqs K^{\{-\}}_{p n}$ (where, we recall, $K^{\{-\}}_{p n}$ is
the sup of the function $\ell \mapsto K^{\{-\}}_{p n}(\ell)$). \par
In conclusion, for fixed $n > d/2$ and all $p$ sufficiently large
we have $K^{[-]}_p < K^{\{-\}}_{p n}$, which proves the last
statement in Remark \ref{remtoprove}.
\end{rema}
\section{Lower bounds for the sharp constant $\boma{G_{p n}}$
of the inequality \rref{katineqa}}
\label{seclowg}
Let $p,n$ be as in \rref{41}.
Similarly to the case of $K_{p n}$, to obtain a lower bound on $G_{p n}$ we use
the tautological inequality
\beq G_{p n} \geqs
{2 | \la \PPP(v, w) | w \ra_p | \over
(\| v \|_p \| w \|_n + \| v \|_n \| w \|_p) \| w \|_p }
\quad \mbox{for $v \in \HM{p} \setminus \{0 \}$, $w \in \HM{p+1} \setminus \{0 \}$},
\label{gneq} \feq
choosing for $v$ and $w$ two suitable
non zero trial vector fields. Hereafter
we present a simple choice of $v,w$,
depending on a discrete parameter
and on three continuous parameters,
giving rise to useful results. In the
sequel we keep from Section \ref{seclowk} the notations
$$ a := (1,0,...,0),~~ b := (0,1,0,...0);~~ c := (0,0,1,0,...,0)~\,\mbox{if $d \geqs 3$};
~~\naturali_0 := \{1,2,3,...\}. $$
\begin{prop}
\label{progmeno}
\textbf{Proposition.} (i) Let $d \geqs 3$. For all $\X = (\r,\A,\B,\F) \in \naturali_0 \times
\reali^3$, it is
\beq G_{p n} \geqs G^{-}_{p n}(\X)~, \quad
G^{-}_{p n}(\X) :=
{2 \sqrt{2} \over (2 \pi)^{d/2}} {|S_p(\X)| \over (\ell^ p N_n(\X) + \ell^ n N_p(\X)) N_p(\X)}
~, \label{tesgn} \feq
where
\beq S_p(\X) := -\A + (1 + \r^2)^p \A(1-\B) + (1 + 4 \r^2)^p \B (\A-\F) + (1 + 9 \r^2)^p \B \F~, \label{sn} \feq
\beq {~} \hspace{-0.5cm}
N_m(\X) \!:= \!\sqrt{1 + 2 (1 + \r^2)^m \A^2 + 2 (1 + 4 \r^2)^m \B^2
+ 2 (1 + 9 \r^2)^m \F^2}~\mbox{for $\hspace{-0.2cm}$ all $m \in \reali$}. \label{nm} \feq
(ii) Let $d=2$. For all $\X = (\ell,\A,\B,\F) \in \naturali_0 \times \reali^3$, it is
\beq G_{p n} \geqs G^{-}_{p n}(\X)~, \quad
G^{-}_{p n}(\X) :=
{\sqrt{2} \over \pi} {|\SS_p(\X)| \over (\ell^p N_n(\X) + \ell^n N_p(\X)) N_p(\X)}
~, \label{tesgn2} \feq
where $N_n(\X)$, $N_p(\X)$ are defined following Eq.\,\rref{nm} and
\beq \SS_p(\X) := - {\A \r \over \sqrt{1 + \r^2}} +
2 \A \r \Big(1 - {3 \B \r \over \sqrt{4 + \r^2}}\Big) (1 + \r^2)^{p - 3/2} \label{sn2} \feq
$$ + {5 \B \r^2 \over \sqrt{4 + \ell^2}}
\Big({3 \A \over \sqrt{1 + \r^2}} - {7 \F \over \sqrt{9 + \r^2}}\Big) (1 + 4 \r^2)^{p-1} +
{70 \B \F \r^2 \over \sqrt{(4 + \r^2)(9 + \r^2)}} (1 + 9 \r^2)^{p-1}. $$
(iii) For any $d \geqs 2$, Eqs.\,\rref{tesgn}\,\rref{tesgn2} imply
\beq G_{p n} \geqs G^{-}_{p n} \qquad G^{-}_{p n} := \sup_{\X \in \naturali_0 \times \reali^3} G^{-}_{p n}(\X)~.
\label{tesgnb2} \feq
\end{prop}
\textbf{Proof.}
(i) Let $d \geqs 3$. For $\X := (\r,\A,\B,\F) \in \naturali_0 \times \reali^3$ we set
\beq v \equiv v(\X) := i b (e_{\r a} - e_{-\r a}), \label{weset} \feq
$$ w \equiv w(\X) := c\Big(\,  e_{b} + e_{-b} + \A  (e_{\r a+ b} + e_{-\r a- b} - e_{\r a-b} - e_{-\r a + b}) $$
$$+ \B  (e_{2 \r a+b} + e_{-2 \r a-b} +  e_{2 \r a-b} + e_{-2 \r a + b})
+ \F  (e_{3 \r a+b} + e_{-3 \r a-b} -  e_{3 \r a-b} - e_{-3 \r a + b})\, \Big)~; $$
incidentally, we note that $v$ is as in Eq.\,\rref{defvi}.
Like the already considered vector field $v$, the vector field $w$ has vanishing mean and divergence
($k \sc w_k=0$ for each $k$), and belongs to
$\HM{m}$ for any real $m$. For each $m$ we have
\beq \| v \|_m = \sqrt{2} \,\r^m~, \qquad  \| w \|_m = \sqrt{2} \, N_m(\X)~, \label{stepg1a} \feq
with $N_m(\X)$ as in \rref{nm}.
Using Eqs.\,\rref{deflpk}\,\rref{infert} we find
\parn
\vbox{
\beq \PPP(v,w) =  {c \over (2 \pi)^{d/2}} \Big(
2 \A (e_{b} + e_{-b}) + (\B-1) (e_{\r a+b} + e_{-\r a-b} - e_{\r a-b} - e_{-\r a+b}) \feq
$$ + (\F-\A) (e_{2 \r a+b} + e_{-2 \r a-b} + e_{2 \r a-b} + e_{-2 \r a+b})
- \B (e_{3 \r a+b} + e_{-3 \r a-b} - e_{3 \r a-b} - e_{-3 \r a+b}) $$
$$ - \F (e_{4 \r a+b} + e_{-4 \r a-b} + e_{4 \r a-b} + e_{-4 \r a+b}) \Big)~;
$$
}
from here we infer
\beq \la \PPP(v,w) | w \ra_{p} = - {4 \over (2 \pi)^{d/2}} S_p(\X)~, \label{stepg1b} \feq
with $S_p(\X)$ as in \rref{sn}.
Now, using Eqs.\,\rref{stepg1a}\,\rref{stepg1b} and \rref{gneq} we readily infer
statement \rref{tesgn}. \parn
(ii) Let $d=2$. We set
\beq v \equiv v(\X) := i b (e_{\r a} - e_{- \r a}), \label{4562} \feq
$$ w \equiv w(\X) := a (e_{b} + e_{-b}) + {\A \over \sqrt{1 + \r^2}} \Big( (- \r a + b)
\,(e_{\r a+b} + e_{-\r a-b} + (\r a + b) \, (e_{\r a-b} + e_{-\r a + b}) \Big) $$
$$ + {\B \over \sqrt{4 + \r^2}} \Big( (\r a - 2 b) \, (e_{2 \r a+b} + e_{-2 \r a-b})
+ (\r a + 2 b) \, (e_{2 \r a - b} + e_{-2 \r a + b}) \Big) $$
$$ + {\F \over \sqrt{9 + \r^2}} \Big( (- \r a + 3 b) \, (e_{3 \r a+b} + e_{-3 \r a-b})
+ (\r a + 3 b) \, (e_{3 \r a - b} + e_{-3 \r a + b}) \Big)~, $$
where $\X = (\r,\A,\B,\F) \in \naturali_0 \times \reali^3$; again, $v$ has the structure
\rref{defvi}. $v$, $w$ have vanishing mean and divergence and belong to
$\HM{m}$ for any real $m$; their norms of any order $m$ are given again by
\rref{stepg1a} \rref{nm}.
Using Eqs.\,\rref{lpk}\,\rref{infert} we find
\parn
\vbox{
\beq \PPP(v,w) = - {\A \r \over \pi \sqrt{1 + \r^2}} a (e_b + e_{-b}) \feq
$$ - {1 \over 2 \pi (1 + \r^2)} \Big(1 - {3 \B \r \over \sqrt{4 + \r^2}}\Big)
\, \Big( ( a- \r b) (e_{\r a + b} + e_{- \r a - b}) -
(a + \r b) (e_{\r a - b} + e_{-\r a + b}) \Big) $$
$$ + {\r \over 2 \pi(1 + 4 \r^2)}
\Big( {3 \A \over \sqrt{1 + \r^2}} - {7 \F \over \sqrt{9 + \r^2}} \Big)
 \, \Big( ( a- 2 \r b) (e_{2 \r a + b} + e_{- 2 \r a - b})
+ (a + 2 \r b) (e_{2 \r a - b} + e_{-2 \r a + b}) \Big) $$
$$ + {7 \B \r \over 2 \pi (1 + 9 \r^2) \sqrt{4 + \r^2}} \,
\Big( ( - a + 3 \r b) (e_{3 \r a + b} + e_{- 3 \r a - b})
+ (a + 3 \r b) (e_{3 \r a - b} + e_{-3 \r a + b}) \Big) $$
$$ + {13 \F \r \over 2 \pi (1 + 16 \r^2) \sqrt{9 + \r^2}} \,
\Big( ( a - 4 \r b) (e_{4 \r a + b} + e_{- 4 \r a - b})
+ (a + 4 \r b) (e_{4 \r a - b} + e_{-4 \r a + b}) \Big)~; $$
}
from here we infer
\beq \la \PPP(v,w) | w \ra_{p} = {2 \over \pi} \,\SS_p(\X)~, \label{stepg1b2} \feq
with $\SS_p(\X)$ as in \rref{sn2}.
Now, using Eqs.\,\rref{stepg1a}\,\rref{stepg1b2} and \rref{gneq} we readily infer
statement \rref{tesgn2}. \parn
(iii) Obvious. \fine
\textbf{Some numerical examples.} Let us consider
the lower bounds $G^{-}_{p n}(\X)$ of Eq.\,\rref{tesgn},
depending on $\X = (\r,\A,\B,\F) \in \naturali_0 \times \reali^3$. For given $(p,n)$, one
can try to maximize this bound with respect to $\X$
using the routines for numerical optimization of {\tt{Mathematica}};
these predict the maximum to be located at some point
$\X_{p n} = (\ell_{p n}, \A_{p n}, \B_{p n}, \F_{p n})$. Even though this is not the actual point
of absolute maximum, the number
\beq G^{(-)}_{p n} := G^{-}_{p n}(\X_{p n}) \feq
is in any case a lower bound for $G_{p n}$. In the forthcoming
Table \TF we report, for the ($d=3$) cases \rref{casestre},
the point $\X_{p n}$ provided by {\tt{Mathematica}} and
the value of $G^{(-)}_{p n}$
({\footnote{More precisely, in the table
we write $G^{(-)}_{p n}$ for the rounddown
to three digits of $G^{-}_{p n}(\X_{p n})$.}}).
The numerical values $G^{(-)}_{p n}$ have been anticipated in Table \TB
of the Introduction.
\spazio
\vbox{
\noindent
\textbf{Table \TF. Maximizing points $\X_{p n}$ and lower bounds $\boma{G^{(-)}_{p n}}$
in the cases \rref{casestre}}
\vskip 0.2cm \noindent
\vbox{
$$ \begin{tabular}{|c||c|| c|}
\hline
$(p,n)$ & $\X_{p n}$ & $G^{(-)}_{p n}$ \\ [0.1cm]
\hline \hline
$(3,3)$ & (1, 0.388104..., 0.084359..., 0.0135851...) & 0.121 \\ \hline
$(4,3)$ & (1, 0.370907..., 0.0628525..., 0.00811876...) & 0.235  \\ \hline
$(5,3)$ & (1, 0.361597..., 0.0415026..., 0.00365302...) & 0.408    \\ \hline
$(6,3)$ & (1, 0.352601..., 0.0256793..., 0.00147754...) & 0.674   \\ \hline
$(7,3)$ & (1, 0.348117..., 0.0157944..., 0.000588218...) & 1.08    \\ \hline
$(8,3)$ & (1, 0.352603..., 0.00994449..., 0.000238632...) & 1.74    \\ \hline
$(9,3)$ & (1, 0.367597..., 0.00645276..., 0.0000993805...) & 2.77    \\ \hline
$(10,3)$ & (1, 0.392975..., 0.00430116..., 0.0000423667...) & 4.40  \\ \hline
\end{tabular} $$
}}
\vskip 0.2cm \noindent
In the table we always have $\ell_{p n}=1$.
This is no more the case for larger values of $p$: for example,
{\tt{Mathematica}} gives $\ell_{p n} = 2$ for ($d=3$ and) $(p,n) = (20,3)$.
\salto
\textbf{The large $\boma{p}$ limit of the previous lower bounds.}
We know that $G_{p n} \geqs G^{-}_{p n}(\X)$ for all $\X \in \naturali_0 \times \reali^3$;
in this section we choose for $\X$ the quadruple
\beq \Xx_{p n} := (\rr_{p n},1,2^{-p},0)~, \qquad \rr_{p n} := \lceil \sqrt{p/n} \, \rceil
\label{xpn} \feq
where, as in the previous section, $\lceil ~\rceil$ denotes the upper integer part. This
choice gives
\beq G_{p n} \geqs G^{-}_{p n}(\Xx_{p n}) \label{wehavesg} \feq
for all $p \geqs n > d/2+1$. We make the choice \rref{xpn} because some
numerical experiments seem to indicate that, for large $p$, $G^{-}_{p n}(\X)$
attains its maximum for $\X$ close to $\Xx_{p n}$ (both for $d=2$,
and for $d \geqs 3$). Independently of
these experiments, Eq.\,\rref{wehavesg} yields
the rigorous statement presented hereafter.
\begin{prop}
\label{lowlimg}
\textbf{Proposition.}  (i) Let $d \geqs 3$. For all real $p, n$ with $p \geqs n > d/2 +1$ one has
\beq G_{p n} \geqs G^{\la - \ra}_{p n}~, \feq
where
\beq G^{\la - \ra}_{p n} :=
{\sqrt{2} \over (2 \pi)^{d/2}} \Big({n \over p}\Big)^{n/2}
{\dd{\Big(1 + {n \over 4 p}\Big)}^p \cdot 2^p + \dd{\Big(1 - {1 \over 2^p}\Big)}
\dd{\Big(1 + {n \over p}\Big)}^p - \dd{\Big({n \over p}\Big)}^p
\over
\Big(\dd{\Big(1 + {\sqrt{n \over p}}\Big)}^p \Upsilon_{p n}^{1/2} +
\dd{\Big(1 + {\sqrt{n \over p}}\Big)}^n \Sigma_{p n}^{1/2} \Big) \Sigma_{p n}^{1/2}  }
~, \label{935g} \feq
\parn
\vbox{
\beq \Sigma_{p n} := \Big(1  + 2 {\sqrt{n \over p}} + {5 n \over 4 p}\Big)^p +
\Big(1  + 2 {\sqrt{n \over p}} + {2 n \over p}\Big)^p + {1 \over 2} \Big({n \over p}\Big)^p, \label{sigma} \feq
$$ \Upsilon_{p n} := \Big(1  + 2 {\sqrt{n \over p}} + {2 n \over p}\Big)^n +
 {1 \over 2^{2 p - 2 n}} \Big(1 + 2 \sqrt{n \over p}
 + {5 n \over 4 p}\Big)^n + {1\over 2}\Big({n \over p}\Big)^n. $$
}
(ii) Let again $d \geqs 3$. For fixed $n$ and $p \vain + \infty$, one has
\beq G^{\la - \ra}_{p n} \sim~{\sqrt{2} \over (2 \pi)^{d/2}}\,
{e^{9 n/8} \over (1 + e^{n/8} \sqrt{1 + e^{3 n/4}}) \sqrt{1 + e^{3 n/4}} }\,
\Big({n \over p}\Big)^{n/2}
{2^{p} \over e^{2 \sqrt{n p}} } \label{936g} \feq
and
\beq (G^{\la - \ra}_{p n})^{1/p} \vain 2~. \label{937g} \feq
(iii) Let $d =2$. For all real $p, n$ with $p \geqs n > 2$ one has
\beq G_{p n} \geqs G^{\la - \ra}_{p n}~, \feq
where
\parn
\vbox{
\beq G^{\la - \ra}_{p n} :=
{\Big(\dd{n \over p}\Big)^{n/2 + 1} \over \sqrt{2} \pi \sqrt{\dd{\Big(1 + {n \over p}\Big)
\Big(1 + {4 n \over p}\Big)}} }  \label{935g2} \feq
$$ \times \, {15 \dd{\Big(1 + {n \over 4 p}\Big)}^{p-1} \, 2^{p-2} +
\dd{\Big(1 + {n \over p}\Big)}^{p-1} \dd{\Big( 2 \sqrt{1 + {4 n \over p}} } -
 \, {3 \over 2^{p-1}} \Big) - \dd{\sqrt{1 + {4 n \over p}}} \, \Big({n \over p}\Big)^{p-1}
\over
\Big(\dd{\Big(1 + {\sqrt{n \over p}}\Big)}^p \Upsilon_{p n}^{1/2} +
\dd{\Big(1 + {\sqrt{n \over p}}\Big)}^n \Sigma_{p n}^{1/2} \Big)~ \Sigma_{p n}^{1/2}  }
$$
}
and $\Sigma_{p n}, \Upsilon_{p n}$ are as in Eq.\,\rref{sigma}. \parn
(iv) Let again $d =2$. For fixed $n$ and $p \vain + \infty$, one has
\beq G^{\la - \ra}_{p n} \sim~{15 \over 4 \sqrt{2} \pi}\,
{e^{9 n/8} \over (1 + e^{n/8} \sqrt{1 + e^{3 n/4}}) \sqrt{1 + e^{3 n/4}} }\,
\Big({n \over p}\Big)^{n/2 + 1}
{2^{p} \over e^{2 \sqrt{n p}} } \label{936g2} \feq
and
\beq (G^{\la - \ra}_{p n})^{1/p} \vain 2~. \label{937g2} \feq
\end{prop}
\textbf{Proof.} (i) Let $d \geqs 3$; for $\gamma \in [1,+\infty)$, we
put
\beq \Ss_p(\gamma) := -1 + (1 - 2^{-p}) (1 + \gamma^2)^p + 2^{-p} (1 + 4 \gamma^2)^p~, \feq
\beq \Nn_m(\gamma) := \sqrt{1 + 2 (1 + \gamma^2)^m + 2^{1 - 2 p} (1 + 4 \gamma^2)^m}
\qquad \mbox{($m=n,p$)}~; \label{defnm} \feq
$\Ss_p, \Nn_n, \Nn_p$ are positive, strictly increasing functions on $[1,+\infty)$. We have
\parn
\vbox{
\beq G_{p n} \geqs_{(1)}
{2 \sqrt{2} \over (2 \pi)^{d/2}} {\SS_p(\rr_{p n}) \over
(\rr_{p n}^ p \NN_n(\rr_{p n}) + \rr_{p n}^n \NN_p(\rr_{p n})) \NN_p(\rr_{p n})}
\label{sameg} \feq
$$ \geqs_{(2)} {\sqrt{2} \over (2 \pi)^{d/2}} \,
{\SS_p(\sqrt{p/n}) \over
\big(\,(1 + \sqrt{p/n})^p \NN_n(1 + \sqrt{p/n}) + (1 + \sqrt{p/n})^n \NN_p(1 + \sqrt{p/n})\,\big)\NN_p(1 + \sqrt{p/n})} $$
$$ =_{(3)} G^{\la - \ra}_{p n}~\mbox{as in \rref{935g}}. $$
}
In the above: the inequality $\geqs_{(1)}$ is just the relation
\rref{wehavesg} with the explicit expression coming from \rref{tesgn} for
$G^{\{-\}}_{p n}(\X)$ and from \rref{xpn} for $\Xx_{p n}$;
the inequality $\geqs_{(2)}$ is obtained noting that \rref{xpn} implies $\sqrt{p/n} \leqs \rr_{p n}
< 1 + \sqrt{p/n}$; the equality $=_{(3)}$ follows
performing some elementary manipulations. \parn
(ii) Let again $d \geqs 3$. To prove Eq.\rref{936g}  one starts from the explicit expression
\rref{935g} of $G^{\la - \ra}_{p n}$ and performs the $p \vain +\infty$ limit in
an elementary way, noting for example that
\beq
\Big(1  + {\sqrt{n \over p}}\Big)^p =
e^{- n/2} e^{\sqrt{n p}} \big(1 + O({1 \over \sqrt{p}})\big), \label{using} \feq
$$ \Big(1  + 2 {\sqrt{n \over p}} + {5 n \over 4 p}\Big)^p =
e^{- 3 n/4} e^{2 \sqrt{n p}} \big(1 + O({1 \over \sqrt{p}})\big)~.
$$
Eq.\,\rref{937g} is a consequence of \rref{936g}. \parn
(iii) Let $d=2$; for $\gamma \in [1,+\infty)$ we put
\beq \Ss_p(\gamma) := - {\gamma \over \sqrt{1 + \gamma^2}} +
2 \gamma \Big(1 - {3 \cdot 2^{-p} \gamma \over \sqrt{4 + \gamma^2}}\Big) (1 + \gamma^2)^{p - 3/2} \label{psn2} \feq
$$ + {15 \cdot 2^{-p} \gamma^2 \over \sqrt{(4 + \gamma^2)(1 + \gamma^2)}} (1 + 4 \gamma^2)^{p-1}~,
$$
and we define $\Nn_{m}(\gamma)$ ($m=n,p$) as in Eq.\,\rref{defnm};
$\Ss_p, \Nn_n, \Nn_p$ are positive, strictly increasing functions on $[1,+\infty)$.
We have
\parn
\vbox{
\beq G_{p n} \geqs_{(1)}
{\sqrt{2} \over \pi} {\Ss_p(\rr_{p n}) \over
(\rr_{p n}^p \Nn_n(\rr_{p n}) +
\rr_{p n}^n \Nn_p(\rr_{p n})) \Nn_p(\rr_{p n})}
\label{sameg2} \feq
$$ \geqs_{(2)} {\sqrt{2} \over \pi} \,
{\Ss_p(\sqrt{p/n}) \over
\big(\,(1 + \sqrt{p/n})^p \Nn_n(1 + \sqrt{p/n}) + (1 + \sqrt{p/n})^n \Nn_p(1 + \sqrt{p/n})\, \big)
\Nn_p(1 + \sqrt{p/n})} $$
$$ =_{(3)} G^{\la - \ra}_{p n}~\mbox{as in \rref{935g2}}. $$
}
In the above: the inequality $\geqs_{(1)}$ is just the relation
\rref{wehavesg} with the explicit expression coming from \rref{tesgn2} for
$G^{\{-\}}_{p n}(\X)$ and from \rref{xpn} for $\Xx_{p n}$, $\rr_{p n}$;
the inequality $\geqs_{(2)}$ is obtained recalling that \rref{xpn} implies $\sqrt{p/n} \leqs \rr_{p n}
< 1 + \sqrt{p/n}$; the equality $=_{(3)}$ follows
performing some elementary manipulations. \parn
(iv) Let again $d =2$. To prove Eq.\,\rref{936g2} one starts from the explicit expression
\rref{935g2} of $G^{\la - \ra}_{p n}$ and performs the $p \vain +\infty$ limit,
using Eq.\,\rref{using} and similar elementary relations. Eq.\,\rref{937g2} is a consequence of \rref{936g2}. \fine
\section{On the large $\boma{p}$ behavior of $\boma{K_{p n}}$ and $\boma{G_{p n}}$}
\label{secinf}
For all $p \geqs n > d/2$, Corollary \ref{corrough} and Proposition \ref{lowlim}
give for the sharp constants $K_{p n}$ the bounds
\beq K^{\la-\ra}_{p n} \leqs K_{p n} \leqs K^{\la+\ra}_{p n}~, \label{cant} \feq
with explicit expressions provided by Eqs.\,\rref{kpp} and \rref{935}\,\rref{835}.
For fixed $d, n$ and $p \vain + \infty$ the upper and lower bounds have
similar, but not coinciding asymptotics: in fact, from Eqs.\,\rref{kpp} and \rref{936} we
know that $K^{\la +\ra}_{p n} = C^{\la+\ra}_n \cdot 2^p$
while $K^{\la - \ra}_{p n} \sim C^{\la-\ra}_n$ $\cdot 2^p p^{-(n+1)/2} e^{-\sqrt{2 (n+1) p}}$,
for suitable coefficients $C^{\la \pm \ra}_n$ (depending also on $d$). \par
Due to these different behaviors, Eq.\,\rref{cant} does not determine
the precise $p \vain + \infty$ asymptotics of $K_{p n}$; however,
we can obtain from it a weaker result on the large $p$ limit.
In fact, as already mentioned after Eq.\,\rref{kpp} and in Eq.\,\rref{937}, we have
$(K^{\la \pm \ra}_{p n})^{1/p} \vain 2$
in this limit; thus, Eq.\,\rref{cant} yields the following result.
\begin{prop}
\textbf{Proposition.} For any fixed $n > d/2$, one has
\beq (K_{p n})^{1/p} \vain 2 \qquad \mbox{for $p \vain + \infty$}~. \feq
\end{prop}
One can make similar considerations about the sharp constants $G_{p n}$. For all $p \geqs n > d/2+1$,
Corollary \ref{corroughg} and Proposition \ref{lowlimg}
give the bounds
\beq G^{\la-\ra}_{p n} \leqs G_{p n} \leqs G^{\la+\ra}_{p n}~, \label{cantg} \feq
with explicit expressions provided by Eqs.\,\rref{gpp} and \rref{935g}\,\rref{935g2}.
For fixed $d, n$ and $p \vain + \infty$ the upper and lower bounds have
similar, but not coinciding asymptotics (see again Eq.\,\rref{gpp}, and compare it
with Eqs.\,\rref{936g}\,\rref{936g2}); however, as already indicated after Eq.\,\rref{gpp}
and in Eqs.\,\rref{937g}\,\rref{937g2}, one has $(G^{\la \pm \ra}_{p n})^{1/p} \vain 2$
and these facts, combined with \rref{cantg}, yield the following result.
\begin{prop}
\textbf{Proposition.} For any fixed $n > d/2+1$, one has
\beq (G_{p n})^{1/p} \vain 2 \qquad \mbox{for $p \vain + \infty$}~. \feq
\end{prop}
\vskip 1.5cm
\textbf{Acknowledgments.}
This work has been partly supported by INdAM, INFN and by MIUR, PRIN 2010
Research Project  ``Geometric and analytic theory of Hamiltonian systems in finite and infinite dimensions''.
\vfill \eject \noindent
\appendix
\section{Appendix. The norm of the bilinear map $\boma{P_{h \ell}}$}
\label{appehel}
Let $h, \ell \in \reali^d \setminus \{0 \}$. We consider the bilinear map defined by Eq.\,\rref{phel}, i.e.,
$$ P_{h \ell} : h^{\perp} \times \ell^{\perp} \vain ({h + \ell})^{\perp}, \quad
(a, b) \mapsto P_{h \ell}(a, b) := {a \sc \ell \over |\ell|} \, \LP_{h + \ell} \, b~; $$
we recall that ${~}^{\perp}$ indicates the orthogonal complement in $\complessi^d$,
and that $\LP_{h+\ell}$ is the orthogonal projection of $\complessi^d$ onto $(h+\ell)^{\perp}$.
We are interested in the norm $| P_{h \ell} | := \min \{ Q \in [0,+\infty)~|~|P_{h \ell}(a, b) |
\leqs Q |a| |b |~\mbox{for all $a \in h^\perp$, $b \in \ell^\perp$} \}$.
\begin{prop}
\textbf{Lemma.} The above norm is given by Eq.\,\rref{eqnorm}, i.e.,
$$ |P_{h \ell}| = \left\{ \barray{ll} \sin \te_{h \ell}  & \mbox{if $d \geqs 3$}\,, \\
\sin \te_{h \ell} \cos \te_{h + \ell, \ell} & \mbox{if $d = 2$}\,.
\farray \right. $$
\end{prop}
\textbf{Proof.} Our argument is closely related to the slightly weaker one employed
in \cite{cok}; in the sequel we use the abbreviations
\beq \te_{h \ell} \equiv \te~, \qquad \te_{h + \ell, \ell} \equiv \te'. \feq
Let $S$ denote a two-dimensional subspace of $\reali^d$ containing $h$ and
$\ell$ (of, course, this is unique if $h, \ell$ are linearly independent).
We choose in $S$ an orthonormal basis $\eta_1, \eta_2$
so that $h$ be a positive multiple of $\eta_1$ and $\ell \sc \eta_2 \geqs 0$; then
\beq h = |h|\, \eta_1~,
\qquad \ell = |\ell| (\cos \te  \, \eta_1 + \sin \te \, \eta_2)~. \label{espu} \feq
In $S$ we also consider a second orthonormal basis $\eta'_1, \eta'_2$, chosen
so that $h+\ell$ be a nonnegative multiple of $\eta'_1$ and $\ell \sc \eta'_2 \geqs 0$; then
\beq h + \ell = |h+\ell|  \, \eta'_1~, \qquad \ell = |\ell| (\cos \te'  \, \eta'_1 + \sin \te' \, \eta'_2)~.
\label{espuu} \feq
Finally let $\eta_3,\eta_4, ... ,\eta_d$ be $d-2$ vectors of $\reali^d$,
forming an orthonormal basis for the orthogonal complement of $S$ in $\reali^d$
(obviously enough, this family is empty if $d=2$).
The orthogonal complements of $h$ and $\ell$ in $\complessi^d$ have
the following representations, involving three orthonormal bases:
\beq h^{\perp} = \langle \eta_2, ..., \eta_d \rangle~, \feq
\beq \ell^{\perp} = \langle - \sin \te \, \eta_1 + \cos \te \, \eta_2, \eta_3, ... , \eta_d \rangle =
\langle - \sin \te' \, \eta'_1 + \cos \te' \, \eta'_2, \eta_3, ... , \eta_d \rangle~. \feq
To go on we consider the cases $d \geqs 3$ and $d =2$, separately.
\vskip 0.2cm\noindent
\textsl{Case $d \geqs 3$.} Let us consider any two vectors $a \in h^{\perp}, b \in \ell^{\perp}$;
we can write
\beq a = a_{(2)} \eta_2 + ... + a_{(d)} \eta_d~, \feq
with $a_{(i)} \in \complessi$ for all $i$. From here and from Eq.\,\rref{espu} for $\ell$ we infer
\beq {a \sc \ell \over |\ell|} = a_{(2)} \sin \te~, \feq
whence
\beq {|a \sc \ell| \over |\ell|} = |a_{(2)}| \sin \te \leqs |a | \sin \te~; \feq
moreover, by a general property of orthogonal projections,
\beq | \LP_{h + \ell} b| \leqs |b| ~.\feq
The last two inequalities imply
\beq |P_{h \ell}(a, b)| \leqs \sin \te \,|a| \,|b|~; \label{imp1} \feq
moreover \rref{imp1} holds as an equality for suitable, nonzero
choices of $a, b$. In fact, setting
\beq a_{*} := \eta_2~, \qquad b_{*} := \eta_3 \feq
we have
\beq {a_{*} \sc \ell \over |\ell|} = \sin \te~, \qquad \LP_{h + \ell} \, b_{*} = b_{*} \feq
(because $h+\ell$ is in the subspace $S$, spanned by $\eta_1, \eta_2$,
and $b^{*}$ is orthogonal to it); this implies
\beq  P_{h \ell}(a_{*}, b_{*}) = \sin \te \, b_{*}~, \feq
so that
\beq  |P_{h \ell}(a_{*}, b_{*})| = \sin \te = \sin \te |a_{*}| |b_{*}|~. \label{imp2} \feq
From Eqs.\,\rref{imp1}\,\rref{imp2} we infer, as desired,
\beq |P_{h \ell}| = \sin \te~. \feq
\textsl{Case $d=2$}.  Let us consider any two vectors $a \in h^{\perp}, b \in \ell^{\perp}$; these
can be written as
\beq a = |a| e^{i \phi} \eta_2~, \qquad b = |b| e^{i \psi}( - \sin \te' \, \eta'_1 + \cos \te' \, \eta'_2)~
\quad (\phi, \psi \in \reali)~. \feq
From here and from the representations \rref{espu} for $\ell$, \rref{espuu} for $h + \ell$ we infer
\beq {a \sc \ell \over |\ell|} = |a| e^{i \phi} \sin \te~, \quad
\LP_{h + \ell} \,b = |b| e^{i \psi} \cos \te' \eta'_2~, \feq
so that
\beq P_{h \ell}(a,b) = \sin \te \cos \te' e^{i (\phi + \psi)} |a| \,|b| \, \eta'_{2}~, \qquad
|P_{h \ell}(a,b)| = \sin \te \cos \te'|a| \,|b| ~; \feq
this trivially yields the desired conclusion
\beq |P_{h \ell}| = \sin \te \cos \te'~. \feq


\vskip 0cm \noindent
\begin{thebibliography}{99}
\bibitem{BKM} J. T. Beale, T. Kato, A. Majda, \textsl{Remarks on the breakdown of smooth solutions
for the 3D Euler equations}, Commun. Math. Phys. \textbf{94} (1984), 61-66.
\vsm
\bibitem{Che} S.I. Chernyshenko, P. Constantin, J.C. Robinson, E.S. Titi,
\textsl{A posteriori regularity of the three-dimensional Navier-Stokes
equations from numerical computations}, J. Math. Phys. \textbf{48} (2007), 065204/10.
\vsm
\bibitem{CoFo} P. Constantin, C. Foias, \textsl{Navier Stokes equations}, Chicago University Press (1988).
\vsm
\bibitem{Ham}
R.S. Hamilton,
\textsl{The inverse function theorem of Nash and Moser},
Bull. Amer. Math. Soc. (N.S.) \textbf{7} (1982), 65-222.
\vsm
\bibitem{Kato} T. Kato, \textsl{Nonstationary flows of viscous and ideal fluids in $\reali^3$},
J.Funct.Anal. \textbf{9} (1972), 296-305.
\vsm
\bibitem{bnw}
C. Morosi, M. Pernici, L. Pizzocchero, \textsl{On power series
solutions for the Euler equation, and the
Behr-Ne$\check{\mbox{c}}$as-Wu initial datum}, ESAIM Math. Model.
Numer. Anal. \textbf{47} (2013), 663-688.
\vsm
\bibitem{padova} C. Morosi, M. Pernici, L. Pizzocchero, \textsl{A posteriori estimates
for Euler and Navier-Stokes equations}, in: F. Ancona, A. Bressan, P. Marcati,
A. Marson (Eds.), Hyperbolic Problems:
Theory, Numerics and Applications, Proceedings
of the XIV International Conference held in Padova  (June 25-29, 2012), in:
AIMS Series on Applied Mathematics \textbf{8} (2014), 847-855.
\vsm
\bibitem{reylarge} C. Morosi, M. Pernici, L. Pizzocchero, \textsl{Large order
Reynolds expansions for the Navier-Stokes equations}, Appl. Math. Lett. \textbf{49} (2015) 58-66.
For an extended version, see arXiv:1402.0487\,.
\vsm
\bibitem{accau} C. Morosi, L. Pizzocchero,
\textsl{An $H^1$ setting for the Navier-Stokes equations: Quantitative estimates},
Nonlinear Anal. \textbf{74} (2011), 2398-2414.
\vsm
\bibitem{appeul} C. Morosi, L. Pizzocchero, \textsl{On approximate solutions for
the Euler and Navier-Stokes equations}, Nonlinear Anal. \textbf{75} (2012), 2209-2235.
\vsm
\bibitem{cog} C. Morosi, L. Pizzocchero, \textsl{On the constants
in a Kato inequality for the Euler and Navier-Stokes equations}, Commun. Pure
Appl. Analysis \textbf{11}(2012), 557-586.
\bibitem{cok} C. Morosi, L. Pizzocchero, \textsl{On the constants
in a basic inequality for the Euler and Navier-Stokes equations}, Appl. Math.
Lett. \textbf{26} (2013), 277-284. For an extended version, see arXiv:1007.4412\,.
\vsm
\bibitem{apprey}
C. Morosi, L. Pizzocchero,
\textsl{On the Reynolds number expansion for the Navier-Stokes equations},
Nonlinear Analysis \textbf{95} (2014), 156-174.
\vsm
\bibitem{smo}
C. Morosi, L. Pizzocchero,
\textsl{Smooth solutions of the Euler and Navier-Stokes
equations from the a posteriori analysis of approximate solutions}, Nonlinear Analysis \textbf{113} (2015), 298-308.
\vsm
\bibitem{RSS}
J. C. Robinson, W. Sadowski, R. P. Silva,
\textsl{Lower bounds on blow up solutions of the three-dimensional Navier–Stokes equations
in homogeneous Sobolev spaces}, J. Math. Phys. \textbf{53} (2012), 115618, 15pp.
\vsm
\bibitem{Tem} R. Temam, \textsl{Local existence of $C^\infty$ solutions of the Euler equation
of incompressible perfect fluids},
in ``Turbulence and Navier Stokes equation'',
Proceedings of the Orsay Conference,
Lecture Notes in Mathematics \textbf{565} (1976), 184-193.
\vsm
\bibitem{arb} F. Yohannson, the {\tt{Arb}} library, see {\tt{http://fredrikj.net/arb/}}.
\end{thebibliography}
\end{document}